\newif\ifspringer\springerfalse

\ifspringer
  \documentclass[default]{sn-jnl}

  \usepackage{graphicx}%
  \usepackage{multirow}%
  \usepackage{amsmath,amssymb,amsfonts}%
  \usepackage{amsthm}%
  \usepackage{mathrsfs}%
  \usepackage[title]{appendix}%
  \usepackage{textcomp}%
  \usepackage{manyfoot}%
  \usepackage{booktabs}%
  \usepackage{algorithm}%
  \usepackage{algorithmicx}%
  \usepackage{algpseudocode}%
  \usepackage{listings}%

  \usepackage{tabularray}
  \usepackage{tabularx}
  \UseTblrLibrary{diagbox}
  \NewColumnType{M}[1]{Q[m,c,#1]}  
  \usepackage{makecell}
  \usepackage{makebox}
  \usepackage{calrsfs} 
  \usepackage{jabbrv}
  \usepackage[numbers,sort&compress]{natbib}
\else
  \documentclass[USenglish,10pt,a4paper]{article}
\fi

\PassOptionsToPackage{usenames,dvipsnames,table}{xcolor}

\PassOptionsToPackage{noend}{algpseudocode}
  \PassOptionsToPackage{capitalize,nameinlink}{cleveref}
  \PassOptionsToPackage{breakspaces}{clevethm}
  \PassOptionsToPackage{usenames,dvipsnames,table}{xcolor}

  \usepackage[cal=boondoxo, scr=euler]{mathalfa}
  \ifspringer
    \usepackage[
      alglinenumber,
      eqreset=section,
      thmreset=section,
    ]{myPreamble}
  \else
    \usepackage[bottom=1.18in,left=1.1in,right=1.1in]{geometry}
    \usepackage[
      alglinenumber,
      eqreset=section,
      thmreset=section,
    ]{myPreamble}
  \fi

\makeatletter
  \foreach \x in {gray, blue, orange, red, black} {%
    \edef\temp{%
      \noexpand\expandafter\noexpand\gdef\noexpand\csname\x\noexpand\endcsname{\noexpand\color{\x}}
    }%
    \temp
  }

  \let\OLDand\and
  \def\and{\texorpdfstring{\OLDand}{, }}%
  \ifspringer
    \newcommand{\Sep}{,}
  \else
    \newenvironment{keywords}{\par\noindent{\bf Keywords. }}{}
    \newenvironment{AMS}{\par\noindent{\bf AMS subject classifications. }}{}
    \newcommand{\Sep}{ \(\cdot\)\ }
  \fi

    \newenvironment{assumption}{\begin{ass}}{\end{ass}}
    \newenvironment{corollary}{\begin{cor}}{\end{cor}}
    \newenvironment{definition}{\begin{defin}}{\end{defin}}
    \newenvironment{example}{\begin{es}}{\end{es}}
    \newenvironment{lemma}{\begin{lem}}{\end{lem}}
    \newenvironment{proposition}{\begin{prop}}{\end{prop}}
    \newenvironment{remark}{\begin{rem}}{\end{rem}}
    \newenvironment{theorem}{\begin{thm}}{\end{thm}}
    \newenvironment{stepsize rule}{\begin{stepsize}}{\end{stepsize}}

  \allowdisplaybreaks

  \makeatletter
    \renewcommand{\clevethm@proofsectiontitle}{Omitted proofs of }
  \makeatother
  \Crefname{assumption}{Assumption}{Assumptions}
  \Crefname{example}{Example}{Examples}
  \Crefname{fact}{Fact}{Facts}
  \Crefname{remark}{Remark}{Remarks}
  \crefname{ALG@line}{step}{steps}
  \crefname{enumeratpropi}{property}{properties}
  \crefname{enumeratpropii}{property}{properties}
  \counterwithout{algorithm}{section}
  \counterwithout{table}{section}
  \renewcommand{\thealgorithm}{\arabic{algorithm}}
  
  \crefformat{equation}{(#2#1#3)}

  \crefname{table}{table}{tables}
  \Crefname{table}{Table}{Tables}

\usepackage{caption}
\usepackage{hyphenat}
\usepackage{cases} 
\usepackage{diagbox}
\usepackage{tabularray}
\UseTblrLibrary{diagbox}
\NewColumnType{M}[1]{Q[m,c,#1]}  
\usepackage{makecell}
\usepackage{makebox}
\usepackage{calrsfs} 
\DeclareMathAlphabet{\pazocal}{OMS}{zplm}{m}{n}

\makeatletter
\newcommand{\qindef}{\@ifstar{\@qindefb}{\@qindefa}}
\newcommand{\@qindefb}[2]{\inner{#1,#1}_{#2}}
\newcommand{\@qindefa}[2]{\operatorname{q}_{#2}(#1)}
\makeatother

\def\A{A}
\def\B{B}
\def\mon{\mu}
\def\com{\rho}
\def\monA{\mon_\A}
\def\monB{\mon_\B}
\def\comA{\com_\A}
\def\comB{\com_\B}
\def\etamin{\eta_\mathrm{min}}

\def\M{P}                           
\def\projQ{\proj_{\range{\M}}}      
\def\projR{Q}                       

\def\DRSrho{\rho}
\def\DRSRho{V}
\def\DRSrhocom{\beta_{\rm P}}
\def\DRSrhomon{\beta_{\rm D}}

\newcommand\optional[1]{[#1]}

\newcommand\tp[1]{#1^\top}
\newcommand\inv[1]{#1^{-1}}

\newcommand\tinprod[2]{\langle {#1},\, {#2} \rangle}
\newcommand\norm[1]{\left\| {#1} \vphantom{X} \right\|}
\newcommand\tnorm[1]{\| {#1} \|}

\newcommand\tnormsq[1]{\tnorm{#1}^2}
\newcommand\defset[2]{\left\{ {#1} \;\middle|\; {#2} \right\}}

\def\Tpd{T_{\!\,\mathrm{PD}}}
\def\Tp{T_{\!\,\mathrm{P}}}
\def\Td{T_{\!\,\mathrm{D}}}

\newcommand\other[1]{{#1}^\prime} 

\newcommand{\x}{\bar x}
\newcommand{\y}{\bar y}

\newcommand{\D}{\mathcal{D}}

\renewcommand\thealgorithm{\arabic{algorithm}}

\makeatletter
  \renewcommand{\appendixproof@toc}{subsubsection}
\makeatother

\newcommand{\gph}{\graph}

\DeclarePairedDelimiter{\nrm}{\lVert}{\rVert}
\DeclarePairedDelimiter{\inner}{\langle}{\rangle}

\newcommand{\overbar}[1]{\mkern 4mu\overline{\mkern-4mu#1\mkern-0mu}\mkern 0mu}
\newcommand{\overhat}[1]{\mkern 4mu\widehat{\mkern-4mu#1\mkern-0mu}\mkern 0mu}

\let\oldBox\Box
\renewcommand{\Box}{\mathbin{\oldBox}}
\newcommand\sym[1]{\mathop{\mathbb{S}}^{#1}}

\pgfplotsset{
  compat=1.17,
  colormap={paper}{rgb255(0cm)=(230,230,230); rgb255(1cm)=(200,200,200); rgb255(2cm)=(185,185,185)}
  }
\usepgfplotslibrary{fillbetween, colormaps}
\usetikzlibrary{arrows.meta, calc, shapes.geometric}
\usetikzlibrary{decorations.pathreplacing,decorations.markings,arrows,shapes.misc}
\pgfplotscolorbarCMYKworkaroundfalse
\showgrayfalse

\definecolor{PaperBlue}{HTML}{185477}
\definecolor{PaperGreen}{HTML}{228B22}
\definecolor{PaperOrange}{HTML}{E69C24}
\definecolor{PaperRed}{HTML}{902A3C}

\newif\iftables\tablesfalse

\newcommand{\PRES}{(\M + T)^{-1}\circ\M}

\newif\iftables\tablesfalse
\newif\ifshort\shorttrue

\newcommand{\TheShortTitle}{Convergence of the Preconditioned Proximal Point Method and Douglas--Rachford Splitting in the Absence of Monotonicity}
\newcommand{\TheTitle}{Convergence of the Preconditioned Proximal Point Method and Douglas--Rachford Splitting in the Absence of Monotonicity}

\newcommand{\TheFunding}{%
  This work was supported by:
  the Research Foundation Flanders (FWO) PhD grants No. 1183822N, No. 11M9523N;
  postdoctoral grant No. 12Y7622N and research projects G081222N, G033822N, G0A0920N;
  European Union's Horizon 2020 research and innovation programme under the Marie Skłodowska-Curie grant agreement No. 953348.%
}
\newcommand{\TheKeywords}{%
Convex/nonconvex optimization\Sep
monotone/nonmonotone variational inequalities\Sep
inclusion problems\Sep
preconditioned proximal point algorithm\Sep
Douglas--Rachford splitting\Sep
semimonotone operators
}
\newcommand{\TheAMSsubj}{%
  47H04\Sep 
  49J52\Sep 
  49J53\Sep
  65K05\Sep 
  65K15\Sep 
  90C26. 
}
\newcommand{\TheAbstract}{%
  The proximal point algorithm (PPA) is the most widely recognized method for solving inclusion problems and 
  serves as the foundation for many numerical algorithms. 
  Despite this popularity, its convergence results have been largely limited to the monotone setting. 
  In this work, we study the convergence of (relaxed) preconditioned PPA for a class of nonmonotone problems that satisfy an oblique weak Minty condition.
  Additionally, we study the (relaxed) Douglas-Rachford splitting (DRS) method in the nonmonotone setting by establishing a connection between DRS and the preconditioned PPA with a positive semidefinite preconditioner.  
  To better characterize the class of problems covered by our analysis, we introduce the class of semimonotone operators, offering a natural extension to (hypo)monotone and co(hypo)monotone operators, and describe some of their properties. 
  Sufficient conditions for global convergence of DRS involving the sum of two semimonotone operators are provided. 
  Notably, it is shown that DRS converges even when the sum of the involved operators (or of their inverses) is nonmonotone. Various example problems are provided,
  demonstrating the tightness of our convergence results and highlighting the wide range of applications our theory is able to cover.
}

\ifspringer
\else
  \author{%
    Brecht Evens
    \and
    Pieter Pas
    \and
    Puya Latafat
    \and
    Panagiotis Patrinos
    \thanks{
      \TheAddressKU\newline
      \textit{E-mail:}
      {\tt
        \{%
          \href{mailto:brecht.evens@kuleuven.be}{brecht.evens},%
          \href{mailto:pieter.pas@kuleuven.be}{pieter.pas},%
          \href{mailto:puya.latafat@kuleuven.be}{puya.latafat},%
          \href{mailto:panos.patrinos@kuleuven.be}{panos.patrinos}%
          \}\texttt{@kuleuven.be}.
      }\newline
      \indent
      \TheFunding
    }
  }
\fi

\ifspringer
  \let\OLDthebibliography\thebibliography
  \renewcommand\thebibliography[1]{
    \OLDthebibliography{#1}
    \setlength{\parskip}{0pt}
    \setlength{\itemsep}{0pt plus 0.2em}
  }
\fi

\begin{document}
  \ifspringer
    \title{
      \TheTitle%
      }

    \author*{\fnm{Brecht} \sur{Evens}}\email{brecht.evens@kuleuven.be}
    \author{\fnm{Pieter} \sur{Pas}}\email{pieter.pas@kuleuven.be}
    \author{\fnm{Puya} \sur{Latafat}}\email{puya.latafat@kuleuven.be}
    \author{\fnm{Panagiotis} \sur{Patrinos}}\email{panos.patrinos@kuleuven.be}

    \affil{\orgdiv{STADIUS}, \orgname{KU Leuven}, \orgaddress{\street{Kasteelpark Arenberg 10}, \city{Leuven}, \postcode{3000}, \country{Belgium}}}

    \abstract{\TheAbstract}
    \pacs[MSC Classification]{\TheAMSsubj}

    \maketitle
    
    \setlength\abovedisplayskip{5pt plus 2pt minus 2pt}
    \setlength\belowdisplayskip{5pt plus 2pt minus 2pt}  
  \else
    \renewcommand\footnotemark{}

    \title{\texorpdfstring{\TheTitle}{\TheShortTitle}}
    \date{}

    \maketitle

    \begin{abstract}
      \TheAbstract
    \end{abstract}

  \fi

  \section{Introduction}

The preconditioned proximal point algorithm (PPPA) is one of the most fundamental methods for solving inclusion problems arising in optimization and variational analysis. It aims to find a solution to the inclusion problem 
\begin{equation} \label{prob:P1} \tag{G-I}
    \text{find} \quad z \in \R^n \quad \text{such that} \quad 
    0 \in T z,
\end{equation}
where  \(T:\R^n \rightrightarrows \R^n\) is a set-valued operator which is usually assumed to be monotone.
The iterates of the (relaxed) preconditioned proximal point algorithm (PPPA) can be expressed as follows:
\begin{equation}\label{eq:PPPA-intro}\tag{PPPA}
    \begin{cases}
        \bar{z}^k&\in(\M+T)^{-1}\M z^k\\
        z^{k+1}&=z^k+\lambda_k(\bar{z}^k-z^k).
    \end{cases}
\end{equation}
Here, \(\M \in \R^{n\times n}\) is a symmetric positive semidefinite matrix and \(\seq{\lambda_k}\) is a sequence of  
strictly positive relaxation parameters. 
The classic proximal point algorithm (PPA)
\cite{martinet1970regularisation} 
corresponds to the setting where
the preconditioner $\M$ takes the form of $\nicefrac1\gamma \I$, where $\gamma$ is a positive stepsize parameter,
    and the relaxation parameter $\lambda_k$ is equal to one.

Several
methods are known to be instances of \ref{eq:PPPA-intro}, each with a specific choice of preconditioner aimed at exploiting structure present in the problem. These include well-known techniques such as Douglas-Rachford splitting (DRS)
\cite[Thm. 6]{eckstein1992DouglasRachford},
the alternating direction method of multipliers (ADMM) 
\cite[\S 8]{eckstein1988LionsMercier},
the (proximal) augmented Lagrangian method (ALM) 
\cite[\S 4 \& 5]{rockafellar1976augmented},
the primal-dual hybrid gradient (PDHG) method
\cite[Lem. 2.2]{he2012Convergence}, \cite{esser2010general}
(also known as Chambolle-Pock (CP) \cite{chambolle2011firstorder}),
iterative refinement \cite[\S 4.1.2]{parikh2014proximal},
the method of partial inverses \cite[\S 5]{eckstein1992DouglasRachford},
progressive hedging 
\cite[Thm. 1]{rockafellar2019solving},
and progressive decoupling
\cite[Thm. 4]{rockafellar2019progressive}.

The convergence analysis of the methods mentioned above typically relies upon an underlying monotonicity/convexity assumption \cite{rockafellar1976augmented,rockafellar1976monotone,rockafellar2019progressive}. See also \cite{spingarn1982Submonotone,luque1987nonlinear,eckstein1993Nonlinear,burke1999Variable,farouq2001Pseudomonotone} for related works and generalizations. Departing from this classical setting, we aim to establish convergence of relaxed \ref{eq:PPPA-intro} for a class of nonmonotone operators that are defined through an \textit{oblique weak Minty condition} (see \cref{def:SWMVI}).
This definition involves a symmetric matrix $\DRSRho$ that is not required to be positive or negative definite and controls the level of (non)monotonicity.
When $V$ is equal to the zero matrix, this condition reduces to the Minty variational inequality (MVI) \cite{minty1962, giannessi1998Minty}, sometimes also referred to as variational coherence%
, satisfied for instance by all pseudoconvex and star-convex functions \cite{zhou2017stochastic}.
When \(V\) takes the form of \(\rho \I\) for negative \(\rho\), one obtains the so-called weak MVI 
that has been employed in the context of the extragradient method and forward-backward-forward splitting 
for solving nonconvex-nonconcave min-max problems 
    \cite{diakonikolas2021Efficient,pethick2023solving,pethick2021escaping,bohm2022solving,alacaoglu2023beyond,alacaoglu2024extending}.
Weak MVI
is closely related to the notion of cohypomonotonicity \cite{bauschke2021Generalized,lee2021fast,tran2023extragradient}, for which a convergence analysis of classic PPA was already performed in \cite{pennanen2002local, iusem2003Inexact,combettes2004proximal}. The key difference between the weak MVI and cohypomonotonicity is that the governing inequality is not considered between any two points in the graph of $T$, but instead only between any point in the graph of $T$ and (a subset of) the zeros of $T$ (see also \Cref{ex:pppa:toy}).

In contrast to the traditional convergence analysis relying on firm nonexpansiveness of the resolvent mapping, 
our convergence analysis leverages a projective interpretation that was introduced in the monotone case in \cite{solodov1996Modified,solodov1999hybrid,konnov1997Class}. We show that the same idea can be applied in the nonmonotone setting for the preconditioned proximal point algorithm provided that a different halfspace is considered. More specifically, each iteration of 
\ref{eq:PPPA-intro}
is interpreted as a projection onto a certain halfspace, separating the current iterate from a subset of solutions that satisfy the oblique weak Minty assumption.  
A second key aspect of our convergence analysis is to allow for positive semidefinite preconditioners (as opposed to positive definite ones), which is crucial for establishing a connection with splitting techniques such as (relaxed) Douglas-Rachford splitting (see \cref{sec:equivDRS}). We note that in the monotone setting convergence of \ref{eq:PPPA-intro} with semidefinite preconditioners has been studied in \cite[Thm. 3.4]{latafat2017Asymmetric}, \cite[\S 2.1]{bredies2022degenerate}. In the nonmonotone setting, the aforementioned halfspace will be constructed using so-called \textit{shadow sequences}, obtained by projecting the sequences $(z^k)_{k\in\N}$ and $(\bar{z}^k)_{k\in\N}$ onto the range of the preconditioner \(\M\). 

In the second part of this paper, we consider the class of Douglas/Peaceman-Rachford splitting methods and shift our attention towards structured inclusion problems of the form
\begin{equation}\label{prob:composite}\tag{P-I}
    \text{find} \quad x \in \R^n \quad \text{such that} \quad 
    0 \in \Tp x \coloneqq Ax + B x, 
\end{equation}
where $A : \R^n \rightrightarrows \R^n$ and $B : \R^n \rightrightarrows \R^n$ are two (possibly nonmonotone) operators.  
Given stepsize $\gamma > 0$, strictly positive relaxation parameters $\lambda_k > 0$, and a certain initial guess $s^0\in\R^n$, the (relaxed) Douglas\hyp{}Rachford splitting method consists of the following iterates:
\begin{equation}
\begin{cases}\tag{DRS}\label{eq:DRS}
    \displaystyle   
    u^k&\in J_{\gamma \A}(s^k)\\
    v^k&\in J_{\gamma \B}(2u^k-s^k)\\
    s^{k+1}&=s^k+\lambda_k(v^k-u^k).
\end{cases}
\end{equation}
For $\lambda_k=1$, this method corresponds to classic Douglas-Rachford splitting \cite{douglas1956Numerical}, while for $\lambda_k=2$, it
reduces to
Peaceman\hyp{}Rachford splitting (PRS) \cite{peaceman1955Numerical}. Initially introduced to solve systems of linear equations emerging in heat conduction problems, these methods were later extended in \cite{lions1979Splitting} 
for finding zeros of the sum of two maximally monotone operators, and convergence of the sequence $\seq{u^k}$ to a solution was shown \cite{lions1979Splitting,svaiter2011Weak}.

In the optimization setting, 
global convergence of 
classic DRS (corresponding to $\lambda_k = 1$) in the nonconvex setting was obtained in \cite{li2016Douglas} assuming that one function is Lipschitz differentiable.
Similarly, 
PRS (corresponding to $\lambda_k = 2$) was considered in \cite{li2017Peaceman} under the additional requirement that the smooth function is strongly convex.
These results were then unified in \cite{themelis2020DouglasRachford} and tight stepsize ranges for (relaxed) \ref{eq:DRS} in the nonconvex setting were provided.
When restricting to their setting (one function being Lipschitz differentiable), we will show that when the nonsmooth term has a so-called \emph{semimonotone} subdifferential, our theory allows for a larger stepsize range compared to \cite{themelis2020DouglasRachford} (see \cref{rem:comparison:DRS:optimization}).

A convergence analysis of \ref{eq:DRS} for minimizing the sum of an $\alpha$-convex and a $\beta$-convex function has been provided in \cite{guo2017Convergence} under the assumption that $\alpha + \beta > 0$, i.e. in the convex case.
Similar results were obtained for \ref{eq:DRS} applied to inclusion problems
involving the sum of an $\alpha$\hyp{} and a $\beta$\hyp{}monotone operator \cite{dao2019Adaptive,giselsson2021compositions}
and the sum of 
an $\alpha$\hyp{} and a $\beta$\hyp{}comonotone operator \cite{bartz2022Conical}, both under the assumption that $\alpha + \beta > 0$. 
In the former work, it is obvious that the so-called primal inclusion \eqref{prob:composite} is monotone, while in the latter, the so-called dual inclusion \eqref{eq:dual} is monotone due to \cite[Lem. 2.6]{bauschke2021Generalized}.
Therefore, the achieved results in this research direction have been limited to the case where either the primal or the dual problem is monotone.

    \subsection{Contributions}

Our main contributions can be summarized as follows:
\begin{enumerate}[label={\arabic{enumi}.}]
    \item We 
        establish convergence
    of \ref{eq:PPPA-intro} for a class of nonmonotone operators that 
        have
    oblique weak Minty 
        solutions
    at one (or more) of its zeros (see \cref{def:SWMVI,thm:pppa}). 
    %
        Our convergence result hinges upon mere positive \emph{semidefiniteness} of the preconditioning matrix (see \Cref{tab:PPPA:connections}),
        and tightness of our stepsize range is demonstrated in \cref{ex:pppa:toy}.
    %
    Additionally, sublinear last-iterate convergence rates are obtained in \Cref{thm:pppa:lastiter} under a more stringent oblique comonotonicity assumption.
    \item 
    Leveraging
    a primal-dual connection between \ref{eq:DRS} and \ref{eq:PPPA-intro},
    convergence results of \ref{eq:DRS} are 
    deduced 
        whenever
    the associated primal-dual operator
    has oblique weak Minty solutions, as well as in the more restrictive, comonotone setting
    (see \cref{sec:convergence DRS}).
    Additionally, we demonstrate the tightness of this result in \cref{it:example5:2}.
    \item Despite the underlying nonmonotonicity, (local) linear convergence of \ref{eq:PPPA-intro} is established in \cref{thm:PPPA:linear} under an additional metric subregularity assumption. In turn, this result leads to sufficient conditions for (local) linear convergence of \ref{eq:DRS} for (possibly) nonmonotone piecewise polyhedral mappings (see \cref{cor:DRS:CP:lin}).
    \item We introduce and develop calculus rules for the class of $(\mon, \com)$\hyp{}semimonotone operators,
    which can be viewed as a natural extension of (hypo)monotone and co(hypo)monotone operators (see \cref{rem:relationship:WMVI}). In addition,
        various examples of well-known (nonconvex) function classes with semimonotone subdifferentials
    are provided%
. 
    \item We provide sufficient conditions for the convergence of \ref{eq:DRS}, based on the semimonotonicity of the operators $A$ and $B$ (see \cref{cor:DRS:semi}). 
        T%
    he stepsize region obtained through our analysis encompasses and extends existing results on relaxed Douglas/Peaceman\hyp{}Rachford splitting methods such as \cite{dao2019Adaptive,bartz2022Conical}, see \cref{rem:comparison:DRS}.
    Up to the knowledge of the authors, this is the first work that is able to cover the convergence of \ref{eq:DRS} in the case where neither the primal, nor the dual, nor the primal-dual inclusion is monotone.
    We provide several example problems that demonstrate the applicability of our theory
    beyond the standard monotone setting.
\end{enumerate}

    \subsection{Organization}

The paper is structured as follows.
\Cref{subsec:notation} introduces some notation and recalls some standard definitions.
In \Cref{sec:PPA} the concept of oblique weak Minty solutions (cf. \Cref{def:SWMVI}) is formally introduced, and the convergence of \ref{eq:PPPA-intro} is studied for a class of operators that admit at least one such solution. 
\Cref{sec:convergence DRS} establishes a primal-dual connection between \ref{eq:PPPA-intro} and \ref{eq:DRS}, which leads to convergence results for \ref{eq:DRS} under an oblique weak Minty assumption on the associated primal-dual operator.
In \Cref{sec:semi} we introduce the class of semimonotone operators, for which we provide several calculus rules and examples.
Sufficient conditions for the convergence of \ref{eq:DRS}, based on the semimonotonicity of the underlying operators are provided in \Cref{sec:drs:semi}.
In \Cref{subsec:semi:examples} various example problems covered by our theory for \ref{eq:DRS} are presented. Finally, \Cref{sec:conclusion} concludes the paper.
For the sake of readability, several proofs and auxiliary results are deferred to the Appendix.

    \subsection{Notation}\label{subsec:notation}

The set of natural numbers including zero is denoted by \(\N\coloneqq\set{0,1,\hdots}\).
The set of real and extended-real numbers are denoted by \(\R\coloneqq(-\infty,\infty)\) and \(\Rinf\coloneqq\R\cup\set\infty\), while the positive and strictly positive reals are \(\R_+\coloneqq[0,\infty)\) and \(\R_{++}\coloneqq(0,\infty)\).
We denote the positive part of a real number by $[\cdot]_{+} \coloneqq \max\{0, \cdot\}$ and the negative part by $[\cdot]_{-} \coloneqq \min\{0, \cdot\}$.
With \(\id\) we indicate the identity function \(x\mapsto x\) defined on a suitable space.
We denote by $\R^n$ the standard $n$-dimensional Euclidean space with inner product $\langle\cdot,\cdot\rangle$ and induced norm $\|\cdot\|$. 
For a vector $w = (w_1, \hdots, w_N) \in \R^n$, $w_i \in \R^{n_i}$ is used to denote its $i$-th (block) coordinate.
The identity matrix is denoted by $\I_n\in\R^{n\times n}$; we write $\I$ when no ambiguity occurs. 
Given a matrix $P \in \R^{n \times n}$, we denote the range of $A$ by $\range{A}$ and the kernel of $A$ by $\ker{A}$.
The sets of symmetric, symmetric positive semidefinite and symmetric positive definite $n$-by-$n$ matrices are denoted by $\sym{n}$, $\sym{n}_{+}$ and $\sym{n}_{++}$, respectively. We also write $P\succeq0$ and $P\succ0$ for $P\in \sym{n}_{+}$ and $P\in \sym{n}_{++}$, respectively. 
We say that matrices $P, Q \in \R^{n \times n}$ are similar if there exists a nonsingular matrix $X \in \R^{n \times n}$ such that $P = X^{-1} Q X$.
    For any matrix $P \in \R^{n \times n}$ with real eigenvalues, we denote its smallest eigenvalue by $\lambda_{\rm min}(P)$.
Given a symmetric matrix $P\in\sym{n}$ we define the scalar product $\langle x,y \rangle_P=\langle x,Py \rangle$ and the quadratic function \(\qindef{x}{P} \coloneqq \inner{x, x}_P\). If $P\in\sym{n}_{++}$ then its associated induced norm is defined as 
$\|x\|_P
= \sqrt{\qindef{x}{P}}$.  
We denote the Kronecker product between two matrices of arbitrary size by \(\otimes\).

We use the notation 
$
    \seq{
        z^k
        }[k\in I]
$
to denote a sequence with indices in the set $I\subseteq \N$. When dealing with scalar sequences we use the subscript notation $\seq{\gamma_k}[k\in I]$.
We say that \(\seq{z^k}\) converges (sublinearly) with 
    a big-$\pazocal{O}\bigl(g(k)\bigr)$ rate if there exists 
    an index $k_0 \in \N$
    and 
    a positive scalar $c$ such that
    \(
        |z^k| \leq c g(k)
    \)
    for all $k \geq k_0$
    and we say that \(\seq{z^k}\) converges (sublinearly) with
    a little-$o\bigl(g(k)\bigr)$ rate if
    \(
        \nicefrac{z^k}{g(k)} \rightarrow 0    
    \).
We say that 
\(
    \seq{z^k}
\)
converges to a point 
\(
    z^\star
\)
(at least) \DEF{\(Q\)-linearly} with \(Q\)-factor given by \(c\in(0,1)\) if there exists an index $k_0 \in \N$ such that
\(
    \|z^{k+1}-z^\star\|\leq c\|z^k-z^\star\|
\)
for all \(k\geq k_0\).
We say that \(\seq{z^k}\) converges to a point \(z^\star\) (at least) \DEF{\(R\)-linearly} if there exists a sequence of nonnegative scalars \(\seq{v^k}\) such that 
\(
    \|z^{k}-z^\star\|\leq v^k
\)
and \(\seq{v^k}\) converges \(Q\)-linearly to zero.


An operator or set-valued mapping $A:\R^n\rightrightarrows\R^d$ maps each point $x\in\R^n$ to a subset $A(x)$ of $\R^d$. We will use the notation $A(x)$ and $Ax$ interchangeably. 
 We denote the domain of $A$ by $\dom A\coloneqq\{x\in\R^n\mid Ax\neq\emptyset\}$,
its graph by $\graph A\coloneqq\{(x,y)\in\R^n\times \R^d\mid y\in Ax\}$, and 
the set of its zeros by $\zer A\coloneqq\{x\in\R^n \mid 0\in Ax\}$. 
The inverse of $A$ is defined through its graph,
i.e.,
$\graph A^{-1}\coloneqq\{(y,x)\mid (x,y)\in\graph A\}$,
and we denote the range of $A$ by $\range A \coloneqq \dom A^{-1}$.
The \emph{resolvent} of $A$ is defined 
as
$J_A\coloneqq(\id+A)^{-1}$.
The composition of two operators \(A\) and \(B\) is denoted by \(A \circ B\).



\begin{definition}[parallel sum of operators]
    The parallel sum between operators $\A,\B : \R^n\rightrightarrows\R^n$ is defined as $\A \Box \B \coloneqq (\A^{-1} + \B^{-1})^{-1}$.
\end{definition}

\begin{definition}[parallel sum of extended-real numbers]
    Let $a,b \in \Rinf$. We say that $a$ and $b$ are parallel summable if either $a = b = 0$ or $a + b \neq 0$ and their parallel sum is defined as 
    \[
        a \Box b \coloneqq
        \begin{cases}
            0, & \qquad\text{if } a = b = 0,\\
            \frac{ab}{a+b}, & \qquad\text{otherwise},
        \end{cases}
    \]
    where we use the convention that $a \Box \infty = a$.
\end{definition}

\begin{definition}[(co)monotonicity]
An operator $A:\R^n\rightrightarrows\R^n$ is said to be $\mon$-monotone for some $\mon\in \R$ if
\[
    \langle x-\other{x},y-\other{y}\rangle \geq \mon\|x-\other{x}\|^2, \qquad \text{for all $(x,y),(\other{x},\other{y})\in\graph A$},
\]
and it is said to be $\com$-comonotone for some $\com\in \R$ if
\[
    \langle x-\other{x},y-\other{y}\rangle \geq \com\|y-\other{y}\|^2, \qquad \text{for all $(x,y),(\other{x},\other{y})\in\graph A$}.
\]
The operator $A$ is said to be maximally 
(co-)monotone if its graph is not strictly contained in the graph of another 
(co-)monotone operator.
\end{definition}
We say that \(A\) is \emph{outer semicontinuous (osc)} at $\bar x\in\dom A$ if 
\ifspringer
    \(
        \limsup_{x\to \bar x} Ax \coloneqq \{y\mid \exists x^k \to \bar x, \exists y^k \to y \textrm{ with } y^k\in Ax^k\} \subseteq A\bar x.
    \)
\else
    \begin{equation}
        \limsup_{x\to \bar x} Ax \coloneqq \{y\mid \exists x^k \to \bar x, \exists y^k \to y \textrm{ with } y^k\in Ax^k\} \subseteq A\bar x.
    \end{equation}
\fi
Outer semicontinuity of $A$ everywhere is equivalent to its graph being a closed subset of $\R^n\times \R^d$.



The domain of an extended real-valued function $f : \R^n \rightarrow \Rinf$ is the set $\dom f \coloneqq \set{x \in \R^n}[f(x) < \infty]$. We say that $f$ is proper if $\dom f \neq \emptyset$ and that $f$ is lower semicontinuous (lsc) if the epigraph $\epi f \coloneqq \set{(x,\alpha) \in \R^n \times \R}[f(x) \leq \alpha]$ is a closed subset of $\R^{n+1}$. We denote the limiting subdifferential of $f$ by $\partial f$. We say that $f$ is $\ell$\hyp{}smooth to indicate that $f$ is continuously differentiable and $\nabla f$ is Lipschitz continuous with modulus $\ell$.

The indicator function of a set $E\subseteq\R^n$ is denoted by \(\indicator_E\), namely \(\indicator_E(x)=0\) if \(x\in E\) and \(\infty\) otherwise. We denote the normal cone of $E$ by \(N_E\). 
The projection onto and the distance from $E$ with respect to $\|\cdot\|_{\projR}$, \(Q\in\sym{n}_{++}\), are denoted by 
\ifspringer
    \(
        \proj^Q_E(x)\coloneqq{}\argmin_{z\in E}\{\|z-x\|_{\projR}\}
    \)
    and
    \(
        \dist_{\projR}(x,E)\coloneqq{} \inf_{z\in E}\{\|z-x\|_{\projR}\},
    \)
\else
    \begin{align*}
        \proj^Q_E(x)\coloneqq{}\argmin_{z\in E}\{\|z-x\|_{\projR}\}, \quad \dist_{\projR}(x,E)\coloneqq{} \inf_{z\in E}\{\|z-x\|_{\projR}\},
    \end{align*}
\fi    
respectively. The absence of super/subscript $Q$ implies the same definitions with respect to the canonical norm.

  \section{The preconditioned proximal point method}
  \label{sec:PPA}

In this section, the convergence of \ref{eq:PPPA-intro} with positive semidefinite preconditioning is studied for a class of nonmonotone operators, defined through an oblique weak Minty assumption.

We begin by introducing this class of nonmonotone operators and detailing our underlying assumptions. Then, leveraging a projective interpretation of the proximal point algorithm, we provide the corresponding convergence results of \ref{eq:PPPA-intro}. Subsequently, we demonstrate the tightness of our convergence results through a simple example. Finally,
we establish (local) linear convergence under an additional metric subregularity assumption.

      \subsection{Convergence analysis of PPPA under oblique weak Minty}\label{subsec:convergence PPA}


Consider the following class of operators with oblique weak Minty solutions, which generalizes the class of operators with weak Minty solutions that has been employed in the context of the extragradient method and forward-backward-forward splitting \cite{diakonikolas2021Efficient,pethick2021escaping,bohm2022solving}. 

\begin{definition}[$\DRSRho$\hyp{}oblique weak Minty solutions]\label{def:SWMVI}
    An operator \(T : \R^n \rightrightarrows \R^n\) is said to have $\DRSRho$\hyp{}oblique weak Minty solutions at (a nonempty set) \(\pazocal{S}^\star\subseteq\zer T\) for some symmetric matrix $\DRSRho\in\sym{n}$ if
    \begin{equation}\label{def:WMVI}
         \langle 
			v, z-z^\star
		\rangle
			\geq
		\qindef{v}{\DRSRho},
		\qquad
		\text{for all $z^\star\in\pazocal{S}^\star, (z,v)\in\graph T$}, 
    \end{equation}
    where the quadratic form \(\qindef{v}{\DRSRho} \coloneqq \langle v, Vv \rangle\).
	Whenever $V = \rho \I$ for some $\rho \in \R$, we
		refer to them as
	$\rho$\hyp{}weak Minty solutions.
\end{definition}
The generalization of the weak MVI to include a general symmetric matrix $\DRSRho$ instead of a scalar value $\rho$
permits a more detailed characterization of problem classes, which in turn will lead to
tight results for the convergence of \ref{eq:PPPA-intro} and \ref{eq:DRS} (cf. \Cref{sec:convergence DRS}).

The above definition can be further relaxed by requiring \eqref{def:WMVI} to hold instead on \((z,v)\in \graph T \cap (\range{(\M+T)^{-1}\M} \times \range{\M})\), where $\M$ is the preconditioner from \ref{eq:PPPA-intro}. Under this weaker condition, the results from \cref{thm:pppa,thm:PPPA:linear} are still satisfied since \eqref{def:WMVI} is invoked only at points in this restricted set (see \eqref{eq:PPPA-VWMVI-application}). We refrain from 
using
this relaxation given that it couples the problem class with the preconditioner and does not appear to lead to tighter sufficient conditions in terms of semimonotonicity (cf. \Cref{sec:drs:semi}).

As previously mentioned in the introduction, a key aspect of our forthcoming convergence analysis for \ref{eq:PPPA-intro} is to consider positive semidefinite preconditioners, largely inspired by 
the convergence analysis of \cite[Thm. 3.4]{latafat2017Asymmetric} and \cite[Thm. 3.3]{condat2013primal}. 
Unlike the aforementioned works, we do not assume that \(T\) is the sum of a maximally monotone and a linear skew-symmetric operator. Instead, we require that it admits at least a single oblique weak Minty solution.
More specifically, we work under the following assumptions.
\begin{assumption} \label{ass:PPPA}
The operator \(T\) in \eqref{prob:P1} and the preconditioner \(\M\) in \eqref{eq:PPPA-intro} satisfy the following properties. 
	\begin{enumeratass}
		\item \label{ass:PPPA:0} \(T:\R^n\rightrightarrows\R^n\) is outer semicontinuous.
		\item \label{ass:PPPA:0.5} The preconditioned resolvent \((\M+T)^{-1}\M\) has full domain. 
		\item \label{ass:PPPA:1} 
		There exists a nonempty set \(\pazocal{S}^\star\subseteq \zer T\) and a symmetric, possibly indefinite matrix \(V \in \sym{n}\) such that \(T\) has $\DRSRho$\hyp{}oblique weak Minty solutions at \(\pazocal{S}^\star\) for \(V\). 
		\item \label{ass:PPPA:2} \(\M \in \sym{n}_{+}\) is a symmetric positive semidefinite matrix
		such that
	    \begin{equation}\label{eq:PPPA:eigcond}
			\etamin
				{}\coloneqq{}
			1 + \lambda_{\mathrm{min}}(\tp U \DRSRho \M U)
				>
			0,
	    \end{equation}
		where \(U\) is any orthonormal basis for the range of \(\M\).
	\end{enumeratass}
\end{assumption}
Note that the preconditioned resolvent 
\(\PRES\) can be seen as a particular instance of the warped/nonlinear resolvent 
considered in the monotone setting in
\cite{bui2020warped,giselsson2021nonlinear}.

\begin{remark}\label{rem:ass:PPPA:eigcond}
		\cref{ass:PPPA:2} provides a restriction on the preconditioner $\M$, depending on the oblique weak Minty matrix $\DRSRho$. To better understand this condition, consider the matrix \(X \coloneqq \tp U \M U\), which is positive definite by construction (see \cref{it:lem:PQ:properties:decomposition}). By similarity transformation using the similarity matrix \(X^{\nicefrac{1}{2}}\) \cite[Cor. 1.3.4]{horn2012matrix}, it follows that $\tp U \DRSRho \M U$ is similar to the symmetric matrix $X^{\nicefrac{1}{2}}\tp U \DRSRho U X^{\nicefrac{1}{2}}$, i.e., that
		\begin{align*}
			\etamin
				{}={}
			1 + \lambda_{\mathrm{min}}(X^{\nicefrac{1}{2}}\tp U \DRSRho U X^{\nicefrac{1}{2}}).
			\numberthis\label{eq:rem:PPPA:eigcond:similarity}
		\end{align*}
		As a result, \cref{ass:PPPA:2} holds if and only if
		\begin{align*}
			(U^\top \M U)^{-1} + U^\top \DRSRho U \succ 0.
			\numberthis\label{eq:rem:PPPA:eigcond:alt}
		\end{align*}
		This condition is vacuously satisfied in the case where $U^\top \DRSRho U \succeq 0$, which includes the WMI setting.
		The preconditioner $\M$ needs to be selected properly only if $U^\top \DRSRho U$ is not positive semidefinite. 
		For instance, when $\M = \nicefrac1\gamma \I$ where $\gamma > 0$, \cref{ass:PPPA:2} reduces to $\gamma > [-\lambda_{\mathrm{min}}(\DRSRho)]_+$.
		Moreover, when $\DRSRho = \DRSrho \I$, it matches the stepsize condition $\gamma > [-\DRSrho]_+$ for cohypomonotone operators (see e.g. \cite[Thm. 3.1]{combettes2004proximal} and \Cref{tab:PPPA:connections} for more details).
	\end{remark}

Our convergence analysis relies on a projective interpretation of the preconditioned proximal point algorithm in the monotone setting that dates back to \cite{solodov1996Modified,solodov1999hybrid,konnov1997Class}.
In particular, it was shown that when the preconditioner is positive definite,
each iteration of \ref{eq:PPPA-intro} can be interpreted as a projection onto a certain halfspace, constructed 
using the sequences $(z^k)_{k\in\N}$ and $(\bar{z}^k)_{k\in\N}$.
In particular, under the assumption that $T$ is maximally monotone and $P = \gamma^{-1} \I$, \cite{solodov1999hybrid} considers the halfspace
	\begin{equation}
		\mathcal{D}_{z, \bar z} \coloneqq \set{r \in \R^n}[
			\langle \gamma^{-1} (z - \bar z), \bar z - r \rangle
			\ge 0].
	\end{equation}
We will show that this idea can be extended to our more general setting provided that a different halfspace is considered.
More specifically, each iteration of \ref{eq:PPPA-intro} is interpreted as a projection onto a 
halfspace
separating the
current iterate from a subset of solutions that satisfy the oblique weak Minty assumption. 
This halfspace will be constructed using the so-called shadow sequences $(w^k)_{k\in\N} \coloneqq (\projQ z^k)_{k\in\N}$ and $(\bar{w}^k)_{k\in\N} \coloneqq (\projQ\bar{z}^k)_{k\in\N}$ to deal with the fact that in our work the preconditioner is only positive \emph{semidefinite}.

We proceed to present this separating halfspace and its main properties in the following lemma, for which a supporting visualization is provided in \Cref{fig:hyperplane-interpretation}. Note that the introduction of the matrix $\projR \coloneqq \M + \proj_{\ker{\M}}$ is inspired by the convergence analysis performed in \cite[Thm. 3.4]{latafat2017Asymmetric} and \cite[Thm. 3.3]{condat2013primal}.

\begin{figure}
    \centering
    \includetikz{Hyperplane/hyperplane-interpretation-I}
	\ifspringer
		\hspace{-0.6cm}%
	\fi
    \includetikz{Hyperplane/hyperplane-interpretation}
	\ifspringer
		\hspace{-0.6cm}%
	\fi
    \includetikz{Hyperplane/hyperplane-interpretation-indef}
    \caption{
        Visualization of the involved variables in \Cref{lem:hyperplane} for different types of preconditioners, where $\DRSRho = -\tfrac32\I$, $\alpha$ is as in 
		\ifspringer
			\Cref{half:xbar:2}
		\else
			\eqref{eq:hyperplane:alpha}
		\fi
		and $\lambda = \nicefrac\alpha2$. Here, we use the shorthand notation $w_\Pi = \proj_{\mathcal{D}_{w, \bar w}}^\projR(w)$ and $z_\Pi = \proj_{\mathcal{D}_{z, \bar z}}^\M(z)$.
        In these figures, the vector $\bar z$ satisfies the update rule \eqref{eq:PPPA-intro} (see \eqref{eq:PPPA-proof-vbar-incl}).
        Note that the precise location of the halfspace is governed by the oblique weak Minty matrix $V$, due to the dependence of $\alpha$ on $V$. For example, if $V = 0$, then $\alpha = 1$ and the corresponding hyperplane would pass through $\bar{w}$.
        (left and middle) Positive definite preconditioner, such that $\projQ = \I$, $\M = Q$ and the variables $z, \bar{z}, z^+, z_\Pi$ are equal to the variables $w, \bar{w}, w^+, w_\Pi$.
        (left) Multiple of identity preconditioner $\M = \tfrac13 \I$. (middle) $\M = \diag(\nicefrac13, \nicefrac16)$.
        (right) Positive semidefinite preconditioner $\M = \diag(\nicefrac13, 0)$.
        }
    \label{fig:hyperplane-interpretation}
\end{figure}

\begin{lemma}[halfspace interpretation]\label{lem:hyperplane}
    Suppose that \cref{ass:PPPA} holds. Consider the points \(z \in \R^n\), \(\bar z \in \inv{(\M+T)} \M z\) and \(z^+ = z + \lambda (\bar z - z)\) satisfying update rule \eqref{eq:PPPA-intro}, and let $w \coloneqq \projQ z$, $\bar w \coloneqq \projQ\bar z$ and $w^+ \coloneqq \projQ z^+$. Define 
    \ifspringer\else
        the matrix
    \fi
    $\projR \coloneqq \M + \proj_{\ker{\M}}$ and the halfspace
    \begin{equation}\label{eq:hyperplane}
        \mathcal{D}_{w, \bar w} \coloneqq \set{r \in \R^n}[
            \tinprod{w - \bar w}{\bar w - r}_{\projR}
            \ge \qindef{\projR(w - \bar w)}{\DRSRho}].
    \end{equation}
    Then, the following hold.
    \begin{enumerate}
        \item\label{half:QS} The set \(\projQ\pazocal{S}^\star\) is a subset of \(\mathcal{D}_{w, \bar w}\).
        \item\label{half:xbar} If \(w \in \mathcal{D}_{w, \bar w}\) then 
        \(\M(z - \bar z) = 0\), and in particular
\(\bar z \in \zer T\).
        \def\wtil{\tilde w}
        \item\label{half:xbar:2} 
        If \(\bar z \notin \zer T\), then \(\qindef{z - \bar z}{\M} \neq 0\) and 
        \ifspringer
        \(
            \alpha 
                {}\coloneqq{}
            1 + \frac{\qindef{\M(z - \bar z)}{\DRSRho}}{\qindef{z - \bar z}{\M}} 
                {}\geq{}
            \etamin > 0.
        \)
        \else
            \begin{align*}
                \alpha 
                {}\coloneqq{}
                1 + \frac{\qindef{\M(z - \bar z)}{\DRSRho}}{\qindef{z - \bar z}{\M}} 
                    {}\geq{}
                \etamin > 0.\numberthis\label{eq:hyperplane:alpha}
            \end{align*}
        \fi
        \item\label{half:update} 
                $
            w^+
                {}={}
            (1 - \nicefrac{\lambda}{\alpha})w + \nicefrac{\lambda}{\alpha} \proj_{\mathcal{D}_{w, \bar w}}^\projR(w)
        $
        for any $\lambda \in \R$.
    \end{enumerate}

    \begin{proof}

\begin{proofitemize}
    \item \ref{half:QS}: By update rule \eqref{eq:PPPA-intro} it holds that
    \begin{equation} \label{eq:PPPA-proof-vbar-incl}
                \M(z - \bar z) \in T \bar z.
    \end{equation}
    Therefore, owing to the $\DRSRho$-oblique weak Minty assumption (\cref{ass:PPPA:1}), for all $z^\star \in \pazocal{S}^\star$ 
    \begin{equation}\label{eq:PPPA-VWMVI-application}
        \tinprod{\M(z - \bar z)}{\bar z - z^\star} \ge 
        \qindef{\M(z - \bar z)}{\DRSRho}.
    \end{equation}
    As a result, for all \(w^\star \coloneqq \projQ z^\star \in \projQ\pazocal{S}^\star\subseteq \projQ\zer T\)
    \begin{equation*}
        \begin{aligned}
            \tinprod{\projR(w - \bar w)}{\bar w - w^\star}
                {}={}&
            \tinprod{\M(z - \bar z)}{\bar z - z^\star}
            \stackrel{\eqref{eq:PPPA-VWMVI-application}}{\geq}
            \qindef{\M(z - \bar z)}{\DRSRho}
                {}={}
            \qindef{\projR(w - \bar w)}{\DRSRho},
        \end{aligned}
    \end{equation*}
    where we used \Cref{it:lem:PQ:properties:equivalences} to relate $\M$ and $\projR$ in the first and final equality.
    It is now evident that by construction, \(\projQ\pazocal{S}^\star \subseteq \mathcal{D}_{w, \bar w}\). 
    \item \ref{half:xbar}:
    Suppose that \(w\in\mathcal{D}_{w, \bar w}\)
        and
    define the shorthand notation \(\tilde{w} \coloneqq w - \bar w\), \(\tilde{z} \coloneqq z - \bar z\)%
    .
        Owing to \cref{it:lem:PQ:properties:proj,it:lem:PQ:properties:equivalences}, it holds that
    \begin{equation}\label{eq:PX}
        \M 
            {}={}
        \projR\projQ
            {}={}
        \projQ\projR\projQ
            {}={}
        U \tp U \projR U \tp U 
            {}={}
        U X \tp U.
    \end{equation}
        Consequently, by definition of \(\mathcal{D}_{w, \bar w}\), it holds that
    \def\sqrtlam{X^{\nicefrac12}}
    \begin{align*}
        0 \geq \nrm{\tilde{w}}_{\projR}^2 + \qindef{\projR\tilde{w}}{\DRSRho}
            {}={}&
        \tinprod{\projR\projQ\tilde{w}}{(\I + \DRSRho\projR\projQ)\tilde{w}}
        \\
            {}={}&
        \tinprod{U X \tp U \tilde{w}}{(\I + \DRSRho U X \tp U)\tilde{w}}
            {}={}
        \tinprod{X^{\nicefrac{1}{2}} \tp U \tilde{w}}{X^{\nicefrac{1}{2}} \tp U(\I + \DRSRho U X \tp U)\tilde{w}}\\
            {}={}&
        \tinprod{X^{\nicefrac{1}{2}} \tp U \tilde{w}}{(\I + X^{\nicefrac{1}{2}}\tp U\DRSRho U X^{\nicefrac{1}{2}}) X^{\nicefrac{1}{2}} \tp U\tilde{w}}.
        \numberthis \label{eq:PPPA-ineq-posdef-X}
    \end{align*}
    Furthermore, due to \cref{ass:PPPA:2} 
        and \cref{rem:ass:PPPA:eigcond} it follows that
    \(\I + X^{\nicefrac{1}{2}}\tp U\DRSRho U X^{\nicefrac{1}{2}}\) is positive definite, which combined with \eqref{eq:PPPA-ineq-posdef-X} and \(X \succ 0\) implies that \(\tp U\tilde{w} = 0\).
    In turn, using \cref{it:lem:PQ:properties:equivalences,it:lem:PQ:properties:proj} we have
    \(
                \M\tilde{z} = \M \projQ\tilde{w} = \M U \tp U \tilde{w} = 0,
    \)
    which by \eqref{eq:PPPA-proof-vbar-incl} implies that \(\bar z \in \zer T\). 

    \item \ref{half:xbar:2}:
    If \(\bar{z} \notin \zer T\), then inclusion \eqref{eq:PPPA-proof-vbar-incl} implies that \(\M(z - \bar z)\neq 0\). Combined with \(\projR(w-\bar w) = \M(z- \bar z)\) (where we used \Cref{it:lem:PQ:properties:equivalences}) and positive definiteness of \(\projR\) (see \Cref{it:lem:PQ:properties:pd}) it follows that 
    \(
        \qindef{z - \bar z}{\M} 
            = 
                \|w- \bar w\|_{\projR}^2
            > 
        0
    \),
    ensuring that
    \(
        \alpha 
            {}={}
        1 + \tfrac{\qindef{\projR\tilde w}{\DRSRho}}{\tnormsq{\tilde w}_{\projR}}
    \)
    is finite-valued.  
    Consequently, using \eqref{eq:PPPA-ineq-posdef-X} and that \(\tnormsq{\tilde w}_{\projR} = \tinprod{\projR\projQ\tilde{w}}{\tilde{w}} = \tnormsq{\sqrtlam \tp U \tilde{w}}\) by \eqref{eq:PX}, \(\alpha\) may be written as a Rayleigh quotient \cite[Thm. 4.2.2]{horn2012matrix} of \(\I + X^{\nicefrac{1}{2}}\tp U\DRSRho U X^{\nicefrac{1}{2}}\) and \(\sqrtlam \tp U \tilde{w}\), i.e.,
    \begin{align*}
        \alpha 
            {}={}
        1 + \tfrac{\qindef{\projR\tilde w}{\DRSRho}}{\tnormsq{\tilde w}_{\projR}}        
            {}={}
        \tfrac{\tinprod{X^{\nicefrac{1}{2}} \tp U \tilde{w}}{(\I + X^{\nicefrac{1}{2}}\tp U\DRSRho U X^{\nicefrac{1}{2}}) X^{\nicefrac{1}{2}} \tp U\tilde{w}}}{\tnormsq{\sqrtlam \tp U \tilde{w}}} 
            {}\geq{}&
        \lambda_{\mathrm{min}}(\I + X^{\nicefrac{1}{2}}\tp U\DRSRho U X^{\nicefrac{1}{2}})
        \stackrel{\eqref{eq:rem:PPPA:eigcond:similarity}}{=}
        \etamin.
    \end{align*}

    \item \ref{half:update}:
    For any \(u\notin \mathcal{D}_{w, \bar w}\), it holds that 
    \(
        \proj^{\projR}_{\mathcal{D}_{w, \bar w}}(u)
            {}={}
        u 
            {}+{}  
        \tfrac{%
            \langle w-\bar{w},u - \bar{w}\rangle_{\projR} + \qindef{\projR(w-\bar{w})}{\DRSRho}
            }{%
            \|w-\bar{w}\|^2_{\projR}
            }%
            (\bar{w}-w)
    \)
    \cite[Ex. 29.20]{bauschke2017Convex}.
    Therefore, 
    \(
        \proj^{\projR}_{\mathcal{D}_{w, \bar w}}(w)
            {}={}
        w+\alpha(\bar{w}-w)
    \)
        and for any $\lambda \in \R$ it holds that
        \ifspringer
            \(
                (1-\nicefrac{\lambda}{\alpha})w+\nicefrac{\lambda}{\alpha}\proj_{\mathcal{D}_{w, \bar w}}^\projR(w)
                    {}={}
                w+\lambda(\bar{w}-w) 
                    {}={}
                \projQ \big(z
                + 
                \lambda(\bar{z}-z)
                \big) 
                {}={}
                    \projQ z^+
                    {}={}
                w^+.
            \)
        \else
            \begin{align*}
                (1-\nicefrac{\lambda}{\alpha})w+\nicefrac{\lambda}{\alpha}\proj_{\D_{w, \bar w}}^\projR(w)
                    {}={}
                w+\lambda(\bar{w}-w) 
                {}={}
                \projQ \big(z
                + 
                \lambda(\bar{z}-z)
                \big) 
                {}={}
                    \projQ z^+
                    {}={}
                w^+.
            \end{align*}
        \fi
        \qedhere
\end{proofitemize}
    \end{proof}
\end{lemma}
Most notably, \Cref{half:update} establishes that the update rule for the (shadow) sequence generated by \ref{eq:PPPA-intro} can be interpreted as a relaxed projection onto the halfspace defined in \eqref{eq:hyperplane}, while \Cref{half:QS,half:xbar} imply that if at any iteration the current iterate belongs to the halfspace containing the set of projected oblique weak Minty solutions \(\projQ\pazocal{S}^\star\), this implies its optimality.
Owing to the non-expansiveness of the projection operator, this leads to the following convergence result for the preconditioned proximal point algorithm.

\begin{theorem}[convergence of \ref{eq:PPPA-intro}]\label{thm:pppa}
    Suppose that \cref{ass:PPPA} holds, and consider a sequence $\seq{z^k, \bar z^k}$ generated by \ref{eq:PPPA-intro} starting from $z^0 \in \R^n$ with relaxation parameters $\lambda_k \in \left(0,2\alpha_k\right)$ such that $\liminf_{k\to\infty}\lambda_k(2\alpha_k-\lambda_k)>0$, where
    \begin{equation}\label{eq:pppa:alpha_k}
        \alpha_k
            {}={}
        1+\frac{\qindef{\M(z^k-\bar{z}^k)}{\DRSRho}}{\qindef{z^k-\bar{z}^k}{\M}}.
    \end{equation}
    Then, either  a point $\bar z^k \in \zer T$ is reached in a finite number of iterations or 
    the following hold for the sequence $\seq{z^k, \bar z^k}$.
    \begin{enumerate}
            \item\label{it:pppa:v} $\bar{v}^k:=\M(z^k-\bar{z}^k)\in T\bar{z}^k$ for all \(k\) and $(\bar{v}^k)_{k\in N}$ converges to zero.
            \item\label{it:pppa:subseq} Every limit point (if any) of $\seq{\bar z^k}$ belongs to $\zer T$.
            \item\label{it:pppa:bounded} The shadow sequences $(\projQ z^k)_{k\in\N}$, $(\projQ\bar{z}^k)_{k\in\N}$ are bounded
            and their limit points belong to $\projQ\zer T$.
            \item\label{it:pppa:rate} If $\lambda_k(2\alpha_k-\lambda_k)\geq \kappa >0$ uniformly for all \(k\) then%
                        \begin{align*}
                \min_{k=0,1,\ldots,N} 
                                        \qindef{z^k - \bar z^k}{\M}
                    {}\leq{}
                                        \frac1{(N+1)\kappa}
                    \qindef{z^0-z^\star}{\M},
                \qquad\text{for all $z^\star\in\pazocal{S}^\star$}.
            \end{align*}
        \end{enumerate}%
                    \begin{enumerate}[resume]
                       \item\label{it:pppa:full} 
            If in \cref{ass:PPPA:1}
            $
                \projQ
                \pazocal{S}^\star
                =
                \projQ
                \zer T
            $, then
            $\seq{\projQ z^k}$ 
                        converges
            to some element of $\projQ\zer T$.
            If additionally $\PRES$ is (single-valued) continuous,
            then $\seq{\bar{z}^k}$ converges to some 
                        $
                z^\star \in
                \zer T
            $.
            Finally, if $\lambda_k$ is additionally uniformly bounded in the interval $(0,2)$, then $\seq{z^k}$ also converges to 
                        $
                z^\star \in
                \zer T
            $.
        \end{enumerate}
    \begin{proof}

Note that as long as $\bar{z}^k \notin \zer T$ it follows from \cref{half:xbar:2} that $\alpha_k \geq \etamin > 0$. 
Consider the shadow sequences $(w^k)_{k\in\N} \coloneqq (\projQ z^k)_{k\in\N}$ and $(\bar{w}^k)_{k\in\N} \coloneqq (\projQ\bar{z}^k)_{k\in\N}$ with corresponding update rule
\begin{equation}\label{eq:PPPA-proof-w-update}
    w^{k+1} = w^k + \lambda_k(\bar w^k - w^k).
\end{equation}
If $w^k \in \mathcal{D}_{w^k, \bar w^k}$ for some $k \in \N$ then by \Cref{half:xbar} the algorithm has reached a point $\bar{z}^k \in \zer T$ in a finite number of iterations. Therefore, we will consider the case when $w^k \notin \mathcal{D}_{w^k, \bar w^k}$ in the remainder of the proof.
\begin{proofitemize}
    \item \ref{it:pppa:v}:
    From the update rule~\eqref{eq:PPPA-intro}, it directly follows that
    \(
        \bar v^k \coloneqq \M(z^k - \bar{z}^k) \in T \bar{z}^k    
    \)
    for all $k$. It remains to show that $\seq{\bar v^k}$ converges to zero.
    By \Cref{half:update}, update step~\eqref{eq:PPPA-proof-w-update} is equivalent to the relaxed projection
    $
        \tilde \proj_k \coloneqq
        (1-\zeta_k)\id+\zeta_k
                \proj^\projR_{\mathcal{D}_{w^k, \bar w^k}},
    $
    where
    $\zeta_k \coloneqq \nicefrac{\lambda_k}{\alpha_k}
    $%
        . Owing to firm nonexpansiveness of $\proj^\projR_{\mathcal{D}_{w^k, \bar w^k}}$ \cite[Prop. 4.16]{bauschke2017Convex}, the relaxed projection $\tilde \proj_k$ is $\nicefrac{\zeta_k}{2}$-averaged
        in the space with inner product \(\langle \cdot,\cdot \rangle_{\projR}\) \cite[Cor. 4.41]{bauschke2017Convex}.
    Note that for any $w^\star\in \projQ\pazocal{S}^\star$ it holds by \cref{half:QS} that $w^\star\in \mathcal{D}_{w^k, \bar w^k}$. Therefore, it follows from~\cite[Prop. 4.35(iii)]{bauschke2017Convex} for any $w^\star\in \projQ\pazocal{S}^\star$ that
    \begin{align}\label{eq:iFejer}
        \|w^{k+1}-w^\star\|_{\projR}^2
        {}={}
        \|\tilde \proj_k(w^{k})-\tilde \proj_k(w^\star)\|_{\projR}^2
        &\leq
        \|w^k-w^\star\|_{\projR}^2-\tfrac{1-\nicefrac{\zeta_k}{2}}{\nicefrac{\zeta_k}{2}}\|w^{k+1}-w^k\|_{\projR}^2\nonumber\\
        &=\|w^k-w^\star\|_{\projR}^2-\tfrac{2\alpha_k-\lambda_k}{\lambda_k}\,\|\lambda_k(w^k-\bar{w}^k)\|_{\projR}^2\nonumber\\
        &=\|w^k-w^\star\|_{\projR}^2-\lambda_k(2\alpha_k-\lambda_k)\,\|w^k-\bar{w}^k\|_{\projR}^2,
    \end{align}
    establishing that $(w^k)_{k\in\N}$ is Fej\'er monotone with respect to $\projQ\pazocal{S}^\star$ \cite[Def. 5.1]{bauschke2017Convex}.
    By telescoping \eqref{eq:iFejer}, and since
    $\liminf_{k\to\infty}\lambda_k(2\alpha_k-\lambda_k)>0$,
    it follows that 
    \(\sum_{k = 0}^\infty \|w^k - \bar w^k\|^2_{\projR} < +\infty\), which implies that 
    \ifspringer
        \(
            \|w^k - \bar w^k\|_{\projR} \to 0.
        \)
    \else
        \begin{equation}\label{eq:wres0}
            \|w^k - \bar w^k\|_{\projR} \to 0.
        \end{equation}
    \fi
    The convergence of $(\bar{v}^k)_{k\in N}$ to zero is established by noting that
    \begin{align*}
        \|\bar{v}^k\|^2
        {}={}
        \|\M(w^k-\bar{w}^k)\|^2
      {}\leq{}
        \|\M^{\nicefrac12}\|^2\, \|\M^{\nicefrac12}(w^k-\bar{w}^k)\|^2
        {}={}
        \nrm{\M}\, \|w^k-\bar{w}^k\|^2_{\projR}.
        \numberthis\label{eq:vw}
    \end{align*}

    \item \ref{it:pppa:subseq}: Since \((\bar v^k)_{k\in\N}\) converges to zero and $\bar{v}^k\in T\bar{z}^k$, it follows from outer semicontinuity of $T$ (due to \Cref{ass:PPPA:0}) that any limit point of $(\bar{z}^k)_{k\in\N}$ belongs to $\zer T$.


    \item \ref{it:pppa:bounded}:
    It follows from \eqref{eq:iFejer} that \(\seq{\|w^k - w^\star\|_{\projR}}\) converges, and in particular that $(w^k)_{k\in\N}$ is bounded. In turn, using 
    \ifspringer
        that $\|w^k - \bar w^k\|_{\projR} \to 0$
    \else
        \eqref{eq:wres0}
    \fi
    and the triangle inequality 
    \(
        \|\bar{w}^k\|_{\projR} \leq \|\bar{w}^k-w^k\|_{\projR}+\|w^k\|_{\projR},
    \)
    it follows that \(\seq{\bar{w}^k}\) is bounded 
        and thus that it has at least one limit point.
        Take a subsequence $(\bar{w}^k)_{k\in K}$ converging to some limit point $w^\infty$. Since $(\bar{w}^k-w^k)_{k\in\N}$ converges to zero
        \ifspringer
            as shown in \Cref{it:pppa:v},
        \else
            by \eqref{eq:wres0},
        \fi
        we have that $(w^k)_{k\in K}$ also converges to the same limit point $w^\infty$.
        The claim is established by noting that $w^\infty \in \projQ \zer T$ owing to \cref{it:pppa:subseq}.

    \item \ref{it:pppa:rate}:
    By telescoping \eqref{eq:iFejer} and from the uniform bound
    \(\forall k \in \N : \lambda_k(2\alpha_k-\lambda_k) \ge \kappa > 0\), it follows that
    \begin{equation}\label{eq:iFejer-telescoped}
        \begin{aligned}
            \min_{0\le k\le N} \|w^k - \bar w^k\|^2_{\projR} &\le \tfrac{1}{(N+1)\kappa}\|w^0 - w^\star\|^2_{\projR}.
        \end{aligned}
    \end{equation}
    The claimed upper bound
        thus follows by noting that
        $\|w^k-\bar w^k\|_{\projR}^2=\qindef{z^k-\bar z^k}{\M}$ and
        $\|w^0-w^\star\|_{\projR}^2=\qindef{z^0-z^\star}{\M}$.

    \item \ref{it:pppa:full}:
        If
        $\projQ \pazocal{S}^\star=\projQ \zer T$,
    then the convergence of $(\projQ z^k)_{k\in\N}$
        to some element of $\projQ\zer T$ follows by~\cite[Thm. 5.5]{bauschke2017Convex}. 
    Let $G = \PRES$ be (single-valued) continuous, so that
        $\bar{z}^k = G z^k = G w^k$.
        Take a subsequence $(w^k)_{k\in K}$ converging to some limit point $w^\infty$. Then,
    by continuity of $G$, the sequence $(\bar{z}^k)_{k\in\N}$ converges to $\bar{z}^\infty = Gw^\infty$, and by \cref{it:pppa:subseq} it holds that \(\bar{z}^\infty\in \zer T\). 

    It remains to show that $(z^k)_{k\in\N}$ also converges to 
        $\bar{z}^\infty$
    if $\lambda_k$ is uniformly bounded in the interval $(0,2)$. Note that
    \ifspringer
        \(
            \nrm{z^{k+1} - \bar{z}^\infty} \leq \lambda_k \nrm{\bar{z}^k - \bar{z}^\infty} + |1-\lambda_k| \nrm{z^k - \bar{z}^\infty}.
        \)
    \else
        \[
            \nrm{z^{k+1} - \bar{z}^\infty} \leq \lambda_k \nrm{\bar{z}^k - \bar{z}^\infty} + |1-\lambda_k| \nrm{z^k - \bar{z}^\infty}.
        \]
    \fi
    Consequently, by uniform boundedness of $\lambda_k$, i.e., $(\lambda_k)_{k\in\N} \subseteq [\epsilon, 2-\epsilon]$ for some $\epsilon$ > 0, and by 
        \cite[Lem. 3 of \S 2.2]{polyak2020Introduction}
    it follows that $\lim_{k \rightarrow \infty} \nrm{z^{k+1} - \bar{z}^\infty} = 0$, i.e., that $(z^k)_{k\in\N}$ also converges to $\bar{z}^\infty \in \zer T$.
    \qedhere
\end{proofitemize}

    \end{proof}
\end{theorem}

\begin{table}
    \centering
    \caption{
        Connection between \Cref{thm:pppa} and existing convergence results for \ref{eq:PPPA-intro}.
    }
    \label{tab:PPPA:connections}
    \footnotesize
    \begin{tabular}{l}
            {
            \ifspringer
                \begin{tblr}{
                    hlines,
                    vlines,
                    vline{1-1} = {1-2}{0pt},
                    hline{1-2} = {1-1}{0pt},
                    hline{5-5} = {2-2}{0pt},
                    hline{7-7} = {2-2}{0pt},
                    hline{9-9} = {2-2}{0pt},
                    hline{11-11} = {2-2}{0pt},
                    hline{7-9} = {4-4}{0pt},
                    hline{11-11} = {4-4}{0pt},
                    vline{3-3} = {2-12}{dotted},
                    vline{5-5} = {2-12}{dotted},
                    hline{4-4} = {2-3}{dotted},
                    hline{8-8} = {2-3}{dotted},
                    colspec={M{0.7cm}M{1.7cm}M{3.1cm}M{2.2cm}M{2.6cm}},
                    cell{1-2}{2-5}={c,gray!20},
                    cell{3-12}{1-1}={c,gray!20}
                    }
                    &\SetCell[c=2]{c}$\M = \gamma^{-1}\I$&&
                    \SetCell[c=2]{c}$\M \succeq 0$\\
                    &Conditions&Theorems&Conditions&Theorems\\
                    \SetCell[r=3]{c}$\DRSRho = 0$&$\gamma > 0, \lambda_k=1$&
                    \cite{martinet1970regularisation}\textsuperscript{1},
                    \cite{rockafellar1976monotone}\textsuperscript{12},
                    \cite[Thm. 3.2]{allevi2006proximal}&
                    \SetCell[r=3]{c}$\lambda_k \in (0,2)$&
                    \SetCell[r=3]{c}\cite[Thm. 3.4]{latafat2017Asymmetric}\textsuperscript{1}, \cite[Thm. 2.9]{bredies2022degenerate}\textsuperscript{1}\\
                    &$\gamma > 0$,&\SetCell[r=2]{c}\cite[Thm. 3]{eckstein1992DouglasRachford}\textsuperscript{12}&\\
                    &$\lambda_k \in (0,2)$&&\\
                    \SetCell[r=4]{c}$\DRSRho = \DRSrho \I$&$\gamma > [-2\DRSrho]_+$,&
                    \SetCell[r=2]{c}\cite[Thm. 9]{pennanen2002local}\textsuperscript{12},
                    \cite[Thm. 1]{iusem2003Inexact}\textsuperscript{12},
                    \cite[Thm. 3.1]{gorbunov2023convergence}&&\SetCell[r=4]{c}\cref{thm:pppa}\\
                    &$\lambda_k=1$&&$1 + \lambda_{\mathrm{min}}(\DRSrho \tp U \M U) > 0$,&\\
                    &$\gamma > [-\DRSrho]_+$,&\SetCell[r=2]{c}\cite[Thm. 3.1]{combettes2004proximal}\textsuperscript{12}&$\lambda_k \in (0, 2\alpha_k)$&\\
                    &$\lambda_k \in \Bigl(0,2\bigl(1+\tfrac{\DRSrho}{\gamma}\bigr)\Bigr)$&&\\
                    \SetCell[r=2]{c}$\DRSRho \in \sym{}$&$\gamma > [-\lambda_{\rm min}(\DRSRho)]_+$,&\SetCell[r=2]{c}\cref{thm:pppa}&$1 + \lambda_{\mathrm{min}}(\tp U \DRSRho \M U) > 0$,&\SetCell[r=2]{c}\cref{thm:pppa}\\
                    &
                    $\lambda_k \in (0, 2\alpha_k)$
                    &&$\lambda_k \in (0, 2\alpha_k)$&
                \end{tblr}
            \else
                \begin{tblr}{
                    hlines,
                    vlines,
                    vline{1-1} = {1-2}{0pt},
                    hline{1-2} = {1-1}{0pt},
                    hline{5-5} = {2-2}{0pt},
                    hline{7-7} = {2-2}{0pt},
                    hline{9-9} = {2-2}{0pt},
                    hline{11-11} = {2-2}{0pt},
                    hline{7-9} = {4-4}{0pt},
                    hline{11-11} = {4-4}{0pt},
                    vline{3-3} = {2-12}{dotted},
                    vline{5-5} = {2-12}{dotted},
                    hline{4-4} = {2-3}{dotted},
                    hline{8-8} = {2-3}{dotted},
                    colspec={M{0.9cm}M{2.1cm}M{3.3cm}M{2.6cm}M{2.8cm}},
                    cell{1-2}{2-5}={c,gray!20},
                    cell{3-12}{1-1}={c,gray!20}
                    }
                    &\SetCell[c=2]{c}$\M = \gamma^{-1}\I$&&
                    \SetCell[c=2]{c}$\M \succeq 0$\\
                    &Conditions&Theorems&Conditions&Theorems\\
                    \SetCell[r=3]{c}$\DRSRho = 0$&$\gamma > 0, \lambda_k=1$&
                    \cite{martinet1970regularisation}\textsuperscript{1},
                    \cite{rockafellar1976monotone}\textsuperscript{1,2},
                    \cite[Thm. 3.2]{allevi2006proximal}&
                    \SetCell[r=3]{c}$\lambda_k \in (0,2)$&
                    \SetCell[r=3]{c}\cite[Thm. 3.4]{latafat2017Asymmetric}\textsuperscript{1}, \cite[Thm. 2.9]{bredies2022degenerate}\textsuperscript{1}\\
                    &$\gamma > 0$,&\SetCell[r=2]{c}\cite[Thm. 3]{eckstein1992DouglasRachford}\textsuperscript{1,2}&\\
                    &$\lambda_k \in (0,2)$&&\\
                    \SetCell[r=4]{c}$\DRSRho = \DRSrho \I$&$\gamma > [-2\DRSrho]_+$,&
                    \SetCell[r=2]{c}\cite[Thm. 9]{pennanen2002local}\textsuperscript{1,2},
                    \cite[Thm. 1]{iusem2003Inexact}\textsuperscript{1,2},
                    \cite[Thm. 3.1]{gorbunov2023convergence}&&\SetCell[r=4]{c}\cref{thm:pppa}\\
                    &$\lambda_k=1$&&$1 + \lambda_{\mathrm{min}}(\DRSrho \tp U \M U) > 0$,&\\
                    &$\gamma > [-\DRSrho]_+$,&\SetCell[r=2]{c}\cite[Thm. 3.1]{combettes2004proximal}\textsuperscript{1,2}&$\lambda_k \in (0, 2\alpha_k)$&\\
                    &$\lambda_k \in \Bigl(0,2\bigl(1+\tfrac{\DRSrho}{\gamma}\bigr)\Bigr)$&&\\
                    \SetCell[r=2]{c}$\DRSRho \in \sym{}$&$\gamma > [-\lambda_{\rm min}(\DRSRho)]_+$,&\SetCell[r=2]{c}\cref{thm:pppa}&$1 + \lambda_{\mathrm{min}}(\tp U \DRSRho \M U) > 0$,&\SetCell[r=2]{c}\cref{thm:pppa}\\
                    &
                    $\lambda_k \in (0, 2\alpha_k)$
                    &&$\lambda_k \in (0, 2\alpha_k)$&
                \end{tblr}
            \fi
            }\vspace{1pt}\\
            \textsuperscript{1}\footnotesize{Considers maximal monotonicity instead of MVI or maximal $\DRSrho$\hyp{}comonotonicity instead of $\DRSrho$\hyp{}weak Minty.}\\
            \textsuperscript{2}\footnotesize{Considers the setting where the stepsize $\gamma_k$ is variable and the resolvent computations are performed inexactly.}
        \end{tabular}
\end{table}

Note that as long as the solution has not yet been reached, i.e., as long as \(\bar z^k \notin \zer T\), then $\alpha_k \geq \etamin > 0$ by \Cref{half:xbar:2}.
    Therefore, all claims of \Cref{thm:pppa} also hold for the (more restrictive) fixed stepsize range $\lambda_k \in (0, 2\etamin)$. 

    Due to the generality of our underlying assumptions,
    \Cref{thm:pppa} unifies and extends existing convergence results for \ref{eq:PPPA-intro} in literature, as summarized in \Cref{tab:PPPA:connections}.
    One notable instance of \Cref{thm:pppa} is when $\DRSRho = \DRSrho \I$ and $\M = \gamma^{-1} \I$, in which case $\alpha_k = \etamin = 1 + \nicefrac{\DRSrho}{\gamma}$.
    For completeness, we present the simplified theorem statement for this setting below. 
    \begin{corollary}\label{cor:pppa:wmi}
        Suppose that $T$ is outer semicontinuous, that there exists a nonempty set \(\pazocal{S}^\star\subseteq \zer T\) such that \(T\) has $\DRSrho$\hyp{}weak Minty solutions at \(\pazocal{S}^\star\) and that the resolvent $J_{\gamma T}$ has full domain given stepsize $\gamma > [-\DRSrho]_+$.
        Consider a sequence $\seq{z^k, \bar z^k}$ generated by \ref{eq:PPPA-intro} starting from $z^0 \in \R^n$ with preconditioner $\M = \nicefrac1\gamma\I$ and relaxation parameters $\lambda_k \in (0,2(1+\nicefrac{\DRSrho}{\gamma}))$ satisfying $\liminf_{k\to\infty}\lambda_k(2(1+\nicefrac{\DRSrho}{\gamma})-\lambda_k)>0$.
        Then, all the claims from \Cref{thm:pppa} hold.
    \end{corollary}

It is worth reiterating that 
all results
from \cref{thm:pppa}
except the final one apply even when only a subset of $\zer T$ are oblique weak Minty solutions.
Most notably, \Cref{it:pppa:subseq} states that
the limit points of $\seq{\bar z^k}$ are zeros of $T$, which do not have to be oblique weak Minty solutions.
The significance of this result is highlighted by the 
second
part of the following example, where this convergence behavior is observed in practice (see also \Cref{fig:ex:pppa:toy}).
The first
part of this example 
demonstrates the tightness of our results 
by considering a simple linear inclusion related to saddle point problems.
    The details are deferred to \Cref{detail:ex:pppa:toy}.
    \ifspringer\else
        Supplementary code for the numerical examples presented in this work can be found on \href{https://github.com/brechtevens/Minty-DRS-examples}{GitHub}.
    \fi

\ifspringer
    \def\widthone{0.465\textwidth}
    \def\widthtwo{0.53\textwidth}
\else
    \def\widthone{0.195\textwidth}
    \def\widthtwo{0.6\textwidth}
\fi

\begin{figure}
    \centering
    \includetikz{Examples/Toy/iterates}
    \caption{
        Visualization of the sequences $\seq{\bar z^k}$ for the relaxed proximal point algorithm 
        with stepsize $\gamma = 1$ and fixed
            relaxation parameter $\lambda$
        applied to \Cref{ex:pppa:toy}, where $a = 2$ and $b = 1$. 
        The zeros of $T$ are marked in green and its gradient flow $-T$ is indicated using gray arrows.
        (left) Sequences generated using $\lambda = 2.3 < \bar \lambda$, where $\bar \lambda$ is the upper bound on $\lambda_k$ from \eqref{eq:ex:pppa:toy:lambda}. The orange sequence converges to the single $\bigl(\nicefrac{b}{a^2 + b^2}\bigr)$\hyp{}weak Minty solution $(0,0)$ while the blue sequence converges to a zero of $T$ which is not the weak Minty solution (cf. \Cref{it:pppa:subseq}).
        (right) Sequences generated using $\lambda = 2.5 > \bar \lambda$. In this setting, it is no longer guaranteed by \Cref{thm:pppa} that the sequence $\seq{\bar z^k}$ converges, as can be seen from the (diverging) blue sequence.
        }
    \label{fig:ex:pppa:toy}
\end{figure}

\begin{example}[toy example]\label{ex:pppa:toy}
    \ifspringer
        Consider the \newline
    \else\fi
    \begin{minipage}[t]{\widthone}
        \ifspringer
            operator
        \else
            Consider the operator
        \fi
        \begin{align*}
            \ifspringer\else
                \hspace{-3cm}
            \fi
            T(z) \coloneqq f\Bigl(\nrm{z}\Bigr) 
            \begin{bmatrix}
                b & a \\
                -a & b
            \end{bmatrix}z,
        \end{align*}
    \end{minipage}\hfill%
    \begin{minipage}[t]{\widthtwo}
        \centering
        \ifspringer
            \vspace{-0.75cm}
        \else
            \vspace{-0.4cm}
        \fi
        \includetikz{Examples/Toy/helper}
        \captionof{figure}{%
                \ifspringer
                    Def.
                \else
                    Definition
                \fi
                of $f : \R_+ \rightarrow [0, 1]$ in \Cref{ex:pppa:toy}.
            }
        \label{fig:ex:pppa:toy:f}
    \end{minipage}

    \noindent
    where
        $f : \R_+ \rightarrow [0, 1]$
        and $a, b \in \R$ are not simultaneously equal to zero.
    \begin{enumerate}
        \item \label{ex:pppa:toy:f1}
        If $f \equiv 1$, then the following hold.
        \begin{enumerate}
            \item \label{it:ex:pppa:toy:f1:minty}
            $T$ is $\bigl(\nicefrac{b}{a^2 + b^2}\bigr)$\hyp{}comonotone, and thus
            $T$ has a single $\bigl(\nicefrac{b}{a^2 + b^2}\bigr)$\hyp{}weak Minty solution at $\pazocal{S}^\star = \zer T = {(0,0)}$. 
            \item \label{it:ex:pppa:toy:f1:convergence} 
            By \Cref{it:pppa:full}, the sequence $\seq{\bar z^k}$ generated by the proximal point algorithm with fixed
            relaxation parameter $\lambda$ converges to $\zer T$ if
                        \(\gamma > \frac{[-b]_+}{a^2+ b^2}\)
                and \(\lambda\)
                lies in the interval
                \begin{align*}
                    \begin{aligned}
                        0 < \lambda < 
                        \bar{\lambda} \coloneqq 2\left(1+\tfrac{b}{
                            \gamma(a^2 + b^2)
                            }\right).
                    \end{aligned}\numberthis\label{eq:ex:pppa:toy:lambda}
                \end{align*}
            By examining the spectral radius of the algorithmic operator, it can be seen that this result is tight.
        \end{enumerate}
        \item 
        If $f$ is defined as in \Cref{fig:ex:pppa:toy:f}
        and $b$ is nonnegative%
        , then the following hold.
        \begin{enumerate}
            \item \label{it:ex:pppa:toy:f2:minty} 
            $T$ has a single $\bigl(\nicefrac{b}{a^2 + b^2}\bigr)$\hyp{}weak Minty solution at the point $\pazocal{S}^\star = {(0,0)} \subset \zer T$ while the other zeros of $T$ are not $\rho$-weak Minty solutions for any $\rho \in \R$.
            \item \label{it:ex:pppa:toy:f2:convergence} 
            By \Cref{it:pppa:subseq}, the limit points of the sequence $\seq{\bar z^k}$ generated by the 
            relaxed
            proximal point algorithm with fixed 
            relaxation parameter $\lambda$ belong to $\zer T$ (not necessarily to $\pazocal{S}^\star$) if 
            $\gamma > 0$ and
            $\lambda$ 
                        is selected according to \eqref{eq:ex:pppa:toy:lambda}.
            This result has been verified numerically in \Cref{fig:ex:pppa:toy}.
        \end{enumerate}
    \end{enumerate}
    \ifspringer
            \else
            \fi
\end{example}

As a second illustrative example, we will consider the application of our theory to the economic equilibrium model of Von Neumann \cite{neumann1945model}.
For this example, it was recently observed in \cite[Prop. 2]{daskalakis2020independent} that for certain choices $R$ and $S$, operator $F$ is nonmonotone but satisfies the weak MVI.
Here, we extend this result by showing that the Nash equilibrium of \eqref{eq:ex:neumann} is a weak Minty solution of the full operator $T \coloneqq F + N_C$.
\ifspringer
    The details can be found in the arXiv preprint \cite[Appendix C]{evens2023convergence}.
\else
    The details are deferred to \Cref{detail:ex:VonNeumann}.
\fi

\begin{figure}
    \centering
    \includetikz{Examples/VonNeumann/iterates_pppa}%
    \caption{
        (left)
            Visualization of the sequences $\seq{z^k}$ and $\seq{\bar z^k}$ for the relaxed proximal point algorithm, applied to the minimax problem from \Cref{ex:VonNeumann}, with
            problem parameter $\epsilon = \nicefrac1{10}$, weak Minty constant $\DRSrho = \tfrac{\epsilon-1}{4} = -\nicefrac{9}{40}$, stepsize $\gamma = \nicefrac14 > [-\DRSrho]_+$ and relaxation parameter $\lambda = \nicefrac13 \in \bigl(0,2(1+\nicefrac{\DRSrho}{\gamma})\bigr)$.
            These choices comply with \cref{cor:pppa:wmi}.
            The domain $z \in \Delta_2 \times \Delta_2$ is parametrized by $z = (s, 1-s, t, 1-t)$.
            The heat map indicates the value of
            $
                \varrho_\epsilon\bigl(z\bigr)
                    =
                \infi_{v \in N_C(z)}
                \nicefrac{
                    \inner*{F(z) + v,z-z^\star}
                }
                {
                    \nrm{F(z) + v}^2
                }
            $,
            which can be thought of as a local measure of nonmonotonicity of the operator (see also \Cref{detail:ex:VonNeumann}).
            The region where $\varrho_\epsilon\bigl(z\bigr)$ is positive is marked in green.
            Even though the iterates pass through the nonmonotone region,
            both sequences converge to the Nash equilibrium.
            (right)
            Visualization of the distance between the sequences $\seq{z^k}$ and $\seq{\bar z^k}$ and the Nash equilbrium $z^\star$. 
            Note that in the first iterations, when passing through the nonmonotone region, the sequence $\|\bar{z}^k - z^\star\|$ increases. On the other hand, the sequence $\|z^k - z^\star\|$ is monotonically decreasing (see also \eqref{eq:iFejer}).
    }
    \label{fig:VonNeumann}
\end{figure}

\begin{example}[Von Neumann's economic equilibrium model]\label{ex:VonNeumann}
    Using the standard definition of the simplex
    \(
    \Delta_d \coloneqq \set{z \in \R^d}[ z \geq 0, \sum_{i=1}^d z_i = 1]
    \),
    consider an economy where there are 
    $n$ goods with relative prices $y \in \Delta_n$,
    which can be produced by $m$ processes with relative intensities $x \in \Delta_m$.
    Let $R_{ij} \geq 0$ and $S_{ij} \geq 0$ denote the number of units of the $j$th good produced and consumed by the $i$th process, respectively.
    By \cite[\S 5]{neumann1945model}, a pair $(x^\star, y^\star)$ for which this economy is expanding corresponds to a minimax solution of
    \begin{align*}
        \minimize_{x \in \Delta_m} 
        \;
        \maximize_{y \in \Delta_n} 
        \;
        f(x,y) \coloneqq \frac{\inner{x, Ry}}{{\inner{x, Sy}}}.
        \numberthis\label{eq:ex:neumann}
    \end{align*}
    Defining $z \coloneqq (x, y)$, the problem of finding stationary points of \eqref{eq:ex:neumann} corresponds to finding a zero of
    $
        T(z)
            {}\coloneqq{}
        F(z)
        +
        N_C(z)
    $,
    where
    $
        F(z)
            {}={}
        \begin{psmallmatrix}
            \nabla_x f(x,y) \\
            - \nabla_y f(x,y)
        \end{psmallmatrix}
    $
    and 
    $
        C = \Delta_2 \times \Delta_2
    $.
    Let $\epsilon \in (0,1)$ and let
    \begin{align*}
        R 
            {}={} 
        \begin{bmatrix}
            0 & 1 + \nicefrac{\epsilon}2\\
            2 - \nicefrac{\epsilon}2 & 2
        \end{bmatrix}
        \quad\text{and}\quad
        S
            {}={}
        \begin{bmatrix}
            \nicefrac12 & \nicefrac12\\
            1 & 1
        \end{bmatrix}.
        \numberthis\label{eq:prop:ex:neumann:matrices}
    \end{align*}
    Then, the unique Nash equilibrium 
    $
        z^\star \coloneqq (x^\star, y^\star) =
        \begin{bsmallmatrix}
            0 & 1 & 0 & 1
        \end{bsmallmatrix}^\top
    $
    of \eqref{eq:ex:neumann} is a $\Bigl(\tfrac{\epsilon-1}{4}\Bigr)$\hyp{}weak Minty solution of $T$.
    A numerical experiment leveraging this result can be found in \Cref{fig:VonNeumann}.
\end{example}

      \subsection{Last-iterate convergence of PPPA under comonotonicity}\label{subsec:convergence PPA:comon}

The convergence result from \Cref{it:pppa:rate} is expressed in terms of a best-so-far $\pazocal{O}\bigl(\nicefrac{1}{(N+1)}\bigr)$ rate for $\seq{\qindef{z^k - \bar z^k}{\M}}$. To obtain stronger last-iterate convergence rates, governing inequality \eqref{def:WMVI} needs to hold between any two points in the graph of $T$ instead of only between any point in the graph of $T$ and (a subset of) $\zer T$. 
This observation leads us to introduce the class of $\DRSRho$-comonotone operators.
\begin{definition}[$\DRSRho$\hyp{}comonotone operator]\label{def:matrix-comon}
    An operator \(T : \R^n \rightrightarrows \R^n\) is said to be $\DRSRho$\hyp{}comonotone for some symmetric matrix $\DRSRho\in\sym{n}$ if
    \begin{equation}\label{eq:matrix-comon}
         \langle x - \other{x},y-\other{y}\rangle\geq \qindef{y-\other{y}}{\DRSRho},\qquad \text{for all $(x,y), (\other{x},\other{y})\in\graph T$}.
    \end{equation}
    It is said to be maximally \(\DRSRho\)\hyp{}comonotone if its graph is not strictly contained in the graph of another \(\DRSRho\)\hyp{}comonotone operator.
\end{definition}

This definition naturally reduces to the notion of (maximally) $\rho$-comonotone operators when $\DRSRho = \rho\I$.
The class of operators with $\DRSRho$\hyp{}oblique weak Minty solutions is significantly larger than the class of (maximally) $\DRSRho$\hyp{}comonotone operators.
For instance, the following proposition establishes that the class of maximally $\DRSRho$\hyp{}comonotone operator only consist of operators which do not have isolated zeros.
\begin{proposition}
    If an operator \(T : \R^n \rightrightarrows \R^n\) is maximally $\DRSRho$\hyp{}comonotone, then $\zer T$ is convex.
    \begin{proof}
        First, observe that 
        $
            \zer\, (T^{-1}+\DRSRho)^{-1}
                =
            (T^{-1}+\DRSRho)(0)
                = 
            \zer T
        $.
        Furthermore, 
        note that by definition $T^{-1}+\DRSRho$ is maximally monotone, and thus so is $(T^{-1}+\DRSRho)^{-1}$. Since the set of zeros of a maximally monotone operator is convex \cite[Prop. 23.39]{bauschke2017Convex}, the claim is established.
    \end{proof}
\end{proposition}
For maximally $\DRSRho$\hyp{}comonotone operators, the convergence of \ref{eq:PPPA-intro} depends merely upon the selection of the preconditioner $\M$ and the relaxation parameters $\lambda_k$.
This is a consequence of the following proposition, which allows to directly verify the underlying assumptions of \Cref{thm:pppa}.

\begin{proposition}\label{prop:comon:implication}
    Suppose that \(T\) is maximally $\DRSRho$\hyp{}comonotone, that $\zer T$ is nonempty and that \(\M\) is selected according to \Cref{ass:PPPA:2}. Then, \Cref{ass:PPPA} holds.
    \begin{proof}
        Let $(x^k, y^k) \in \graph T$, where \(x^k \rightarrow \bar{x}\) and $y^k \rightarrow \bar{y}$.
        Note that $x^k \in (T^{-1} - \DRSRho)^{-1}(y^k - \DRSRho x^k)$.
        Since $T$ is maximally $\DRSRho$\hyp{}comonotone, $(T^{-1} - \DRSRho)^{-1}$ is maximally monotone and thus also outer semicontinuous \cite[Ex. 12.8]{rockafellar2009Variational}. Consequently, it follows that $\bar{x} \in (T^{-1} - \DRSRho)^{-1}(\bar{y} - \DRSRho \bar{x})$.
        This implies that \(\bar{x} \in T \bar{y}\), showing outer semicontinuity of \(T\). 
        \Cref{ass:PPPA:0.5} follows directly from 
        \Cref{prop:comon:fulldom}
        and \Cref{ass:PPPA:1} holds by definition of $\DRSRho$-comonotonicity and using that $\zer T$ is nonempty.
    \end{proof}
\end{proposition}

Leveraging \Cref{prop:comon:implication}, sublinear convergence rates can be obtained for the sequence $\seq{\qindef{z^k - \bar z^k}{\M}}$ in the maximally $\DRSRho$\hyp{}comonotone setting by showing that this sequence is monotonically nonincreasing and using the best-iterate result from \Cref{it:pppa:rate}, as detailed below.
\begin{theorem}[last-iterate convergence]\label{thm:pppa:lastiter}
    Suppose that \(T\) is maximally $\DRSRho$\hyp{}comonotone, that $\zer T$ is nonempty and that \(\M\) is selected according to \Cref{ass:PPPA:2}.
    Consider a sequence $\seq{z^k, \bar z^k}$ generated by \ref{eq:PPPA-intro} starting from $z^0 \in \R^n$ with relaxation parameters $\lambda_k \in \left(0,2\etamin\right)$ such that $\lambda_k(2\etamin-\lambda_k)\geq \kappa > 0$ uniformly for all \(k\). Then, for all $z^\star\in\zer T$, the following convergence estimates hold:
    \begin{align*}
        \qindef{z^N - \bar z^N}{\M}
            {}\leq{}
        \frac{\qindef{z^0-z^\star}{\M}}{(N+1)\kappa}
        \quad\text{and}\quad
        \qindef{z^N - \bar z^N}{\M}
            {}={}
        o\bigl(\nicefrac{1}{(N+1)}\bigr).
        \numberthis\label{eq:pppa:last-it:rate}
    \end{align*}
    \begin{proof}
        By update rule \eqref{eq:PPPA-intro} it holds that
        \begin{equation*}
            \M(z^k - \bar z^k) \in T \bar z^k
            \quad\text{and}\quad
            \M(z^{k+1} - \bar z^{k+1}) \in T \bar z^{k+1}.
        \end{equation*}
        Defining $\tilde z^k \coloneqq z^k - \bar z^k$ and using $\DRSRho$-comonotonicity of $T$ between the points $\bar z^k$ and $\bar z^{k+1}$, this implies that
        \begin{align*}
            \inner{\M(\tilde z^k - \tilde z^{k+1}), z^k - z^{k+1}}
                {}\geq{}
            \qindef{\tilde z^k - \tilde z^{k+1}}{\M + \M \DRSRho \M},
        \end{align*}
        which by update rule \eqref{eq:PPPA-proof-w-update} reduces to
        \begin{align*}
            \inner{2\M(\tilde z^k - \tilde z^{k+1}), \tilde z^k}
                {}\geq{}
            \tfrac{2}{\lambda_k}\qindef{\tilde z^k - \tilde z^{k+1}}{\M + \M \DRSRho \M}.
        \end{align*}
        Since
        $\qindef{u}{\M} - \qindef{v}{\M} = \inner{2\M u, u-v} - \qindef{u-v}{\M}$
        for any $u, v \in \R^n$ and using that $\M = U X U^\top$, this implies that
        \begin{align*}
            \qindef{\tilde z^k}{\M} - \qindef{\tilde z^{k+1}}{\M}
                {}={}&
            \inner{2\M \tilde z^k, \tilde z^k - \tilde z^{k+1}} - \qindef{\tilde z^k - \tilde z^{k+1}}{\M}\\
                {}\geq{}&
            \tfrac{2}{\lambda_k}\qindef{\tilde z^k - \tilde z^{k+1}}{\M + \M \DRSRho \M} - \qindef{\tilde z^k - \tilde z^{k+1}}{\M}\\
                {}={}&
            \tfrac{2}{\lambda_k}\qindef{\tilde z^k - \tilde z^{k+1}}{U X^{\nicefrac12}(\I + X^{\nicefrac{1}{2}}\tp U \DRSRho U X^{\nicefrac{1}{2}})X^{\nicefrac12} U^\top} - \qindef{\tilde z^k - \tilde z^{k+1}}{\M}\\			
                \dueto{\eqref{eq:rem:PPPA:eigcond:similarity}}{}\geq{}&
            \Bigl(\tfrac{2\etamin}{\lambda_k} - 1\Bigr)\qindef{\tilde z^k - \tilde z^{k+1}}{\M},
        \end{align*}
        establishing that $\qindef{z^k - \bar z^k}{\M}$ is nonincreasing. Combining this with \Cref{it:pppa:rate} yields the claimed big-$\pazocal{O}$ convergence rate. Finally, the little-$o$ convergence follows from \cite[Lem. 3-(1a)]{davis2015convergence}.
    \end{proof}
\end{theorem}

\Cref{thm:pppa:lastiter} provides a general framework for the sublinear convergence rate of \ref{eq:PPPA-intro}, unifying many existing results in literature.
Among others, we retrieve the big-$\pazocal{O}$ rates for standard PPA (where $\M = \gamma^{-1} \I$ and $\lambda_k = 1$) in the monotone setting (where $\DRSRho = 0$) \cite[Prop. 8]{brezis1978produits} and the weak Minty setting (where $\DRSRho = \DRSrho \I$) \cite[Thm. 3.1]{gorbunov2023convergence}, as well as the rates obtained for \ref{eq:PPPA-intro} in the monotone setting \cite[Thm. 3.4]{latafat2017Asymmetric}.
Moreover, it is worth highlighting that when $\M = \gamma^{-1} \I$, \Cref{thm:pppa:lastiter} indicates that the optimal rate is obtained for $\lambda_k = \etamin$, as opposed to the usual choice of $\lambda_k = 1$.



      \subsection{Linear convergence of PPPA under metric subregularity}\label{subsec:linear PPA}

Metric subregularity is a well-known concept used
to describe the behavior of a set-valued mapping in a neighborhood around a given point and can be seen as a one-point version of metric regularity.
Many algorithms used in optimization, such as splitting methods, have been shown to have good convergence properties when applied to operators that are metrically subregular, see e.g. \cite{tseng2000modified,drusvyatskiy2016error,latafat2017Asymmetric,Latafat2019Randomized,giselsson2021nonlinear}. 
Following \cite[Ex. 3H.4]{dontchev2009implicit}, we define metric subregularity
as follows.

\begin{definition}[metric subregularity]\label{metricsub:Prelim}
  	A set-valued mapping $T:\R^n\rightrightarrows\R^d$ is \emph{metrically subregular} at $\bar{x}$ for $\bar{y}$ if $(\bar{x},\bar{y})\in\graph T$ and there exists a positive constant 
	$
				\zeta
	$
	together with 
  	 a \emph{neighborhood of subregularity} $\mathcal{U}$ of $\bar{x}$ such that 
   	\begin{align}\label{eqnn:metricsubregularity:TriPD}
  		\dist(x,T^{-1}\bar{y})\leq{}&
				\zeta
		\dist(\bar{y},Tx), \qquad \text{for all } x\in \mathcal{U}. 
    \end{align}
\end{definition}

There exist many prominent examples of operator classes that are metrically subregular. For instance, 
all piecewise polyhedral mappings are metrically subregular at all points in their graph \cite[\S 3]{dontchev2009implicit}. Another important example is the subdifferential operator \(\partial f\) of a proper lsc convex function satisfying quadratic growth conditions (cf. \cite{aragon2008characterization}). 
We refer the interested reader to \cite[\S 3]{dontchev2009implicit} and \cite[\S 9]{rockafellar2009Variational} for further discussion.

The main result of this subsection is \Cref{thm:PPPA:linear}, where it is shown that 
local
linear convergence of \ref{eq:PPPA-intro} 
can be attained under an additional metric subregularity assumption.
If the metric subregularity holds globally, i.e., $\mathcal{U} = \R^n$ in \Cref{metricsub:Prelim}, then the linear convergence holds globally.
This result will in turn lead to sufficient conditions for (local) linear convergence of \ref{eq:DRS} for (possibly) nonmonotone piecewise polyhedral mappings (cf. \cref{cor:DRS:CP:lin}).

Interestingly, when $\M = \gamma^{-1}\I$ and $\DRSRho = 0$, \Cref{thm:PPPA:linear} establishes that the sequence \(\seq{\dist(z^k, \zer T)}\) 
	converges \(Q\)-linearly to zero with the same factor as
	the one obtained in the maximally monotone setting \cite[Thm. 3.2(b)]{shen2016linear}. 
	For a more detailed discussion on linear convergence of PPA in the monotone setting, we refer the interested reader to \cite{artacho2007convergence,leventhal2009metric} and the references therein.

\begin{theorem}[linear convergence of \ref{eq:PPPA-intro} under metric subregularity]
\label{thm:PPPA:linear}
	In addition to \cref{ass:PPPA}, suppose that \(\zer T\) is nonempty, closed and convex, that \cref{ass:PPPA:1} holds with \(\pazocal{S}^\star = \zer T\), that $\PRES$ is continuous, and
	that \(T\) is metrically subregular
	with modulus $\zeta \in \R_+$
	at all \(z^\star\in \zer T\) for \(0\).
	Consider a sequence $\seq{z^k, \bar z^k}$ generated by \ref{eq:PPPA-intro} starting from $z^0 \in \R^n$ with relaxation parameters $\lambda_k \in \left(0,2\alpha_k\right)$ such that 
		$\lambda_k(2\alpha_k-\lambda_k) \geq \kappa > 0$ for all \(k\) large enough,
	where $\alpha_k$ is defined as in \eqref{eq:pppa:alpha_k}.
		Then, $\seq{\dist_{\projR}(\projQ z^{k}, \projQ\zer T)}$ converges \(Q\)-linearly to zero with \(Q\)-factor
		$
			\sqrt{ 
				1
					{}-{}
				\nicefrac{\kappa}{(1 + \zeta\nrm{\M})^2} 
			}
		$,
		\(\seq{\projQ z^{k}}\) converges \(R\)-linearly to some element of \(\projQ\zer T\)
		and
		\(\qindef{z^k - \bar z^k}{\M}\) 
		converges \(R\)-linearly to zero.
	\begin{proof}
By \cref{it:pppa:full} the sequence \(\seq{\bar z^k}\) converges to a point \(z^\star\in \zer T\). Up to possibly discarding initial iterates, it holds that \(\bar z^k \in \mathcal U_{z^\star}\), where \(\mathcal U_{z^\star}\) denotes the neighborhood of subregularity associated with \(z^\star\). Metric subregularity of \(T\) at \(z^\star\) for \(0\) along with the fact that \(\M(z^k- \bar z^k)\in T \bar{z}^k\) (cf. \cref{it:pppa:v}) implies that
\ifspringer
	\(
		\dist(\bar z^k, \zer T) 
			{}\leq{}
				\zeta
		\|\M(z^k- \bar z^k)\|, 
	\)
\else
	\begin{equation}\label{eq:metricsub}
		\dist(\bar z^k, \zer T) 
			{}\leq{}
				\zeta
		\|\M(z^k- \bar z^k)\|, 
	\end{equation}
\fi
for some
\(
		\zeta
	\in(0,\infty)
\).
The distance of the shadow sequence from the set \(\projQ\zer T\) may therefore be bounded as
\ifspringer
	\begin{align*}
		\dist(\bar w^k, \projQ\zer T) 
			{}={}
		\inf_{v\in \projQ\zer T} \{\|v - \bar w^k\|\}
			{}={}&
		\inf_{v\in \zer T} \{\|\projQ(v - \bar z^k)\|\}\\
			{}\leq{}&
		\dist(\bar z^k, \zer T) 
			{}\leq{}
				\zeta
		\|\M(z^k- \bar z^k)\|
		\numberthis \label{eq:wMx}
	\end{align*}
\else
	\begin{align*}
		\dist(\bar w^k, \projQ\zer T) 
			{}={}
		\inf_{v\in \projQ\zer T} \{\|v - \bar w^k\|\} 
			{}={}
		\inf_{v\in \zer T} \{\|\projQ(v - \bar z^k)\|\}
			{}\leq{}&
		\|\projQ\|\dist(\bar z^k, \zer T) 
		\\
		\dueto{\eqref{eq:metricsub}}
			{}\leq{}&
				\zeta
		\|\M(z^k- \bar z^k)\|
		\numberthis \label{eq:wMx}
	\end{align*}
\fi
where we used that \(\|\projQ\| = 1\). 
Note that since the set \(\pazocal{S}^\star = \zer T\) is convex, so is \(\projQ\zer T\), implying that 
\(
	\dist(\bar w^k, \projQ\zer T) = \|\bar w^k - \proj_{\projQ\zer T}(\bar w^k)\|
\)
is attained. 
By the triangle inequality 
\ifspringer
	\begin{align*}
		\|w^k - \proj_{\projQ\zer T}(\bar w^k)\|
			{}\leq{}&
		\|\bar w^k - \proj_{\projQ\zer T}(\bar w^k)\| + \|w^k - \bar w^k\|\\
			{}={}&
		\dist(\bar w^k, \projQ\zer T) + \|w^k - \bar w^k\|
		\\
		\dueto{\eqref{eq:wMx}}
			{}\leq{}&
				\zeta
		\|\M(z^k- \bar z^k)\| + \|w^k - \bar w^k\|
		\\
		\dueto{\eqref{eq:vw}}
			{}\leq{}&
		\xi	\|w^k-\bar{w}^k\|_{\projR},
	\end{align*}
\else
	\begin{align*}
		\|w^k - \proj_{\projQ\zer T}(\bar w^k)\|
			{}\leq{}
		\|\bar w^k - \proj_{\projQ\zer T}(\bar w^k)\| + \|w^k - \bar w^k\|
			{}={}&
		\dist(\bar w^k, \projQ\zer T) + \|w^k - \bar w^k\|
		\\
		\dueto{\eqref{eq:wMx}}
			{}\leq{}&
				\zeta
		\|\M(z^k- \bar z^k)\| + \|w^k - \bar w^k\|
		\\
		\dueto{\eqref{eq:vw}}
			{}\leq{}&
		\xi	\|w^k-\bar{w}^k\|_{\projR},
	\end{align*}
\fi
where 
\(
	\xi =
		\zeta
	\nrm{\M}^{\nicefrac12} + \nrm{\M}^{-\nicefrac12}
\). 
It then follows that
\ifspringer
	\begin{align*}
		\dist_{\projR}^2( w^k, \projQ\zer T) 
			{}={}
		\inf_{v \in \projQ\zer T}\{\|w^k-v\|_{\projR}^2\}
			{}\leq{}&
		\|w^k - \proj_{\projQ\zer T}(\bar w^k)\|^2_{\projR}\\
			{}\leq{}&
		\xi^2 \nrm{\M} \|w^k - \bar w^k\|_{\projR}^2.
		\numberthis\label{eq:distw}
	\end{align*}
\else
	\begin{align*}
		\dist_{\projR}^2( w^k, \projQ\zer T) 
			{}={}
		\inf_{v \in \projQ\zer T}\{\|w^k-v\|_{\projR}^2\}
			{}\leq{}
		\|w^k - \proj_{\projQ\zer T}(\bar w^k)\|^2_{\projR}
			{}\leq{}
		\xi^2 \nrm{\M} \|w^k - \bar w^k\|_{\projR}^2.
		\numberthis\label{eq:distw}
	\end{align*}
\fi
Consequently, due to \eqref{eq:iFejer} and \eqref{eq:distw} we obtain for all \(k\) large enough that
\begin{align*}
	\dist_{\projR}^2(w^{k+1}, \projQ\zer T) 
		{}\leq{}&
	\|w^{k+1} - \proj_{\projQ\zer T}^\projR(w^k)\|_{\projR}^2
	\\ 
	\dueto{\eqref{eq:iFejer}}
		{}\leq{}&
	\|w^{k} - \proj_{\projQ\zer T}^\projR(w^k)\|_{\projR}^2 
		{}-{}
	\lambda_k(2\alpha_k-\lambda_k) \|w^k - \bar w^k\|_{\projR}^2
	\\
		{}={}&
	\dist_{\projR}^2(w^{k}, \projQ\zer T) 
		{}-{}
	\lambda_k(2\alpha_k-\lambda_k) \|w^k - \bar w^k\|_{\projR}^2
	\numberthis\label{eq:linnorm}
	\\
	\dueto{\eqref{eq:distw}}
		{}\leq{}&
	\left(%
		1
			{}-{}
		\tfrac{\kappa}{\xi^2\nrm{\M}} 
	\right)
	\dist_{\projR}^2(w^{k}, \projQ\zer T),
\end{align*}
establishing that the sequence \(\seq{\dist_{\projR}^2(w^{k}, \projQ\zer T)}\) converges \(Q\)-linearly to zero.
Therefore, it follows from \eqref{eq:linnorm} that \(\seq{\|w^k - \bar w^k\|_{\projR}}\) converges \(R\)-linearly to zero since
\ifspringer
	\(
		\dist_{\projR}^2(w^{k}, \projQ\zer T) 
			{}\geq{}
		\lambda_k(2\alpha_k-\lambda_k) \|w^k - \bar w^k\|_{\projR}^2
			{}\geq{}
		\kappa\|w^k - \bar w^k\|_{\projR}^2.
	\)
\else
	\begin{align*}
		\dist_{\projR}^2(w^{k}, \projQ\zer T) 
			{}\geq{}
		\lambda_k(2\alpha_k-\lambda_k) \|w^k - \bar w^k\|_{\projR}^2
			{}\geq{}
		\kappa\|w^k - \bar w^k\|_{\projR}^2.
	\end{align*}
\fi
Consequently, by \cref{eq:PPPA-proof-w-update} it follows that
\ifspringer
	\(
		\|w^{k+1} -  w^k\|_{\projR}^2
		{}={}
		\lambda_k^2\|\bar w^{k} -  w^k\|_{\projR}^2
		{}\leq{}
		\tfrac{\lambda_k^2}{\kappa}\dist_{\projR}^2(w^{k}, \projQ\zer T). 
	\)
\else
	\[
		\|w^{k+1} -  w^k\|_{\projR}^2
		{}={}
		\lambda_k^2\|\bar w^{k} -  w^k\|_{\projR}^2
		{}\leq{}
		\tfrac{\lambda_k^2}{\kappa}\dist_{\projR}^2(w^{k}, \projQ\zer T). 
	\]
\fi
Combined with the fact that \(\lambda_k < 2\alpha_k < +\infty\) is bounded (the latter inequality follows directly from \eqref{eq:pppa:alpha_k}), and since the distance on the right-hand side converges \(Q\)-linearly to zero, \(\seq{\|w^{k+1} -  w^k\|_{\projR}}\) converges \(R\)-linearly to zero. 
This implies that \(\seq{\|w^k- w^\star\|_{\projR}}\) (with \(w^\star\in \projQ\zer T\) being the point where the sequence converges owing to \cref{it:pppa:full}) also converges \(R\)-linearly to zero.
Noting that 
\(\xi^2\nrm{\M} = (1 + \zeta\nrm{\M})^2\) and that
\(\projR\) is positive definite completes the proof.
	\end{proof}
\end{theorem}

%

  \section{Douglas-Rachford splitting}\label{sec:convergence DRS}

In the previous section, convergence results of the preconditioned proximal point method 
were
provided for a class of nonmonotone operators.
In this section,
convergence results for \ref{eq:DRS} are obtained in the nonmonotone setting by leveraging the equivalence between \ref{eq:DRS} and \ref{eq:PPPA-intro}. 
    In particular, we provide convergence results whenever the primal-dual operator has oblique weak Minty solution(s),
    as well as in the more restrictive, comonotone setting,
    and
we show (local) linear convergence for piecewise polyhedral mappings.

      \subsection{Equivalence between DRS and PPPA}\label{sec:equivDRS}




In the monotone setting, the convergence of \ref{eq:DRS} has been studied in \cite{condat2013primal} by establishing its connection with the preconditioned proximal point method with a positive semidefinite preconditioner.
In the nonmonotone setting, this equivalence heavily relies upon the abstract duality framework introduced in \cite{attouch1996general}, \cite[\S 6.9]{auslender2006Asymptotic}. In this framework, the inclusion problem \eqref{prob:composite} is referred to as the \emph{primal inclusion}. Associated with the primal inclusion is the \emph{dual inclusion}, given by
\begin{equation}\label{eq:dual}\tag{D-I}
    \text{find} \quad y \in \R^n \quad \text{such that} \quad 
    0\in \Td y 
    {}\coloneqq{}
    -A^{-1}(-y)+B^{-1}(y).
\end{equation}
A key property of the dual inclusion is that the primal inclusion is solvable if and only if the dual inclusion is solvable \cite[Prop. 6.9.1]{auslender2006Asymptotic}.
Finally, the \emph{primal-dual inclusion} is given by
\begin{equation}\label{eq:primaldual}\tag{PD-I}
    \text{find} \quad z = (x,y) \in \R^{2n} \quad \text{such that} \quad 
    \begin{bmatrix}0\\ 0\end{bmatrix}\in \Tpd z \coloneqq \begin{bmatrix}Ax\\ B^{-1}y\end{bmatrix}+\begin{bmatrix}y\\ -x\end{bmatrix}.
\end{equation}
For the primal-dual inclusion, it holds by \cite[Prop. 6.9.2]{auslender2006Asymptotic} that
\begin{align*}
    (x^\star,y^\star) &\in\zer \Tpd \iff (x^\star,- y^\star)\in\graph A,\ (x^\star,y^\star)\in\graph B,
\end{align*}
where \(\zer\Tp = \{x^\star \mid \exists y^\star : (x^\star,y^\star) \in\zer \Tpd\}\) and \(\zer\Td = \{y^\star \mid \exists x^\star : (x^\star,y^\star) \in\zer \Tpd\}\).
Therefore, the primal-dual inclusion encompasses the solutions of both the primal and the dual inclusion, in the sense that $(x^\star,y^\star)$ is a solution of \eqref{eq:primaldual} if and only if $x^\star$ is a solution of \eqref{prob:composite} and $y^\star$ is a solution of \eqref{eq:dual}.
Hence, to obtain a zero of $\Tp$, one may apply \ref{eq:PPPA-intro} to the primal-dual operator $\Tpd$ starting from an initial point $z^0 = (x^0, y^0)$ with iterates $z^k = (x^k, y^k)$ and $\bar z^k = (\bar x^k, \bar y^k)$.
Most notably, \ref{eq:PPPA-intro} applied to the primal-dual operator $\Tpd$ reduces to \ref{eq:DRS} when the preconditioner $\M$ is selected to be equal to
\begin{align}
    \M
        {}={}
    \begin{bmatrix}
        \gamma^{-1}\I_n&-\I_n\\
        -\I_n&\gamma \I_n
    \end{bmatrix},
    \label{eq:M}
\end{align}
where $\gamma \in \R_{++}$ is a strictly positive stepsize parameter. In the following lemma, the equivalence between \ref{eq:PPPA-intro} and \ref{eq:DRS} is formally summarized.
\begin{lemma}[equivalence of \ref{eq:DRS} and \ref{eq:PPPA-intro}]\label{lem:equivalence:PPPA-DRS}
    Let $z^0 = (x^0, y^0) \in \R^{2n}$ and set $s^0 = x^0 - \gamma y^0$. Then, to any sequence \(\seq{u^k, v^k, s^k}\) generated by \ref{eq:DRS} (initialized with $s^0$) for solving \eqref{prob:composite}, there correspond sequences \(\seq{z^k} = \seq{x^k, y^k}\), \(\seq{\bar z^k} = \seq{\bar x^k, \bar y^k}\) generated by \ref{eq:PPPA-intro} (initialized with $z^0$) for solving \eqref{eq:primaldual} with preconditioner \eqref{eq:M} (and vice versa), and the correspondence is as follows.
    \begin{enumerate}
        \item\label{it:equivalence:s} $s^k = x^k - \gamma y^k$.
        \item\label{it:equivalence:xbar} $\bar x^k = u^k$.
        \item\label{it:equivalence:ybar} $\bar y^k = \nicefrac{1}{\gamma}(2u^k - s^k - v^k)$.
    \end{enumerate}
    \begin{proof}
        See \Cref{proof:lem:equivalence:PPPA-DRS}.
    \end{proof}
\end{lemma}

      \subsection{Convergence analysis of DRS under oblique weak Minty}\label{subsec:convergence DRS}

Despite the fact that \ref{eq:DRS} is a primal algorithm, \Cref{lem:equivalence:PPPA-DRS} shows that it can be interpreted as applying \ref{eq:PPPA-intro} to the associated primal-dual operator $\Tpd$.
Based on this equivalence, we proceed to derive convergence results for \ref{eq:DRS}, leveraging the results obtained in \Cref{sec:PPA}.
First, the underlying assumptions for \ref{eq:PPPA-intro} provided in \Cref{ass:PPPA} are translated to conditions on the individual operators $\A$ and $\B$ as follows. 
\begin{assumption} \label{ass:DRS}
In problem \eqref{prob:composite}, the following hold.
    \begin{enumeratass}
        \item \label{ass:DRS:1} Operators \(A\) and \(B\) are outer semicontinuous.
        \item \label{ass:DRS:2} For the selected stepsizes, the resolvents have full domain, i.e., 
        $\dom J_{\gamma A} = \dom J_{\gamma B} = \R^n$.
        \item \label{ass:SWMVIDRS} 
        There exists a nonempty set \(\pazocal{S}^\star \subseteq \zer \Tpd\) such that 
        the primal-dual operator~\eqref{eq:primaldual} has $\DRSRho$\hyp{}oblique weak Minty solutions at \(\pazocal{S}^\star\), where
        \begin{equation}\label{eq:WMSDRS}
            \DRSRho
                {}\coloneqq{}
            \blkdiag(\DRSrhocom \I_n, \DRSrhomon \I_n)
        \end{equation}
        and \(\DRSrhocom,\DRSrhomon, \in\R\) satisfy
        \(
            [\DRSrhomon]_{-}[\DRSrhocom]_{-} < \nicefrac{1}{4}.
        \)
    \end{enumeratass}
\end{assumption}
Note that these assumptions are not very restrictive. For instance, only the existence of a single oblique weak Minty solution is required for the primal-dual inclusion. Furthermore, it can be seen that the range of $\DRSrhocom$ and $\DRSrhomon$ in \Cref{ass:SWMVIDRS} covers cases where neither the primal inclusion \eqref{prob:composite} nor the dual inclusion \eqref{eq:dual} satisfies the Minty variational inequality (MVI), which is an immediate consequence of the following lemma (see also \Cref{tab:MVI}).

\begin{table}
    \centering
    \caption{Table summarizing which operators from the abstract duality framework satisfy the Minty variational inequality (MVI), depending on the sign of $\DRSrhocom$ and $\DRSrhomon$.}
    \label{tab:MVI}
    \begin{tabular}{c}
        {
        \begin{tblr}{
            hlines,
            vlines,
            vline{1-1} = {1-1}{0pt},
            hline{1-1} = {1-1}{0pt},
            colspec={M{2cm}M{2.5cm}M{2.5cm}},
            cell{1-1}{2-3}={c,gray!20},
            cell{2-3}{1-1}={c,gray!20}
            }
            &$\DRSrhomon \geq 0$ & $\DRSrhomon < 0$\\
            $\DRSrhocom \geq 0$&$\Tp$, $\Td$, $\Tpd$&$\Tp$\\
            $\DRSrhocom < 0$&$\Td$&-
        \end{tblr}
        }
    \end{tabular}
\end{table}

\begin{lemma}[weak Minty for primal and dual operator]\label{lem:SWMI:primal-dual}
    Suppose that \Cref{ass:SWMVIDRS} holds and let 
    \ifspringer
        \(
            \pazocal{S}^\star_P \coloneqq \{x^\star \mid \exists y^\star : (x^\star,y^\star) \in \pazocal{S}^\star\}
        \)
        and
        \(
            \pazocal{S}^\star_D \coloneqq \{y^\star \mid \exists x^\star : (x^\star,y^\star) \in \pazocal{S}^\star\}.
        \)
    \else
        \[
            \pazocal{S}^\star_P \coloneqq \{x^\star \mid \exists y^\star : (x^\star,y^\star) \in \pazocal{S}^\star\}
            \quad \text{and} \quad 
            \pazocal{S}^\star_D \coloneqq \{y^\star \mid \exists x^\star : (x^\star,y^\star) \in \pazocal{S}^\star\}.
        \]
    \fi
    Then, the primal operator $\Tp$ has $\DRSrhocom$\hyp{}weak Minty solutions at $\pazocal{S}^\star_P$ and the dual operator $\Td$ has $\DRSrhomon$\hyp{}weak Minty solutions at $\pazocal{S}^\star_D$.
    \begin{proof}
        See \Cref{proof:lem:SWMI:primal-dual}.
    \end{proof}
\end{lemma}

The converse of \Cref{lem:SWMI:primal-dual} does not hold. For instance, consider the linear operators $A = -\I_n$ and $B = \I_n$: in this case both the primal operator 
    $\Tp$
    and the dual operator 
    $\Td$
    reduce to the zero mapping, which
    has $0$-weak Minty solutions everywhere,
    while the primal-dual operator reduces to
    $
        \Tpd =
        \begin{bsmallmatrix}
            -\I_n & \I_n\\
            -\I_n & \I_n
        \end{bsmallmatrix}
    $,
    which does not have any $0$-weak Minty solutions.



The following theorem follows as a consequence of \Cref{thm:pppa} and characterizes the convergence of \ref{eq:DRS} in the nonmonotone setting, assuming the corresponding primal-dual operator has \(V\)-oblique weak Minty solutions.
To prove this result, we rely upon \Cref{lem:implication:assDRS:assPPPA}, which shows that when \Cref{ass:DRS} holds for \ref{eq:DRS}, \Crefrange{ass:PPPA:0}{ass:PPPA:1} are satisfied for the corresponding preconditioned proximal point algorithm. 
Note that \Cref{ass:PPPA:2} is not covered by this lemma.
To ensure satisfaction of this final assumption, the stepsize $\gamma$ of \ref{eq:DRS} needs to be selected appropriately (see condition \eqref{eq:thm:pppaDRS:gamma}).


\iftables
    \begin{table}[H]
        \centering
        \begin{tabular}{c}
            {
            \begin{tblr}{
                hlines,
                vlines,
                vline{1-1} = {1-1}{0pt},
                hline{1-1} = {1-1}{0pt},
                colspec={M{2cm}M{4cm}M{5.5cm}},
                cell{1-1}{2-3}={c,gray!20},
                cell{2-3}{1-1}={c,gray!20}
                }
                &$\DRSrhomon \geq 0$ & $\DRSrhomon < 0$\\
                $\DRSrhocom \geq 0$&$\gamma \in \left(0,+\infty\right)$&$\gamma \in \left(0, -\tfrac{1 + \sqrt{1 - 4\DRSrhomon\DRSrhocom}}{2\DRSrhomon}\right)$\\
                $\DRSrhocom < 0$&$\gamma \in \left(-\tfrac{2\DRSrhocom}{1+\sqrt{1-4\DRSrhomon\DRSrhocom}},+\infty\right)$&$\gamma \in \left(-\tfrac{2\DRSrhocom}{1+\sqrt{1-4\DRSrhomon\DRSrhocom}},-\tfrac{1 + \sqrt{1 - 4\DRSrhomon\DRSrhocom}}{2\DRSrhomon}\right)$
            \end{tblr}
            }
        \end{tabular}
        \caption{Range of $\gamma$ for Douglas-Rachford.}
        \label{tab:gamma}
    \end{table}
\fi

\begin{theorem}[convergence of \ref{eq:DRS}]\label{thm:pppaDRS}
    Suppose that \Cref{ass:DRS} holds and 
        that $\gamma$ lies within the interval
    \begin{align*}
        \left(
            \tfrac{2[-\DRSrhocom]_{+}}{1+\sqrt{1-4\DRSrhomon\DRSrhocom}},\tfrac{1 + \sqrt{1 - 4\DRSrhomon\DRSrhocom}}{2[-\DRSrhomon]_{+}}
        \right)\numberthis\label{eq:thm:pppaDRS:gamma}
    \end{align*}
    (which is nonempty since        
        \(
            [\DRSrhomon]_{-}[\DRSrhocom]_{-} < \nicefrac{1}{4}
        \)).
    Consider a sequence \(\seq{u^k, v^k, s^k}\) generated by \ref{eq:DRS} starting from $s^0 \in \R^n$ with stepsize $\gamma$
    and relaxation parameters $\lambda_k \in \left(0,2\alpha\right)$ such that $\liminf_{k\to \infty}\lambda_k(2\alpha - \lambda_k)>0$, where
    \(
        \alpha \coloneqq 1+\tfrac{1}{\gamma}\DRSrhocom+\gamma\DRSrhomon > 0.
    \)
    Then, either a point
        for which $u^k = v^k$ 
        is reached in a finite number of iterations, implying that $\bigl(u^k,\nicefrac1{\gamma}(u^k - s^k)\bigr) \in \zer \Tpd$, 
        or the following hold for the sequence \(\seq{u^k, v^k, s^k}\).
    \begin{enumerate}
        \item \label{it:DRS:v}  
        $\gamma^{-1}(u^k-v^k)\in Au^k+Bv^k$ for all \(k\), and \(\seq{u^k-v^k}\) converges to zero.
        \item\label{it:DRS:subseq} 
        Every limit point (if any) of 
                $\seq{u^k,\nicefrac1{\gamma}(u^k - s^k)}$
        belongs to $\zer \Tpd$.
        \item\label{it:DRS:bounded} The sequence $(s^k)_{k\in\N}$ is bounded
        and its limit points belong to $\{(x^\star-\gamma y^\star)\mid (x^\star,y^\star)\in\zer \Tpd\}$.
        \item\label{it:DRS:rate} 
        If $\lambda_k(2\alpha-\lambda_k)\geq \kappa >0$ 
        uniformly
        for all \(k\), then, for all 
                $s^\star \in \{(x^\star-\gamma y^\star)\mid (x^\star,y^\star)\in\pazocal{S}^\star\}$
        we have
        \[
            \min_{k=0,1,\ldots,N}\|u^k-v^k\|^2\leq \frac{\|s^0-s^\star\|^2}{(N+1)\kappa}.
        \]
    \end{enumerate}%
        \begin{enumerate}[resume]
                \item\label{it:DRS:full} If in \cref{ass:SWMVIDRS}
        $
            \projQ
            \pazocal{S}^\star
                =
            \projQ
            \zer \Tpd
        $,
        then $(s^k)_{k\in\N}$ converges to some element of $\{(x^\star-\gamma y^\star)\mid (x^\star,y^\star)\in\zer \Tpd\}$.
                If additionally $J_{\gamma A}$ and $J_{\gamma B}$ are (single-valued) continuous, then
            $\seq{u^k,\nicefrac1{\gamma}(u^k - s^k)}$
        converges to some element of $\zer \Tpd$.
    \end{enumerate}

    \begin{proof}

Recall that \ref{eq:DRS} is cast as an instance of \ref{eq:PPPA-intro}, as summarized in \Cref{lem:equivalence:PPPA-DRS}. 
Therefore, all the assertions follow from those of \Cref{thm:pppa} applied to the primal-dual inclusion \eqref{eq:primaldual} with \(P\) as defined in \eqref{eq:M}. 
Moreover, by \Cref{lem:implication:assDRS:assPPPA} it follows immediately that \Crefrange{ass:PPPA:0}{ass:PPPA:1} hold so that it only remains to show that \Cref{ass:PPPA:2} holds whenever $\gamma$ 
\iftables
    is selected according to \Cref{tab:gamma}.
\else
    satisfies \eqref{eq:thm:pppaDRS:gamma}.
\fi

Throughout the proof, \(\Tpd\) is as in \eqref{eq:primaldual}.
The preconditioning matrix \(\M\) has two distinct eigenvalues, namely zero and $\gamma+\gamma^{-1}$, each with a multiplicity of \(n\). Therefore, 
\(
    \M
        {}={}
    (\gamma+\gamma^{-1})U \tp U
\), where \(U\) is an orthonormal basis for the range of \(\M\) as defined in \cref{ass:PPPA:2}. In particular, we have 
\ifspringer
    \(
        U
            {}={}
        \nicefrac{1}{\sqrt{\gamma^2+1}}
        \begin{bsmallmatrix}
            \I\\ -\gamma \I
        \end{bsmallmatrix}
    \)
    and
    \(
        \projQ 
            {}={}
        U \tp U 
            {}={}
        \nicefrac{1}{\gamma^2 + 1}\begin{bsmallmatrix}\I & -\gamma \I\\ -\gamma \I & \gamma^2 \I\end{bsmallmatrix}.
    \)
\else
    \begin{align*}
        U
            {}={}
        \tfrac{1}{\sqrt{\gamma^2+1}}
        \begin{bmatrix}
            \I\\ -\gamma \I
        \end{bmatrix}
        \quad\text{and}
        \quad 
        \projQ 
            {}={}
        U \tp U 
        {}={}
        \tfrac{1}{\gamma^2 + 1}\begin{bmatrix}\I & -\gamma \I\\ -\gamma \I & \gamma^2 \I\end{bmatrix}.
    \end{align*}
\fi
Therefore, substituting
\(
    \DRSRho
        {}\coloneqq{}
    \blkdiag(\DRSrhocom \I_n, \DRSrhomon \I_n)
\)
into \eqref{eq:PPPA:eigcond}, it follows that 
\(
    \etamin
        {}={}
    1+\tfrac1\gamma\DRSrhocom +\gamma\DRSrhomon.
\)
\iftables
    Finally, it follows from \Cref{it:lem:quadratic:solutions} that \(\gamma>0\) satisfies $\alpha > 0$ iff it complies with \Cref{tab:gamma}.
\else
    Owing to \Cref{it:lem:quadratic:solutions}, \(\etamin >0\) iff $\gamma$ satisfies \eqref{eq:thm:pppaDRS:gamma}, and therefore \Cref{ass:PPPA:2} holds.
\fi

Having established that the underlying assumptions of \Cref{thm:pppa} are satisfied, it remains to translate the sequences therein in terms of the variables of \ref{eq:DRS}. 
\ifspringer
    First, leveraging the equivalences established in \Cref{lem:equivalence:PPPA-DRS}, it follows that
\else
        Leveraging
    the equivalences established in \Cref{lem:equivalence:PPPA-DRS}, it follows that
\fi
\begin{align*}
    \M(z^k-\bar{z}^k)
        {}={}&
    \big(-\tfrac1\gamma(u^k-v^k), u^k-v^k\big)\numberthis\label{eq:DRS:equivalence:prod}\\
    \qindef{z^k-\bar{z}^k}{\M}
        {}={}&
    \tfrac1\gamma\|u^k-v^k\|^2\numberthis\label{eq:DRSequivalence:seminorm}\\
    \qindef{\M(z^k-\bar{z}^k)}{\DRSRho}
        {}={}&
    (\tfrac{1}{\gamma^2}\DRSrhocom + \DRSrhomon)
    \|u^k-v^k\|^2.\numberthis\label{eq:DRS:equivalence:seminormoblique}
\end{align*}
Analogous to the proof of \Cref{thm:pppa}, consider the shadow sequences $(w^k)_{k\in\N} \coloneqq (\projQ x^k)_{k\in\N}$ and $(\bar{w}^k)_{k\in\N} \coloneqq (\projQ\bar{x}^k)_{k\in\N}$ with update rule \Cref{eq:PPPA-proof-w-update}.
    Owing to \cref{half:xbar} and \eqref{eq:DRS:equivalence:prod}, it holds that $w^k \in \mathcal{D}_{w^k, \bar w^k}$ if and only if $u^k = v^k$
    Consequently, if $w^k \in \mathcal{D}_{w^k, \bar w^k}$ for some $k \in \N$ then the algorithm has reached a point for which $u^k = v^k$ in a finite number of iterations, which by update rule \eqref{eq:DRS} immediately implies that $\bigl(u^k,\nicefrac1{\gamma}(u^k - s^k)\bigr) \in \zer \Tpd$.
    We will consider the case when $w^k \notin \mathcal{D}_{w^k, \bar w^k}$ in the remainder of the proof, in which case
$\alpha_k$ as defined in~\eqref{eq:pppa:alpha_k} coincides with $\alpha > 0$ as defined in this theorem.
\begin{proofitemize}
    \item \ref{it:DRS:v}:
    By \cref{it:pppa:v} and \eqref{eq:DRS:equivalence:prod}, the sequence \(\seq{u^k-v^k}\) converges to zero and 
    \ifspringer
        it holds that
        \(
            \tfrac1\gamma\big(u^k-v^k\big) \in A\bar x^k + \bar y^k
        \)
        and
        \(
            -(u^k-v^k) \in B^{-1} \bar y^k - \bar x^k, 
        \)
    \else
        \[
            \tfrac1\gamma\big(u^k-v^k\big) \in A\bar x^k + \bar y^k,
                \quad
            -(u^k-v^k) \in B^{-1} \bar y^k - \bar x^k, 
        \]
    \fi
    implying the claimed inclusion. 

    \item \ref{it:DRS:subseq}: 
        By \Cref{lem:equivalence:PPPA-DRS}, it holds that $(\bar{x}^k,\bar{y}^k)_{k\in\N} = \seq{u^k,\nicefrac1{\gamma}(2u^k - s^k - v^k)}$. Consequently, it follows from \cref{it:pppa:subseq} that every limit point (if any) of $\seq{u^k,\nicefrac1{\gamma}(2u^k - s^k - v^k)}$ belongs to $\zer \Tpd$. The claimed result then follows from the fact that \(\seq{u^k-v^k}\) converges to zero.

        \item \ref{it:DRS:bounded}:
        The claim follows directly from \cref{it:pppa:bounded} after noting that 
        \(
            \projQ z^k=\frac{1}{1+\gamma^2} \big(s^k, -\gamma s^k\big).
        \)

    \item \ref{it:DRS:rate}: The claimed rate follows from \cref{it:pppa:rate} since
    \(
        \|\M(z^k - \bar z^k)\|^2 = (1+ \tfrac{1}{\gamma^{2}})\|u^k - v^k\|^2,
    \)
    \(
        \qindef{z^0-\bar{z}^\star}{\M} = \tfrac1\gamma \nrm{s^0 - s^\star}^2
    \)
    and \(\nrm{\M} = \gamma + \tfrac1\gamma\).

    \item \ref{it:DRS:full}:
    The assertion for \(\seq{s^k}\) is of immediate verification based on the assertion for \(\seq{\projQ z^k}\) in \cref{it:pppa:full}. 
    Furthermore, by \Cref{lem:equivalence:PPPA-DRS}, we know that $\bar z^k \in (\M + \Tpd)^{-1}\M z^k$ iff \(\bar{x}^k \in J_{\gamma A}(x^k - \gamma y^k)\) and \(\bar y^k \in \nicefrac{1}{\gamma}(2J_{\gamma A}(x^k - \gamma y^k) - x^k + \gamma y^k - J_{\gamma B}(2J_{\gamma A}(x^k - \gamma y^k) - x^k + \gamma y^k))\).
        Consequently, if $J_{\gamma A}$ and $J_{\gamma B}$ are (single-valued) continuous, then \((\M + \Tpd)^{-1}\M\) is (single-valued) continuous as well.
        Therefore, the claim for 
        $\seq{u^k,\nicefrac1{\gamma}(u^k - s^k)}$
        follows from the assertion for $(\bar{x}^k,\bar{y}^k)_{k\in\N} = \seq{u^k,\nicefrac1{\gamma}(2u^k - s^k - v^k)}$ in \cref{it:pppa:full} and the fact that \(\seq{u^k-v^k}\) converges to zero.
    %
    \qedhere
\end{proofitemize}

    \end{proof}
\end{theorem}

Note that in \Cref{it:DRS:full} we show that under suitable conditions 
$\seq{u^k,\nicefrac1{\gamma}(u^k - s^k)}$
converges to some element of $\zer \Tpd$. By \cite[Prop. 6.9.2]{auslender2006Asymptotic}, this immediately implies that $\seq{u^k}$ converges to some element of $\zer \Tp$ and 
$\seq{\nicefrac1{\gamma}(u^k - s^k)}$
converges to some element of $\zer \Td$. An analogous argument holds for the limit point result of \Cref{it:DRS:subseq}.

            Under an additional $\DRSRho$-comonotonicity assumption on the primal-dual operator, last-iterate convergence rates can be obtained for the sequence $\seq{\|u^k - v^k\|^2}$.
This result follows directly from \Cref{thm:pppa:lastiter} by casting \ref{eq:DRS} as an instance of \ref{eq:PPPA-intro}, analogous to the proof of \Cref{thm:pppaDRS}.
\begin{corollary}[last-iterate convergence]\label{cor:DRS:lastiter}
    Suppose that the primal-dual operator $\Tpd$ is maximally
    $
        \DRSRho
            \coloneqq
        \blkdiag(\DRSrhocom \I_n, \DRSrhomon \I_n)
    $\hyp{}comonotone, that 
    \(
        [\DRSrhomon]_{-}[\DRSrhocom]_{-} < \tfrac{1}{4}
    \),
    that $\zer \Tpd$ is nonempty and that $\gamma$ is selected according to \eqref{eq:thm:pppaDRS:gamma}.
    Consider a sequence \(\seq{u^k, v^k, s^k}\) generated by \ref{eq:DRS} starting from $s^0 \in \R^n$ with stepsize $\gamma$
    and relaxation parameters $\lambda_k \in \left(0,2\alpha\right)$ such that $\lambda_k(2\alpha - \lambda_k)\geq\kappa>0$ uniformly for all $k$, where
    \(
        \alpha \coloneqq 1+\tfrac{1}{\gamma}\DRSrhocom+\gamma\DRSrhomon > 0
    \).
    Then, for all $s^\star \in \{(x^\star-\gamma y^\star)\mid (x^\star,y^\star)\in\zer \Tpd\}$, the following convergence estimates hold:
    \begin{align*}
        \|u^N-v^N\|^2 
            {}\leq{}&
        \frac{\|s^0-s^\star\|^2}{(N+1)\kappa}
        \quad\text{and}\quad
        \|u^N-v^N\|^2
            {}={}
        o\bigl(\nicefrac{1}{(N+1)}\bigr).
    \end{align*}
\end{corollary}

Finally, consider the class of piecewise polyhedral mappings, which trivially includes piecewise linear mappings and the subgradient of all piecewise linear quadratic convex functions such as normal cones of polyhedral sets \cite[\S 12]{rockafellar2009Variational}.

\begin{definition}[piecewise polyhedral mappings]
    A set-valued mapping \(T:\R^n\rightrightarrows\R^n\) is said to be piecewise polyhedral if its graph can be expressed as the union of finitely many polyhedral sets. 
\end{definition}

The following theorem shows that \ref{eq:DRS} achieves (local) linear convergence for (possibly) nonmonotone piecewise polyhedral mappings. The proof is similar to \cite[Lem. IV.4]{Latafat2019Randomized} and relies on the fact that piecewise polyhedral mappings are metrically subregular at all points in their graph \cite[\S 3]{dontchev2009implicit} and on previous results from \Cref{thm:PPPA:linear,lem:equivalence:PPPA-DRS}.
\ifspringer
    It is included in \Cref{proof:cor:DRS:CP:lin} for completeness.
\fi

\begin{theorem}[linear convergence of \ref{eq:DRS} for piecewise polyhedral mappings]
    \label{cor:DRS:CP:lin}
    Suppose that \(A\) and \(B\) are piecewise polyhedral mappings, that \(\zer \Tpd\) is nonempty closed convex, and that \cref{ass:SWMVIDRS} holds with \(\pazocal{S}^\star = \zer \Tpd\). Set the stepsize \(\gamma\) according to \eqref{eq:thm:pppaDRS:gamma} and suppose that the resolvents \(J_{\gamma A}\), \(J_{\gamma B}\) have full domain and are continuous. 
    Consider a sequence \(\seq{u^k, v^k, s^k}\) generated by \ref{eq:DRS} starting from $s^0 \in \R^n$ with stepsize $\gamma$
    and relaxation parameters $\lambda_k \in \left(0,2\alpha\right)$ such that \(\liminf_{k\rightarrow \infty} \lambda_k(2\alpha -\lambda_k) > 0\), where $\alpha$ is defined as in
    \Cref{thm:pppaDRS}.
    Then, \(\seq{s^k}\)
    converges \(R\)-linearly to some element of \(\{(x^\star-\gamma y^\star)\mid (x^\star,y^\star)\in\zer \Tpd\}\)
    and \(\seq{\|u^k-v^k\|}\) converges
    \(R\)-linearly to zero.
    \ifspringer\else
        \begin{proof}
            See \Cref{proof:cor:DRS:CP:lin}.
        \end{proof}
    \fi
\end{theorem}

  \section{Semimonotone operators}\label{sec:semi}

In this section, we introduce the class of semimonotone operators
and provide several examples of well-known function classes with semimonotone subdifferentials.
Sufficient conditions for the convergence of \ref{eq:DRS} applied to the sum of two semimonotone operators will be provided in \Cref{sec:drs:semi}.

We begin by providing the definition of \((\mon,\com)\)-semimonotone operators, which can be seen as a natural extension of (hypo)monotone and co(hypo)monotone operators.
The semimonotone terminology adapted here is inspired by the notion of semimonotonicity introduced in \cite[Def. 2]{otero2011regularity} (which corresponds to the case $\mon = \com$).


\begin{definition}[semimonotonicity]\label{def:semimonotonicity}
    Let $\mon,\com\in\R$.
    An operator 
    $A : \R^n \rightrightarrows \R^n$ is said to be $(\mon,\com)$\hyp{}semimonotone at $(\other{x}, \other{y}) \in \graph A$ if
    \begin{equation}\label{eq:obliquequasisemimonotonicity:matrix}
        \inner*{x - \other{x}, y - \other{y}} \geq \mon\nrm{x - \other{x}}^2 + \com\nrm{y - \other{y}}^2, \qquad \text{for all $(x, y) \in \graph A$}.
    \end{equation}
    An operator \(A\) is said to be $(\mon,\com)$\hyp{}semimonotone if it is $(\mon,\com)$\hyp{}semimonotone at all $(\other{x}, \other{y})\in \graph A$.
    It is said to be maximally \((\mon,\com)\)-semimonotone if its graph is not strictly contained in the graph of another \((\mon,\com)\)-semimonotone operator. 
\end{definition}
Due to the additional degree of freedom this operator class enjoys, it is able to cover a wide array of other commonly used operator classes available in literature. We list the most prominent connections below.
\begin{remark}[relationship with other types of operators]\label{rem:relationship:WMVI}
    \begin{enumerate}
        \item $(0,0)$\hyp{}semimonotonicity is equivalent to the monotonicity of $A$ (and of its inverse $A^{-1}$).
        \item $(\mon,0)$\hyp{}semimonotonicity is equivalent to $\mon$\hyp{}monotonicity, which is also known as $|\mon|$\hyp{}\emph{hypo}monotonicity \cite[Example 12.28]{rockafellar2009Variational} when $\mon < 0$ and as \emph{strong} monotonicity when $\mon > 0$.
        \item Similarly, $(0,\com)$\hyp{}semimonotonicity is equivalent to $\com$\hyp{}comonotonicity, which is also referred to as $|\com|$\hyp{}co\emph{hypo}monotonicity when $\com < 0$ and as $\com$\hyp{}cocoercivity when $\com > 0$. 
        \item $(0,\rho)$\hyp{}semimonotonicity of an operator $T$ at a point $(x^\star, 0) \in \graph T$ is equivalent to $T$ having $\rho$\hyp{}weak Minty solutions at \(x^\star\in\zer T\), sometimes also referred to as $\com$\hyp{}star\hyp{}cocoercivity around $x^\star$ \cite[Def. 2.1]{gorbunov2022Extragradient}.
        \item When $\mon = \com < 0$, our definition of $(\mon,\com)$\hyp{}semimonotonicity is equivalent to the definition of $|\com|$\hyp{}semimonotonicity proposed in \cite[Def. 2]{otero2011regularity}.
        \item $\alpha$\hyp{}averagedness of an operator is equivalent to $\left(1-\nicefrac{1}{(2-2\alpha)},\nicefrac{1}{(2-2\alpha)}\right)$\hyp{}semimonotonicity \cite[Prop. 4.35(iv)]{bauschke2017Convex}. Additionally, an operator is firmly nonexpansive if and only if it is $\nicefrac{1}{2}$\hyp{}averaged \cite[Remark 4.34(iii)]{bauschke2017Convex}, which holds for any $\left(0,\com\right)$\hyp{}semimonotone operator with $\com \geq 1$.
        \item When $\mon > 0$ and $\com < 0$, then the definition of $(\mon, \com)$\hyp{}semimonotonicity coincides with the definition of $(-\com, \mon)$\hyp{}relaxed-cocoercivity \cite[Def. 2.3]{huang2007explicit}.
        \qedhere
    \end{enumerate}
\end{remark}
\Cref{fig:semi:srg} provides the Scaled Relative Graph (SRG) \cite{ryu2021Scaled} of different operator classes and a graphical representation of the connections between them.
In \Cref{fig:semi:srg:a}, the green (resp. red) hatched area represents the region where all operators (resp. no operators) satisfy $(\mon,\com)$\hyp{}semimonotonicity, which is a consequence of the following proposition.

\begin{figure}
    \centering
    \ifspringer
        \begin{subfigure}[b]{0.5\textwidth}
            \centering
            \scalebox{0.875}{
            \includetikz{Semimonotonicity/operator_overview}
            }
            \caption{}
            \label{fig:semi:srg:a}
        \end{subfigure}\hspace{0cm}%
        \begin{subfigure}[b]{0.48\textwidth}
            \footnotesize
            \centering
            \begin{tabular}[b]{@{}c@{\;\;}c@{\;\;}c@{}}
            Monotone & Hypomonotone & Strongly monotone\\
            \includetikz{Semimonotonicity/monotone}&
            \includetikz{Semimonotonicity/hypomonotone}&
            \includetikz{Semimonotonicity/strmonotone}\\
            Cocoercive & FNE & $\alpha$-averaged\\
            \includetikz{Semimonotonicity/cocoercive}&
            \includetikz{Semimonotonicity/FNE}&
            \includetikz{Semimonotonicity/averaged}\\
            Cohypomonotone & Relaxed cocoercive & Semimonotone\\
            \includetikz{Semimonotonicity/cohypomonotone}&
            \includetikz{Semimonotonicity/relaxed}&
            \includetikz{Semimonotonicity/semimonotone}\\
            [-3pt]
            \end{tabular}
            \vspace{1pt}
            \caption{}
            \label{fig:semi:srg:b}
        \end{subfigure}
    \else
        \begin{subfigure}[b]{0.5\textwidth}
            \centering
            \includetikz{Semimonotonicity/operator_overview}
            \caption{}
            \label{fig:semi:srg:a}
        \end{subfigure}\hspace{-1.25cm}%
        \begin{subfigure}[b]{0.58\textwidth}
            \scriptsize
            \centering
            \begin{tabular}[b]{@{}c@{\;\;}c@{\;\;}c@{}}
            Monotone & Hypomonotone & Strongly monotone\\
            \includetikz{Semimonotonicity/monotone}&
            \includetikz{Semimonotonicity/hypomonotone}&
            \includetikz{Semimonotonicity/strmonotone}\\
            Cocoercive & FNE & $\alpha$-averaged\\
            \includetikz{Semimonotonicity/cocoercive}&
            \includetikz{Semimonotonicity/FNE}&
            \includetikz{Semimonotonicity/averaged}\\
            Cohypomonotone & Relaxed cocoercive & Semimonotone\\
            \includetikz{Semimonotonicity/cohypomonotone}&
            \includetikz{Semimonotonicity/relaxed}&
            \includetikz{Semimonotonicity/semimonotone}\\
            [-2pt]
            \end{tabular}
            \caption{}
            \label{fig:semi:srg:b}
        \end{subfigure}\hspace{-1.25cm}
    \fi
    \caption{
        (left)
        Relationship between semimonotonicity and other types of operators. The green hatched area represents the region where all operators satisfy $(\mon,\com)$\hyp{}semimonotonicity and the red hatched area represents the region for which there do not exist any $(\mon,\com)$\hyp{}semimonotone operators (see \Cref{prop:semi:young}).
        (right)
        Scaled Relative Graph (SRG) for different operator classes. For instance, in the top row, the SRG of a $\mon$\hyp{}monotone operator is visualized
        \cite[Prop. 3.3 and Thm. 3.5]{ryu2021Scaled}.
        In the bottom right figure a \((\mon,\com)\)\hyp{}semimonotone operator with negative \(\mon\) and \(\com\) is visualized.
        For more details on the SRG and similar geometric illustrations, we refer the interested reader to \cite{eckstein1989Splitting,giselsson2016Linear,ryu2021Scaled}.
        }
    \label{fig:semi:srg}
\end{figure}

\begin{proposition}[existence of semimonotone operators]\label{prop:semi:young}
    Let $\mon,\com\in\R$. The following hold.
    \begin{enumerate}
        \item\label{it:prop:semi:young:1}
        If $[\mon]_{-}[\com]_{-} \geq \tfrac{1}{4}$,
        then  all operators $A : \R^n \rightrightarrows \R^n$ 
        satisfy the definition of \((\mon, \com)\)\hyp{}semimonotonicity.
        \item\label{it:prop:semi:young:2} 
        If $[\mon]_+[\com]_+ > \tfrac{1}{4}$, 
        then there does not exist any operator $A : \R^n \rightrightarrows \R^n$ that is \((\mon, \com)\)\hyp{}semimonotone.
        \item\label{it:prop:semi:young:3} 
        If $[\mon]_+[\com]_+ = \tfrac{1}{4}$,
        then operator $A : \R^n \rightrightarrows \R^n$
        \begin{enumerate}
            \item\label{it:prop:semi:young:3:regular}
            is \((\mon, \com)\)\hyp{}semimonotone
            iff for some $c \in \R^n$ it holds that $\graph A \subseteq \graph T$, where $T : x \mapsto 2\mon x + c$.
            \item\label{it:prop:semi:young:3:maximal} 
            is maximally \((\mon, \com)\)\hyp{}semimonotone
            iff for some $c \in \R^n$ it holds that $A : x \mapsto 2\mon x + c$.
        \end{enumerate}
    \end{enumerate}
    \begin{proof}
        See \Cref{proof:prop:semi:young}.
    \end{proof}
\end{proposition}


We proceed to examine basic properties and calculus rules for the class of (maximally) semimonotone operators.
For instance, 
their inverses belong to the same class of operators, with the roles of \(\mon\) and \(\com\) reversed.
Additionally, the following proposition analyzes the scaling and shifting of semimonotone operators.
\begin{proposition}[inverting, shifting and scaling]\label{prop:calculus}
        Let operator $A : \R^n \rightrightarrows \R^n$ be (maximally) $(\mon,\com)$\hyp{}semimonotone \optional{at $(\other{x}, \other{y}) \in \gph A$}.
        Then, the following hold.
        \begin{enumerate}
            \item\label{prop:calculus:inverse} The inverse operator $A^{-1}$ is (maximally) $(\com, \mon)$\hyp{}semimonotone [at $(\other{y}, \other{x})\in \graph A^{-1}$].
            \item\label{prop:calculus:scaling}
            For all $\alpha \in \R_{++}$ and $u,w \in \R^n$, operator $T(x) \coloneqq w + \alpha \A (x+u)$ is (maximally) $(\alpha\mon,\alpha^{-1}\com)$\hyp{}semimonotone \optional{at $(\other{x}-u,w + \alpha\other{y})$}.
        \end{enumerate}
    \end{proposition}
Linear operators are a particular instance of semimonotone operators. For these operators, there exists an interplay between $\mon$ and $\com$, which is summarized in the following proposition. This result generalizes \cite[Prop. 5.1]{bauschke2021Generalized} for $\mon$\hyp{}monotone and $\com$\hyp{}comonotone operators. Additionally, a subsequent corollary for multiple of identity operators is provided, including its maximality. 
\begin{proposition}[linear operator]\label{it:prop:semi:lin:1}
    Let $A \in \R^{n \times n}$ and let \(\mon,\com\in\R\).
    Then, $A$ is \((\mon, \com)\)\hyp{}semimonotone if and only if \(\tfrac{1}{2}(A + A^\top) - \com A^\top A \succeq \mon \I\).
    \begin{proof}
        Consider $x \in \R^n$. Then, it holds by \((\mon, \com)\)\hyp{}semimonotonicity that $\inner*{x, Ax} \geq \mon \nrm{x}^2 + \com \nrm{Ax}^2 \Leftrightarrow \inner{x, (A - \mon \I - \com A^\top A)x} \geq 0 \Leftrightarrow \inner*{x, \left(\tfrac{1}{2}(A + A^\top) - \mon \I - \com A^\top A\right)x} \geq 0
        \Leftrightarrow \tfrac{1}{2}(A + A^\top) - \mon \I - \com A^\top A \succeq 0$.
    \end{proof}
\end{proposition}

\begin{corollary}[multiple of identity]\label{it:prop:semi:lin:2}
    If $A = \alpha \I_n$ for some $\alpha \in \R$, then $A$ is \(\left(\alpha(1-c\alpha), c\right)\)\hyp{}semimonotone for any $c \in \R$. Furthermore, $A$ is maximal except when $\alpha < 0$ and $c = \nicefrac{1}{2\alpha}$.
    \begin{proof}
        See \Cref{proof:it:prop:semi:lin:2}.
    \end{proof}
\end{corollary}

The upcoming proposition discusses the sum and the parallel sum of two semimonotone operators, $A$ and $B$. Notice that the semimonotonicity parameters of the resulting operator always involve both a sum and a parallel sum of $A$ and $B$'s semimonotonicity parameters. As expected, the commutative property holds as swapping $A$ and $B$ does not affect the result.

\begin{proposition}[sum and parallel sum]\label{lem:sum}
    Let operator $A : \R^n \rightrightarrows \R^n$ be $(\monA,\comA)$\hyp{}semimonotone \optional{at $(\other{x}_A, \other{y}_A) \in \graph A$} and operator $B : \R^n \rightrightarrows \R^n$ be $(\monB,\comB)$\hyp{}semimonotone \optional{at $(\other{x}_B, \other{y}_B) \in \graph B$}. 
    Then, it holds that
    \begin{enumerate}
        \item \label{lem:sum:1} \optional{if $\other{x} = \other{x}_A = \other{x}_B$, then} $A+B$ is $(\monA+\monB,\comA \Box \comB)$\hyp{}semimonotone \optional{at $(\other{x}, \other{y}_A+\other{y}_B)$} if either $\comA+\comB > 0$ or $\comA=\comB=0$.
        \item \label{lem:sum:2} \optional{if $\other{y} = \other{y}_A = \other{y}_B$, then} $A \Box B$ is $(\monA \Box \monB, \comA+\comB)$\hyp{}semimonotone \optional{at $(\other{x}_A+\other{x}_B, \other{y})$} if either $\monA+\monB > 0$ or $\monA=\monB=0$.
    \end{enumerate}
    \begin{proof}
        See \Cref{proof:lem:sum}.
    \end{proof}
\end{proposition}

Combining the previous proposition for the sum with the characterization of multiple of identity operators from \Cref{it:prop:semi:lin:2}, 
the semimonotonicity of the sum of a semimonotone and a multiple of identity operator
is obtained in
\Cref{it:cor:sum:identity:2}. However, under a particular condition, a 
stronger equivalence relationship can be obtained by exploiting the multiple of identity structure present. This result is provided in \Cref{it:cor:sum:identity:1} and generalizes \cite[Lem. 2.6 \& 2.8]{bauschke2021Generalized}.

\begin{proposition}[sum with identity]\label{cor:sum:identity}
    Let $\monA, \comA \in \R$. Consider the operator $A : \R^n \rightrightarrows \R^n$ and define $T \coloneqq A + \alpha \I$ where $\alpha \in \R$.
    Then, the following hold.
    \begin{enumerate}
        \item \label{it:cor:sum:identity:1} If $1 + 2\comA\alpha > 0$, then $T$ is (maximally) $\left(\tfrac{\monA + \alpha(1 + \comA\alpha)}{1 + 2\comA\alpha}, \tfrac{\comA}{1 + 2\comA\alpha}\right)$\hyp{}semimonotone \optional{at $(\other{x}, \other{y} + \alpha\other{x}) \in \graph T$}
        if and only if $A$ is (maximally) $(\monA,\comA)$\hyp{}semimonotone \optional{at $(\other{x}, \other{y}) \in \graph A$}.
        \item \label{it:cor:sum:identity:2} If $1 + 2\comA\alpha \leq 0$, then $T$ is $\left(\monA+\alpha(1-c\alpha),\comA \Box c\right)$\hyp{}semimonotone \optional{at $(\other{x}, \other{y}+\alpha \other{x}) \in \graph T$} for all $c > -\comA$ if $A$ is $(\monA,\comA)$\hyp{}semimonotone \optional{at $(\other{x}, \other{y}) \in \graph A$}.
    \end{enumerate}
    \begin{proof}
        See \Cref{proof:cor:sum:identity}.
    \end{proof}
\end{proposition}

Leveraging previous results, we now formalize the relationship between the semimonotonicity of an operator $T$ defined 
as a
 function of another (maximally) monotone operator $A$. In practical applications, this result can be leveraged to check (maximal) semimonotonicity of a given operator or to generate semimonotone operators from a given monotone operator.
%
\begin{corollary}[connection to monotone operators]\label{cor:sem2mon}
    Let $\xi, \nu \in \R$. Consider the operator $A : \R^n \rightrightarrows \R^n$ and define $T \coloneqq (A + \xi \id)^{-1} + \nu \id$. Then, the following hold.
    \begin{enumerate}
        \item \label{it:cor:gen:1} If $1 + 2\xi\nu > 0$, then $T$ is (maximally) $\Bigl(\tfrac{\nu(1 + \xi\nu)}{1 + 2\xi\nu}, \tfrac{\xi}{1 + 2\xi\nu}\Bigr)$\hyp{}semimonotone if and only if $A$ is (maximally) monotone.
        \item \label{it:cor:gen:2} 
        If $1 + 2\xi\nu \leq 0$, then
        $T$ is $\Bigl(\nu(1-c \nu), \xi \Box c\Bigr)$\hyp{}semimonotone for all $c > -\xi$ if $A$ is monotone.
    \end{enumerate}

    \begin{proof}
        By applying \Cref{cor:sum:identity,prop:calculus:inverse} it follows that $A$ is (maximally) monotone if and only if $A + \xi \id$ is (maximally) \((\xi, 0)\)\hyp{semimonotone} if and only if $(A + \xi \id)^{-1}$ is (maximally) \((0, \xi)\)\hyp{semimonotone}. 
        The claimed result then follows by using \Cref{cor:sum:identity} again.
    \end{proof}
\end{corollary}
When \(\nu =0\), \cref{it:cor:gen:1} reduces to a known equivalence between comonotonicity of an operator and monotonicity of its Yosida regularization, see \cite[Prop. 6.9.3]{burachik2008setvalued}.
The following 
result
is a direct consequence of \Cref{it:cor:gen:1}, noting that for
\ifspringer
    this particular choice of $\xi$ and $\nu$
\else
    the choice of $\xi$ and $\nu$ from \eqref{eq:XiBeta}
\fi
it holds that
\(
    1+2\xi\nu = \nicefrac{1}{\sqrt{1-4\com\mon}} > 0.        
\)
This result generalizes \cite[Thm. 2]{otero2011regularity}, which only covers the case \(\mon=\com\). 
\begin{corollary}\label{cor:generation semimonotone operator}
    A set-valued operator $T: \R^n \rightrightarrows \R^n$ is (maximally) $(\mon, \com)$\hyp{}semimonotone with $\com\mon < \tfrac{1}{4}$ if and only if $A \coloneqq (T - \nu \id)^{-1} - \xi \id$ is a (maximally) monotone operator where
    \ifspringer
    \(
        \xi 
            {}={}
        \tfrac{\com}{\sqrt{1-4\mon\com}}  
    \)
    and
    \(
        \nu 
            {}={}
        \tfrac{2\mon}{1+\sqrt{1-4\mon\com}}
    \).
    \else
        \begin{equation}\label{eq:XiBeta}
            \xi 
                {}={}
            \tfrac{\com}{\sqrt{1-4\mon\com}}
            \quad \text{and} \quad
            \nu 
                {}={}
            \tfrac{2\mon}{1+\sqrt{1-4\mon\com}}.
        \end{equation}
    \fi
\end{corollary}

In light of the above connection with maximally monotone operators, a $(\mon, \com)$\hyp{}semimonotone operator \(A\) with $\com\mon < \tfrac{1}{4}$ is maximal if \((A-\nu \id)^{-1}\) is continuous (since any continuous monotone mapping is maximally monotone). 
Additionally, the following proposition shows that the graph of a maximally $(\mon, \com)$\hyp{}semimonotone operator is closed.
    \begin{proposition}[outer semicontinuity of a semimonotone operator]\label{prop:semimon:osc}
        Let $A: \R^n \rightrightarrows \R^n$ be a maximally $(\mon, \com)$\hyp{}semimonotone operator with either $\com\mon < \tfrac{1}{4}$ or \([\mon]_+[\com]_+ = \tfrac14\). Then, $A$ is outer semicontinuous.
        \begin{proof}
            See \Cref{proof:prop:semimon:osc}.
        \end{proof}
    \end{proposition}
We conclude this section by showing that if the stepsize is selected appropriately, the resolvent of a semimonotone operator has full domain and is Lipschitz continuous.
This result generalizes \cite[Prop. 2.17]{bauschke2021Generalized} for comonotone operators and \cite[Thm. 3]{otero2011regularity}, which 
is restricted to
 the case $\mon=\com$.
\begin{proposition}[resolvent of a semimonotone operator]\label{cor:semimon:fulldom}
    Let $A: \R^n \rightrightarrows \R^n$ be a $(\mon, \com)$\hyp{}semimonotone operator with either $\com\mon < \tfrac{1}{4}$ or \([\mon]_+[\com]_+ = \tfrac14\).
    \iftables
        Choose $\gamma$ according to \Cref{tab:semimon:domain}.
    \else
        Let 
        \begin{align*}
            \gamma
                {}\in{}
            \left(
                \tfrac{2[-\com]_{+}}{1+\sqrt{1-4\com\mon}},\tfrac{1+\sqrt{1-4\com\mon}}{2[-\mon]_{+}}
            \right).\numberthis\label{eq:cor:fulldom:gamma}
        \end{align*}
    \fi
    Then, the following hold.
    \begin{enumerate}
        \item \label{cor:semimon:resolvent} The resolvent \(J_{\gamma A}\) has full domain if and only if \(A\) is maximal.  
        \item \label{cor:semimon:Lip} The resolvent \(J_{\gamma A}\) is Lipschitz continuous with constant
        \(
            \smash{
                L(\gamma) = \tfrac{|1 + 2\com\gamma^{-1}| + \sqrt{1-4\mon\com}}{2(1 + \mon\gamma + \com \gamma^{-1})}.
            }
        \)
    \end{enumerate}
    \begin{proof}
        See \Cref{proof:cor:semimon:fulldom}.
    \end{proof}
\end{proposition}

\iftables
    \begin{table}
        \centering
        \begin{tabular}{c}
            {
            \begin{tblr}{
                hlines,
                vlines,
                vline{1-1} = {1-1}{0pt},
                hline{1-1} = {1-1}{0pt},
                colspec={M{2cm}M{4cm}M{4cm}},
                cell{1-1}{2-3}={c,gray!20},
                cell{2-3}{1-1}={c,gray!20}
                }
                &$\mon \geq 0$ & $\mon < 0$\\
                $\com \geq 0$&$\gamma \in \left(0,+\infty\right)$&$\gamma \in \left(0, -\tfrac{1+\sqrt{1-4\com\mon}}{2\mon}\right)$\\
                $\com < 0$&$\gamma \in \left(\tfrac{-2\com}{1+\sqrt{1-4\com\mon}},+\infty\right)$&$\gamma \in \left(\tfrac{-2\com}{1+\sqrt{1-4\com\mon}},-\tfrac{1+\sqrt{1-4\com\mon}}{2\mon}\right)$
            \end{tblr}
            }
        \end{tabular}
        \caption{Range of $\gamma$ in \Cref{cor:semimon:fulldom}.}
        \label{tab:semimon:domain}
    \end{table}
\fi

      \subsection{Examples of functions with semimonotone subdifferentials}\label{subsec:semi:subdif}

In what follows, we present various function classes of which the (limiting) subdifferentials belong to the class of semimonotone operators.
\ifspringer
    Proofs of the upcoming results can be found in \Cref{subsec:auxiliary:functions}.

    We start by studying the class of functions with upper and/or lower curvature.
\else
    We start by studying the class of functions with upper and/or lower curvature, which are defined as follows.
\fi
\begin{definition}[upper and lower curvature]
    Let \(\sigma, \ell \in \R\) and consider a proper lsc function \(f:\R^n \to \Rinf\). Then, 
    \begin{enumerate}
        \item (lower curvature) \(f\) is said to have lower curvature with modulus \(\sigma\) (denoted by \(f \in \mathcal F_{\sigma}(\R^n)\)) if  
        \(
            \varphi_{\sigma} \coloneqq f - \tfrac{\sigma}2\|\cdot\|^2
        \) is convex. 
        \item (upper curvature) \(f\) is said to have upper curvature with modulus \(\ell\) 
        (denoted by \(f \in \mathcal F^{\ell}(\R^n)\))
        if  
        \(f\) is finite-valued and
        \(
            \varphi^{\ell} \coloneqq \tfrac{\ell}2\|\cdot\|^2 - f
        \) is convex. 
    \end{enumerate}
    A function with both upper and lower curvatures is denoted by \(\mathcal F_{\sigma}^{\ell}(\R^n)\). 
\end{definition}

Many traditional function classes fit within this curvature framework.
For instance, functions with positive (resp. negative) lower curvatures are widely known as strongly (resp. hypo) convex functions. Another notable example is that of $\ell$\hyp{}smooth functions, which have upper curvature $\ell$ and lower curvature $-\ell$ \cite[Lem. 2.1]{themelis2020DouglasRachford}.

The next proposition presents the semimonotonicity properties of the subdifferential of functions according to their curvature. Additionally, a corollary is provided for the pointwise minimum of functions with upper curvatures.
\begin{proposition}[functions with upper and lower curvature]\label{prop:curvature}
    Let \(\sigma, \ell \in \R\). Then, the following hold. 
    \begin{enumerate}
        \item \label{prop:hypo:cvx} If \(f\in
        \mathcal F_{\sigma}(\R^n)\), \ie, \(f\) is \(\sigma\)-convex, then, 
        $\partial f$ is $\big(\sigma,0\big)$\hyp{}semimonotone.

        \item \label{prop:hypo:ccv} If $f$ is continuous and \(f\in
        \mathcal F^{\ell}(\R^n)\) with \(\ell<0\), then, 
        $\partial f$ is $\big(0, \tfrac1{\ell}\big)$\hyp{}semimonotone.
        \item \label{prop:ccv} If $f$ is continuous and \(f\in
        \mathcal F^{\ell}(\R^n)\) with \(\ell\geq0\), then,  
        $\partial f$ is 
        \(
        (\alpha(1-c\alpha), c \Box \tfrac1{\ell}) 
        \)-semimonotone for all \(\alpha> \ell\) and \(c> \tfrac{1}{\alpha - \ell}\). 
    \end{enumerate}
    If \(f\in \mathcal F_{\sigma}^{\ell}(\R^n)\), then, \(\sigma\leq \ell\), \(f\) is continuously differentiable and  the following hold.  
    \begin{enumerate}[resume]
        \item \label{prop:hypo:dif+} If $\ell + \sigma > 0$, then $\nabla f$ is $\big(\sigma \Box \ell,\tfrac{1}{\ell+\sigma}\big)$\hyp{}semimonotone.
        \item \label{prop:hypo:dif-} 
        If $\ell + \sigma \leq 0$, %
        $\nabla f$ is $\big(\sigma(1-c\sigma), c \Box \tfrac{1}{\ell-\sigma}\big)$\hyp{}semimonotone for any constant \(c\) such that \(1 + c(\ell-\sigma) > 0 \). 
    \end{enumerate}%
    \ifspringer\else
        \begin{proof}
            See \Cref{proof:prop:curvature}.
        \end{proof}
    \fi
\end{proposition}

\begin{corollary}[pointwise minimum]\label{lem:pointwise:semi}
    Let \(f\) be the pointwise min of finite-valued and continuous functions \(f_i \in \mathcal F^{\ell_{f_i}}(\R^n)\) where \(\ell_{f_i}\in \R\), i.e., \(f(x) = \min\{f_1(x),\ldots, f_N(x)\}\). Let \(\bar \ell = \max_i{\ell_{f_i}}\). Then, the following hold.
    \begin{enumerate}
        \item\label{lem:pointwise:semi:1} If $\bar \ell < 0$, then \(\partial f\) is \(\left(0, \tfrac{1}{\bar \ell}\right)\)\hyp{}semimonotone.
        \item\label{lem:pointwise:semi:2} If $\bar \ell \geq 0$, then \(\partial f\) is     
        \(
            \left(\alpha(1-c\alpha), c \Box \tfrac1{\bar \ell}\right) 
        \)\hyp{}semimonotone
        for all $\alpha > \bar \ell$ and $c > \tfrac1{\alpha - \bar \ell}$.
    \end{enumerate}%
    \ifspringer\else
        \begin{proof}
            See \Cref{proof:lem:pointwise:semi}.
        \end{proof}
    \fi
\end{corollary}

One of the most fundamental operations for functions is the infimal convolution, which is the subject of the upcoming proposition. To obtain this result, we first show that the subdifferential of the infimal convolution of two semimonotone operators belongs to the parallel sum of the subdifferentials of these two functions (see \Cref{it:infimal:subdif}), and then we use \Cref{lem:sum:2}.

\begin{proposition}[infimal convolution]\label{lem:infimal}
    For \(i=1,2\), let \(f_i:\R^n\to \Rinf\) be proper lsc functions whose subdifferentials are \((\mon_i,\com_i)\)-semimonotone with \(\mon_i,\com_i\in\R\) such that 
    $\mon_1 + \mon_2 > 0$ or $\mon_1 = \mon_2 = 0$.
    Denote their infimal convolution by
    \ifspringer
        \(
            \varphi : s \mapsto \min_{w \in \R^n} \bigl\{f_1(w) + f_2(s-w)\bigr\}.
        \)
    \else
        \[
            \varphi : s \mapsto \min_{w \in \R^n} \bigl\{f_1(w) + f_2(s-w)\bigr\}.
        \]
    \fi
    Suppose that for each \(r_1\geq 0\) and \(\bar s\in\R^n\) there exists \(r_2>0\) such that the set 
    \begin{align*}
        \textstyle
        \left\{(w_1,w_2)\in\R^n \;\middle|\; \sum_{i=1}^2f_i(w_i) \leq r_1, \; \left\|\sum_{i=1}^2 w_i-\bar s\right\|\leq r_2\right\}\numberthis\label{eq:lem:infimal:uniform-level-bounded}
    \end{align*}
    is bounded. 
    Then, the following hold.
    \begin{enumerate}
        \item \label{it:infimal:subdif} $\partial \varphi(\bar s) \subseteq (\partial f_1 \Box \partial f_2)(\bar s)$.
        \item \label{it:infimal:semi} \(\partial \varphi\) is \((\mon_1 \Box \mon_2,\com_1+\com_2)\)-semimonotone.
    \end{enumerate}%
    \ifspringer\else
        \begin{proof}
            See \Cref{proof:lem:infimal}.
        \end{proof}
    \fi
\end{proposition}

As a final example, we show that the associated saddle operator of a minimax optimization problem is semimonotone under the so-called interaction dominance condition \cite{grimmer2022landscape}.
\newcommand{\SN}{S} 
\begin{proposition}[saddle operator]\label{prop:saddleOp}
    Suppose that
        $\varphi : \R^n\times \R^m\rightarrow \Rinf$
    is twice continuously differentiable and $\sigma$-convex-concave for some $\sigma\in\R$, i.e.,
    \(
        \nabla_{xx}^{2}\varphi(z) \succeq \sigma \I 
    \)
    and
    \(
        -\nabla_{yy}^{2}\varphi(z) \succeq \sigma \I
    \)
    for all \(z \in \R^{n+m}\).
    Moreover, the $\alpha$-interaction dominance condition \cite[Def. 1]{grimmer2022landscape} holds for some $\alpha> -\tfrac{1}\delta$ and $\delta\in(0,\nicefrac1{[-\sigma]_+})$, \ie, for all \(z \in \R^{n+m}\)
    \begin{align*}
     \nabla_{xx}^{2}\varphi(z)+\nabla_{xy}^{2}\varphi(z)\left(\tfrac1{\delta}\I-\nabla_{yy}^{2}\varphi(z)\right)^{-1}\nabla_{yx}^{2}\varphi(z)
        {}\succeq{}&
     \alpha\I,
     \\
     -\nabla_{yy}^{2}\varphi(z)+\nabla_{yx}^{2}\varphi(z)\left(\tfrac1{\delta}\I+\nabla_{xx}^{2}\varphi(z)\right)^{-1}\nabla_{xy}^{2}\varphi(z)
        {}\succeq{}&
    \alpha\I.
    \end{align*}
    Then, the following hold for the associated saddle operator $T_{\varphi}(z)\coloneqq \big(\nabla_x\varphi(z), -\nabla_y\varphi(z)\big)$. 
    \begin{enumerate}
        \item If $\alpha < \tfrac{1}{\delta}$, then \(T_{\varphi}\) is \((\tfrac{\alpha}{1-\alpha \delta},-\tfrac{\delta}{1-\alpha \delta})\)-semimonotone. 
        \item If $\alpha\geq 0$, then $T_{\varphi}$ is $(0, -\delta)$-semimonotone.
    \end{enumerate}%
    \ifspringer\else
        \begin{proof}
            See \Cref{proof:prop:saddleOp}.
        \end{proof}
    \fi
\end{proposition}

  \section{Douglas-Rachford splitting for semimonotone operators}\label{sec:drs:semi}


In practical applications, it might be difficult to determine whether the conditions of \Cref{ass:DRS} for \ref{eq:DRS} are satisfied. 
To address this issue, we provide a set of sufficient conditions for the convergence of \ref{eq:DRS} involving the sum of two semimonotone operators. More specifically, we work under the following assumptions on the underlying operators $A$ and $B$ from problem~\eqref{prob:composite}.

\begin{assumption}\label{ass:SWMVIDRSsemimonotone}
    In problem~\eqref{prob:composite}, there exist parameters $\monA, \comA, \monB, \comB \in \R$ and a nonempty set \(\pazocal{S}^\star\subseteq \zer \Tpd\) such that for every $(x^\star, y^\star) \in \pazocal{S}^\star$, the following hold.
    \begin{enumeratass}
        \item \label{ass:SWMVIDRSsemimonotone:params} It holds that either $\monA + \monB > 0$ or $\monA = \monB = 0$, that either $\comA + \comB > 0$ or $\comA = \comB = 0$ and that $[\monA \Box \monB]_-[\comA \Box \comB]_- < \tfrac{1}{4}$.
        \item \label{ass:SWMVIDRSsemimonotone:A} Operator $\A$ is $(\monA, \comA)$\hyp{}semimonotone at $(x^\star,-y^\star) \in \graph A$.
        \item \label{ass:SWMVIDRSsemimonotone:B} Operator $\B$ is
        $(\monB, \comB)$\hyp{}semimonotone at $(x^\star,y^\star) \in \graph B$.
    \end{enumeratass}
\end{assumption}

Under these assumptions,
the primal-dual operator $\Tpd$ has oblique weak Minty solutions 
as shown below in \Cref{it:prop:semi:WMVI}.
    Additionally, whenever operators \(A\) and \(B\) are maximally semimonotone as in \cref{ass:SWMVIDRSsemimonotone:maximal}, the primal-dual operator is maximally comonotone.
    The proof of \Cref{prop:semi:WMVI} is deferred to \Cref{proof:prop:semi:WMVI}.

\begin{assumption}\label{ass:SWMVIDRSsemimonotone:maximal}
        In problem~\eqref{prob:composite}, 
        $\zer(A + B)$ is nonempty and there exist parameters $\monA, \comA, \monB, \comB \in \R$ satisfying \cref{ass:SWMVIDRSsemimonotone:params} such that operator $A$ is maximally $(\monA, \comA)$\hyp{}semimonotone and operator $B$ is maximally $(\monB, \comB)$\hyp{}semimonotone.
    \end{assumption}

\begin{lemma}[primal-dual operator]\label{prop:semi:WMVI}
    Suppose that
    either
    \Cref{ass:SWMVIDRSsemimonotone}
    or \Cref{ass:SWMVIDRSsemimonotone:maximal}
    holds
        and let $\DRSRho$ be defined as in \eqref{eq:WMSDRS}, where
    \begin{align*}
        \DRSrhocom = \comA \Box \comB
        \quad
        \text{and}
        \quad
        \DRSrhomon = \monA \Box \monB.
        \numberthis\label{eq:DRSRho:semimonotonicity}
    \end{align*}
    Then, the following hold.
        \begin{enumerate}
            \item\label{it:prop:semi:WMVI} If \Cref{ass:SWMVIDRSsemimonotone} holds, then $\Tpd$ has $\DRSRho$\hyp{}oblique weak Minty solutions at \(\pazocal{S}^\star\).
            \item\label{it:prop:semi:WMVI:max} If \Cref{ass:SWMVIDRSsemimonotone:maximal} holds, then $\Tpd$ is maximally $\DRSRho$\hyp{}comonotone.
        \end{enumerate}        
\end{lemma}

We are now in a position to state the main result of this section:
sufficient conditions for the convergence of \ref{eq:DRS},
through the semimonotonicity properties of the operators \(\A\) and \(\B\). 
The proof is deferred to \Cref{proof:cor:DRS:semi}.
\begin{theorem}[convergence of \ref{eq:DRS} under semimonotonicity]\label{cor:DRS:semi}
    Suppose that 
    \begin{enumerate}
        \item \label{it:cor:DRS:semi:1} either \Cref{ass:DRS:1}, \Cref{ass:DRS:2} and \Cref{ass:SWMVIDRSsemimonotone} hold,
        \item \label{it:cor:DRS:semi:2} or 
                \Cref{ass:SWMVIDRSsemimonotone:maximal} 
        hold
        s
.
    \end{enumerate}
    Let $\gamma$ be selected according to \Cref{tab:gamma:DRS:semi}, where
    \begin{equation}\label{eq:eigcondDRS:gamma:semi}
        \begin{aligned}
            \gamma_{-} \coloneqq 
            -\tfrac{2(\comA \Box \comB)}{1+\sqrt{1-4(\monA \Box \monB)(\comA \Box \comB)}}
            \quad \text{and} \quad
            \gamma_{+} \coloneqq -\tfrac{1 + \sqrt{1 - 4(\monA \Box \monB)(\comA \Box \comB)}}{2(\monA \Box \monB)}.
        \end{aligned}
    \end{equation}
    Then, 
    \(
        \alpha \coloneqq 1 + \tfrac{1}{\gamma}(\comA \Box \comB) + \gamma(\monA \Box \monB) > 0.
    \)
    Consider a sequence \(\seq{u^k, v^k, s^k}\) generated by \ref{eq:DRS} starting from $s^0 \in \R^n$ with stepsize $\gamma$
    and relaxation parameters $\lambda_k \in \left(0,2\alpha\right)$ such that $\liminf_{k\to \infty}\lambda_k(2\alpha - \lambda_k)>0$.
    Then, all the claims of \cref{thm:pppaDRS} hold.
    Furthermore, if \Cref{ass:SWMVIDRSsemimonotone:maximal} holds and $\lambda_k(2\alpha - \lambda_k)\geq\kappa>0$ uniformly for all $k$, then the convergence estimates from \Cref{cor:DRS:lastiter} also hold.
    \begin{table}
        \centering
        \caption{Range of $\gamma$ for \ref{eq:DRS} for semimonotone operators.}
        \label{tab:gamma:DRS:semi}
        \begin{tabular}{c}
            {
            \ifspringer
                \begin{tblr}{
                    hlines,
                    vlines,
                    vline{1-2} = {1-2}{0pt},
                    hline{1-2} = {1-2}{0pt},
                    colspec={M{1.6cm}M{1.3cm}M{2.9cm}M{2.2cm}M{2.2cm}},
                    cell{1-2}{3-5}={c,gray!20},
                    cell{3-5}{1-2}={c,gray!20}
                    }
                    &&\SetCell[r=2]{c}$\monA = \monB = 0$ & \SetCell[c=2]{c}$\monA + \monB > 0$\\
                    &&& $\monA\monB \geq 0$ & $\monA\monB < 0$\\
                    \SetCell[c=2]{c}$\comA = \comB = 0$&&\SetCell[c=2, r=2]{c}$\gamma \in \left(0,+\infty\right)$&&$\gamma \in \left(0, - \tfrac{1}{\monA \Box \monB}\right)$\\
                    \SetCell[r=2]{c}$\comA + \comB > 0$&$\comA\comB \geq 0$&&&$\gamma \in \left(0, \gamma_{+}\right)$\\
                    &$\comA\comB < 0$&$\gamma \in \Bigl(- (\comA \Box \comB), +\infty\Bigr)$&$\gamma \in \left(\gamma_{-}, +\infty\right)$&$\gamma \in \left(\gamma_{-}, \gamma_{+}\right)$
                \end{tblr}
            \else
                \begin{tblr}{
                    hlines,
                    vlines,
                    vline{1-2} = {1-2}{0pt},
                    hline{1-2} = {1-2}{0pt},
                    colspec={M{1.8cm}M{1.4cm}M{3.1cm}M{2.2cm}M{2.2cm}},
                    cell{1-2}{3-5}={c,gray!20},
                    cell{3-5}{1-2}={c,gray!20}
                    }
                    &&\SetCell[r=2]{c}$\monA = \monB = 0$ & \SetCell[c=2]{c}$\monA + \monB > 0$\\
                    &&& $\monA\monB \geq 0$ & $\monA\monB < 0$\\
                    \SetCell[c=2]{c}$\comA = \comB = 0$&&\SetCell[c=2, r=2]{c}$\gamma \in \left(0,+\infty\right)$&&$\gamma \in \left(0, - \tfrac{1}{\monA \Box \monB}\right)$\\
                    \SetCell[r=2]{c}$\comA + \comB > 0$&$\comA\comB \geq 0$&&&$\gamma \in \left(0, \gamma_{+}\right)$\\
                    &$\comA\comB < 0$&$\gamma \in \Bigl(- (\comA \Box \comB), +\infty\Bigr)$&$\gamma \in \left(\gamma_{-}, +\infty\right)$&$\gamma \in \left(\gamma_{-}, \gamma_{+}\right)$
                \end{tblr}
            \fi
            }
        \end{tabular}
    \end{table}
\end{theorem}

\Cref{cor:DRS:semi} provides a general framework for the convergence of \ref{eq:DRS}, unifying and extending existing convergence results on relaxed Douglas/Peaceman Rachford splitting methods. We provide several examples.
\begin{remark}[comparison with existing theory]\label{rem:comparison:DRS}
    \begin{enumerate}
        \item         W%
        hen $A$ and $B$ are maximally $(0,0)$\hyp{}semimonotone, i.e., maximally monotone, the convergence result of classic DRS is recovered, see e.g. \cite{lions1979Splitting, eckstein1992DouglasRachford}.
                \item When \(A\) and \(B\) are respectively \(\mu_A\)- and \(\mu_B\)-monotone, we recover several existing results in literature.
            For instance, in the strongly monotone setting where
            $
                \mu \coloneqq \monA = \monB > 0
            $, 
            we recover the convergence result from \cite{monteiro2018complexity}, showing that 
                        \ref{eq:DRS} converges for $\lambda \in (0, 2+\gamma\mon)$ instead of the traditional $\lambda \in (0, 2)$.
            On the other hand, when \(\mu_A+ \mu_B>0\), the stepsize range provided in \Cref{tab:gamma:DRS:semi} matches \cite[Thm. 4.5(ii)]{dao2019Adaptive}.
            Analogously, we recover the results for the sum of a $\comA$\hyp{} and a $\comB$\hyp{}comonotone operator with $\comA + \comB > 0$ from \cite[Cor. 5.5(ii)]{bartz2022Conical}.
        \item\label{rem:comparison:DRS:optimization} In the optimization setting, global convergence of \ref{eq:DRS} has recently been shown in \cite{li2016Douglas,themelis2020DouglasRachford} for structured nonsmooth, nonconvex optimization problems where the cost is split as the sum of one Lipschitz differentiable function $f: \R^n \rightarrow \R$ and one proper lsc function $g : \R^n \rightarrow \R$, i.e., $A = \nabla f$ and $B = \partial g$. These results rely on showing descent of the iterates on the so-called Douglas-Rachford envelope \cite{patrinos2014douglas}, for which only properties of the Lipschitz differentiable function $f$ are exploited to determine the range of admissible stepsizes. 
                For instance, consider the setting where $f$ is $\ell$-smooth and $\sigma$-convex, where $\ell > 0$ and $\sigma \in [-\ell, 0)$. 
                Then,
        it is shown in \cite{themelis2020DouglasRachford} that \ref{eq:DRS} converges globally for $\lambda \in (0,2)$ and $\gamma < \min\bigl\{\tfrac{\lambda-2}{2\sigma}, \tfrac{1}{\ell}\bigr\} < -\tfrac{1}{\sigma}$
        provided that $\argmin_x f(x) + g(x)$ is nonempty.

        Note that when additional semimonotonicity properties of $g$ are given, a larger stepsize range for $\gamma$ can be obtained through \Cref{cor:DRS:semi}.
            To see this, first observe that in this setting $A$ is 
            $(\monA, \comA)$\hyp{}semimonotone for some $\monA < 0$ and $\comA \geq 0$ (see \Cref{prop:hypo:dif+}).
                Therefore, supposing
        that $\partial g$ is $(\monB, \comB)$\hyp{}semimonotone at $(x^\star, y^\star) \in \graph \partial g$ where $\monA + \monB > 0$ and $\comB \geq 0$
                it follows from \Cref{cor:DRS:semi} that
        \ref{eq:DRS} converges for $\gamma \in (0, \gamma_+)$%
                , which is a larger stepsize range since $\gamma_+ > -\tfrac{1}{\sigma}$.
    \end{enumerate}
    \vspace{-10pt}
\end{remark}

  \section{Example problems}\label{subsec:semi:examples}
In this section, we provide several example problems for \ref{eq:DRS}
that go beyond the standard monotone setting.
\ifspringer
        Details of the examples
    presented in this section can be found in the arXiv preprint \cite[Appendix C]{evens2023convergence}.
\else
        Details of the examples presented
    in this section can be found in \Cref{sec:app:examples}.
\fi

We commence by applying \ref{eq:DRS} to a simple linear inclusion problem. This setting is of particular interest for two specific reasons. First, this example demonstrates the tightness of the stepsize range provided in \Cref{thm:pppaDRS} (see \Cref{it:example5:2}). Second, 
it
 highlights that our theory is able to cover cases where neither the primal nor the dual nor the primal-dual inclusion is monotone (see \Cref{it:example5:3}).
\begin{example}[saddle point problem]\label{example:ex5}
    Consider the problem of finding a stationary point of a linear inclusion, with a particular structure that emerges naturally in saddle point problems \cite[\S 1.4.1]{facchinei2003FiniteDimensional}
    \begin{equation} \label{eq:ex5-Tp}
        0 \in \Tp x
            {}={}
        \begin{bmatrix}b & a \\ -a & b\end{bmatrix}x
            {}={}
        \underbrace{\begin{bmatrix}0 & a \\ -a & 0\end{bmatrix}}_Ax + 
        \underbrace{\begin{bmatrix}b & 0 \\ 0 & b\end{bmatrix}}_Bx,
    \end{equation}
    \begin{minipage}{0.65\linewidth}
        where $a, b \in \R \setminus \{0\}$. In particular, any solution to the inclusion problem $0 \in \Tp x$ is a minimax solution of
        \(
            f(x_1,x_2) \coloneqq ax_1x_2 + \tfrac{b}{2}(x_1^2 - x_2^2) 
        \)
        when \(b > 0\) and a
        maximin solution 
        when \(b < 0\). For this problem, the
        following assertions hold.
    \end{minipage}%
    \vspace{0.225\baselineskip}
    \begin{enumerate}
        \item \label{it:example5:1} 
        \begin{minipage}[t]{0.63\linewidth}
            By examining the spectral radius of the algorithmic operator, it can be seen that the sequence $\seq{u^k, v^k, s^k}$ generated by \ref{eq:DRS} with fixed relaxation parameter $\lambda$ converges iff $\lambda$ lies in the interval
            \begin{align*}
                \begin{aligned}
                    0 < \lambda < 
                    \bar{\lambda} \coloneqq 2\frac{(a^2\gamma + b)(1+\gamma b)}{\gamma (a^2 + b^2)},
                \end{aligned}\numberthis\label{eq:example5:lambda}
            \end{align*}
            which is nonempty iff either \(b>0\) or
            \[
                b < 0,\; a^2 \ne b^2 \text{ and }
                \gamma \in \left(\min\{-\tfrac{1}{b}, -\tfrac{b}{a^2}\}, \max\{-\tfrac{1}{b}, -\tfrac{b}{a^2}\}\right).    
            \]
            A numerical experiment verifying this convergence result is provided in \Cref{fig:saddle:convergence}.
        \end{minipage}\hfill%
        \begin{minipage}[t]{0.36\linewidth}
            \centering
            \vspace{-1.75cm}
            \includetikz{Examples/Saddle/iterates-drs-saddle}
            \captionof{figure}{
                Sequences $u^k$ generated by applying \ref{eq:DRS} on \Cref{example:ex5} for $a = 2$, $b = -1$ and stepsize $\gamma = \tfrac12$.
            }
            \label{fig:saddle:convergence}
        \end{minipage}%
        \item \label{it:example5:2} The results from \Cref{thm:pppaDRS} are tight in the sense that the entire range of relaxation parameters $\lambda$ from \eqref{eq:example5:lambda} is covered, when in \eqref{def:WMVI} the vector $v$ is restricted to $\range{\M}$ (see the remark below \Cref{def:SWMVI}).
        \item \label{it:example5:3} The 
        range of \(a\) and \(b\) from \Cref{it:example5:1} for which \ref{eq:DRS} converges covers cases where neither the primal, nor the dual, nor the primal-dual problem is monotone.
    \end{enumerate}
    Additionally, an analysis of the tightness of \Cref{cor:DRS:semi} is provided in 
    \ifspringer
        the arXiv preprint \cite[Appendix C]{evens2023convergence}.
    \else
        \Cref{detail:example:ex5}.
    \fi
    In this analysis, the semimonotonicity parameters $\monA, \monB, \comA, \comB$ are selected according to \Cref{it:prop:semi:lin:1} such that the relaxation parameter in \Cref{cor:DRS:semi} is maximized. Notice that although $A$ is obviously $(0,0)$\hyp{}semimonotone, the maximal relaxation parameter is in general attained for other values of $\monA$ and $\comA$.
\end{example}

Next, we turn our attention to a nonsmooth optimization problem. 
This example highlights the generality of our underlying assumptions, as the convergence of \ref{eq:DRS} for this problem is not covered by any previous works on \ref{eq:DRS} (to the best of the authors' knowledge).
\ifspringer
    \def\widthcol{0.67\textwidth}
    \def\widthfig{0.315\textwidth}
\else
    \def\widthcol{0.59\textwidth}
    \def\widthfig{0.4\textwidth}
\fi
\begin{example}[nonsmooth optimization]\label{ex:nonsmooth:minimization}
    Consider the separable nonsmooth optimization problem\\
    \begin{minipage}[t]{\widthcol}
        \phantom{a}
        \begin{equation*}
            \begin{aligned}
                &\minimize_{x
                \in \R
} && f_A(x) + f_B(x)
            \end{aligned}
        \end{equation*}
        which is visualized in \Cref{fig:nonsmooth:minimization}, where
        \begin{align*}
            f_A(x) 
                {}\coloneqq{}&
            \begin{cases}
                3x^2 - 3x + 8, & \quad\text{if } x \in [-1, 1], \\
                -\tfrac12 x^2 - 3x + \tfrac{23}{2}, & \quad\text{otherwise},
            \end{cases}\\
            f_B(x)
                {}\coloneqq{}&
            \begin{cases}
                2x^2 - x, & \quad\text{if } x < -1, \\
                x^2 + 3x + 5, & \quad\text{if } x \in [-1, 2], \\
                \makebox[\widthof{$-\tfrac12 x^2 - 3x + \tfrac{23}{2},$}][l]{$x^2 + 15x - 19,$}
                & \quad\text{if } x > 2. \\
            \end{cases}
        \end{align*}
        The following hold for the inclusion \(0 \in \Tp x = \partial f_A(x) + \partial f_B(x)\).
    \end{minipage}%
    \raisebox{0pt}[0pt][0pt]{
        \begin{minipage}[t]{\widthfig}
            \ifspringer
                \vspace{-0.5cm}
            \else
                \vspace{-0.175cm}
            \fi
            \centering
            \includetikz{Examples/Nonsmooth/nonsmooth}
            \ifspringer
                \vspace{-0.4cm}
            \else\fi
            \captionof{figure}{
                Visualization of the (nonsmooth) functions involved in \Cref{ex:nonsmooth:minimization}.
            }
            \label{fig:nonsmooth:minimization}
        \end{minipage}%
    }
    \ifspringer
        \vspace{0.225\baselineskip}
    \else\fi
    \begin{enumerate}
        \item 
        \Cref{ass:DRS:1} and \Cref{ass:DRS:2} hold when $\gamma \neq 1$.
        \item \Cref{ass:SWMVIDRSsemimonotone} holds for $\pazocal{S}^\star = \bigl\{(0, 3)\bigr\}$ and \(\monA = -1.2\), \(\comA = 0.2\), \(\monB = 1.6\) and \(\comB = 0.1\). 
        \item 
                The claims from \cref{thm:pppaDRS}
        hold for the sequence $\seq{u^k, v^k, s^k}$ generated by \ref{eq:DRS} with stepsize
        \(
            \gamma
                {}\in{}
            \bigl(0, 0.2614\bigr)
        \)
        and relaxation parameter
        $\lambda\in\bigl(0,2 + \nicefrac{2}{15\gamma} - \nicefrac{48\gamma}{5}\bigr)$.
        \qedhere
    \end{enumerate}
\end{example}

\begin{remark}[comparison with previous works]\label{rem:nonsmoothopt:comparison}
    \begin{enumerate}
        \item Since both $f_A$ and $f_B$ are nonsmooth
        and $f_A$ is nonconvex,
        convergence results for \ref{eq:DRS} in the optimization setting such as \cite{li2016Douglas, themelis2020DouglasRachford} cannot be applied.
        \item Since $\comA > 0$ and $\comB > 0$, it holds by definition that $\partial f_A$ is $\monA$\hyp{}monotone at $(x^\star, -y^\star)$ and $\partial f_B$ is $\monB$\hyp{}monotone at $(x^\star, y^\star)$. However, as the monotonicity of $A$ and $B$ 
        does not hold globally,
        the results from \cite{dao2019Adaptive} are not applicable here.
        \qedhere
    \end{enumerate}
\end{remark}

      Next,
consider the following example, where we apply \ref{eq:DRS} to the subdifferential of a separable nonsmooth function. 
For this example, we show that our analysis of \ref{eq:DRS} covers its convergence behavior towards the stationary points,
and provide numerical results in \Cref{fig:nonsmooth:stationary}. Notice that in the iterates are attracted by the local minimum, while in contrast e.g. descent methods would quickly diverge.

\begin{figure}
    \centering
    \includetikz{Examples/Stationary/stationary}
    \caption{(left) Visualization of the functions involved in \Cref{ex:nonsmooth:stationary}. Additionally, the sequence $(u^k)_{k=0}^{400}$ is indicated using black arrows, which is obtained by applying \ref{eq:DRS} to \(0 \in \partial f_A(x) + \partial f_B(x)\) with 
    stepsize $\gamma = \nicefrac{11}{60} \in (\nicefrac16, \nicefrac15)$ and relaxation 
    $\lambda = \nicefrac{9}{10}\bigl(2 - \nicefrac{4}{15\gamma} - \nicefrac{12\gamma}{5}\bigr) \approx 0.095$. 
    When evaluating the (multi-valued) resolvents $J_{\gamma A}$ and $J_{\gamma B}$, each element is sampled with equal probability.
    (right) Visualization of the sequence $(u^k)_{k=0}^{400}$ for $200$ such experiments, of which the initializations $s^0$ are evenly spaced within the interval $[-5,5]$.}
    \label{fig:nonsmooth:stationary}
\end{figure}

\begin{example}[stationary point]\label{ex:nonsmooth:stationary}
    Let \(\monA = -0.3, \comA = -0.1, \monB = 0.4\) and \(\comB = 0.4\). Consider the problem of finding a stationary point of
    \(
            f_A(x) + f_B(x)
    \),
    visualized in \Cref{fig:nonsmooth:stationary}, where
    \begin{align*}
        f_A(x)
            {}\coloneqq{}&
        \begin{cases}
            \tfrac{1-\sqrt{1-4\monA\comA}}{4\comA}x^2 - x - 3, & \quad\text{if } x \in [-3, 3], \\
            \tfrac{1+\sqrt{1-4\monA\comA}}{4\comA}x^2 - x - \tfrac{9\sqrt{1-4\monA\comA}}{2\comA} - 3, & \quad\text{otherwise},
        \end{cases}\\
        f_B(x)
            {}\coloneqq{}&
        \begin{cases}
            \makebox[\widthof{$\tfrac{1+\sqrt{1-4\monA\comA}}{4\comA}x^2 - x - \tfrac{9\sqrt{1-4\monA\comA}}{2\comA} - 3,$}][l]{$x^2 + x + 2,$} & \quad\text{if } x \in [-1, 1], \\
            \tfrac{1}{4}x^2 + x + \tfrac{11}{4}, & \quad\text{otherwise}.
        \end{cases}
    \end{align*}
    The following hold for the inclusion \(0 \in \Tp x = \partial f_A(x) + \partial f_B(x)\).
    \begin{enumerate}
        \item 
        \Cref{ass:DRS:1} and \Cref{ass:DRS:2} hold for 
        \( 
            \smash{
            \gamma \in \bigl(0, \tfrac{-2\comA}{1+\sqrt{1-4\monA\comA}}\bigr) \cup \bigl(\tfrac{-2\comA}{1+\sqrt{1-4\monA\comA}}, \nicefrac15\bigr)}.
        \)
        \item \Cref{ass:SWMVIDRSsemimonotone} holds for $\pazocal{S}^\star = \bigl\{(0, -1)\bigr\}$ and the given \(\monA, \comA, \monB, \comB\). 
        \item 
                The claims from \cref{thm:pppaDRS}
        hold for the sequence $\seq{u^k, v^k, s^k}$ generated by \ref{eq:DRS} with stepsize
        \(
            \gamma
                {}\in{}
            \bigl(\nicefrac16, \nicefrac15\bigr)
        \)
        and fixed relaxation parameter
        $\lambda\in\bigl(0,2 - \nicefrac{4}{15\gamma} - \nicefrac{12\gamma}{5}\bigr)$ (see \Cref{fig:nonsmooth:stationary} for numerical results).
        \qedhere
    \end{enumerate}
\end{example}

    As a final example, consider the following quadratic programming (QP) problem, with box and linear equality constraints.
    When $Q$ is positive semidefinite, which is the convex setting, convergence of \ref{eq:DRS} for \Cref{ex:constrainedQP} follows directly by monotonicity of $A$ and $B$.
    The following example demonstrates that our results allow to go beyond the monotone setting, leveraging the calculus rules from \Cref{sec:semi}.
    Note that in this example both $f$ and $g$ are nonsmooth, and thus convergence results for \ref{eq:DRS} such as \cite{li2016Douglas, themelis2020DouglasRachford} cannot be applied here.

    \begin{example}[constrained QP]\label{ex:constrainedQP}
        Consider the following quadratic program
        \begin{equation}
            \begin{aligned}
                \minimize_{x\in\R^n} \quad& \tfrac12 x^\top Q x + q^\top x\\
                \stt \quad& Hx = h,\\
                & l \leq x \leq u,
            \end{aligned}
            \label{eq:ex:constrainedQP:problem}
        \end{equation}
        where $Q \in \sym{n}$, $q \in \R^n$, $H \in \R^{m \times n}$ has full row rank, $h \in \R^m$ and scalars $l, u \in \R$ satisfying $l < u$.
        This problem can be equivalently formulated as
        $
            \minimize_{x} f(x) + g(x)
        $
        where
        \[
            f(x) \coloneqq \tfrac12 x^\top Q x + q^\top x + \delta_{\defset{s}{Hs = h}}(x)
                \quad\text{and}\quad
            g(x) \coloneqq \delta_{\defset{s}{l \leq s \leq u}}(x),
        \]
        and the associated first-order optimality condition is given by
        $
            0 \in A(x) + B(x) \coloneqq \partial f(x) + \partial g(x),
        $
        By splitting the quadratic program \eqref{eq:ex:constrainedQP:problem} into $f$ and $g$, the resulting \ref{eq:DRS} iterations can be computed efficiently in closed form. In particular, $f$ corresponds to an equality constrained QP problem, of which the prox can be solved efficiently through its corresponding KKT system \cite[\S 16.1]{nocedal2006numerical}, and computing the resolvent of $\partial g$ corresponds to projecting onto a box.
    
        For instance, let $Q = \diag(-1,\tfrac12)$, $q = \begin{bsmallmatrix} -4 & -1\end{bsmallmatrix}^\top$, 
        $
            H
                {}={} 
            \begin{bsmallmatrix}
                1 & 1
            \end{bsmallmatrix}
        $,
        $h = 1$,
        $l = 0$ and $u = 1$.
        Then, the unique minimizer of \eqref{eq:ex:constrainedQP:problem} is given by $x^\star = \begin{bsmallmatrix} 1 & 0\end{bsmallmatrix}^\top$ and the following assertions hold.
        \begin{enumerate}
            \item Operator $A$ is $(-1, 0)$-semimonotone.
            \item Operator $B$ is $(2, 0)$-semimonotone at 
            $
                (x^\star, y^\star) 
                    = 
                (
                    \begin{bsmallmatrix} 1 & 0\end{bsmallmatrix}^\top,
                    \begin{bsmallmatrix} 2 & -2\end{bsmallmatrix}^\top
                ) 
                    \in 
                \gph B
            $.
            \item The sequence $\seq{u^k,\nicefrac1{\gamma}(u^k - s^k)}$ generated by \ref{eq:DRS} with fixed relaxation parameter $\lambda$ converges to $\zer \Tpd$ for 
            $
                \gamma \in (0, \nicefrac12)
            $
            and
            $
                \lambda \in 
                (0, 2 - 4\gamma).
            $
        \end{enumerate}
        \begin{proof}
            By \Cref{prop:calculus:scaling,it:prop:semi:lin:1}, it follows that the mapping $x \mapsto Qx + q$ is $(-1,0)$\hyp{}semimonotone.
            The claimed semimonotonicity of $A$ then follows by observing that the normal cone $N_{\defset{s}{Hs = h}}$ is monotone \cite[Ex. 20.26]{bauschke2017Convex} and applying \Cref{lem:sum:1}.
            Finally, the claimed semimonotonicity of $B$ follows from \Cref{prop:box} and the convergence result for \ref{eq:DRS} follows from \Cref{cor:DRS:semi}.
        \end{proof}
    \end{example}

  \section{Conclusion}\label{sec:conclusion}

This paper presented a comprehensive study of the (relaxed) preconditioned proximal point algorithm (PPPA) and the (relaxed) Douglas-Rachford splitting (DRS) method for a class of nonmonotone problems that satisfy an oblique weak Minty condition.
To make the achieved results for DRS readily applicable in practice, the class of semimonotone operators was introduced, and sufficient conditions were developed for the convergence of DRS for the sum of two semimonotone operators.
Various relevant example problems were provided thoughout the paper, showing tightness of the achieved convergence results.

Future research directions include investigating inexact variants of PPPA, 
analyzing other splitting methods in nonmonotone settings,
and studying applications in noncooperative game theory and nonconvex distributed optimization.

  \ifspringer
      \backmatter

      \section*{Statements and Declarations}

      \paragraph{Competing Interests}
      The authors declare no conflict of interest.

      \paragraph{Code Availability}
      The code for the numerical examples is available on \href{https://github.com/brechtevens/Minty-DRS-examples}{GitHub}.
  \fi

  \begin{appendix}\label{sec:appendix}
  \section{Auxiliary lemmas}\label{sec:auxiliary}

\begin{lemma}
	Let $\M \in \sym{n}_{+}$
	be a symmetric positive semidefinite matrix, $U$ be an orthonormal basis for $\range{\M}$ and define $\projR \coloneqq \M + \proj_{\ker{\M}}$. Then, the following properties hold.
	\begin{enumerate}
		\item\label{it:lem:PQ:properties:proj} The projection onto the range of $\M$ is given by $\projQ = U \tp U$.
		\item\label{it:lem:PQ:properties:equivalences}
		\(
			\M = \M \projQ = \projQ \M = \projR \projQ
		\).
		\item\label{it:lem:PQ:properties:pd} $Q \succ 0$.
		\item\label{it:lem:PQ:properties:decomposition}
		 $\tp U P U = \tp U \projR U\succ 0$.
	\end{enumerate}
	\begin{proof}
		Statement \ref{it:lem:PQ:properties:proj} follows from \cite[\S 5.13]{meyer2000matrix} and \ref{it:lem:PQ:properties:equivalences} follows directly by definition of $\M$ and $\projR$.
		\begin{proofitemize}
			\item \ref{it:lem:PQ:properties:pd}: 
			Let \(x\in \R^n \backslash \{0\}\) and consider decomposing it as \(x = x_1  + x_2\) where \(x_1 \in \range{\M}\) and \(x_2 \in \ker{\M}\).
						Let $U_1 \Lambda_1 \tp U_1$ be a compact eigendecomposition of $\M$, so that $U_1 \tp U_1 = \projQ$ and $p_{\rm min} \coloneqq \lambda_{\rm min}(\Lambda_1) > 0$.
				Then,
				\ifspringer
					\(
						\inner{x, \projR x} 
							{}={}
						\inner{x, P x}
							{}+{}
						\inner{x, \proj_{\ker{\M}} x}
							{}={}
						\inner{x_1, U_1 \Lambda_1 \tp U_1 x_1}
							{}+{}
						\|x_2\|^2
							{}\geq{}
						p_{\rm min}\inner{x_1, \proj_{\range{\M}} x_1}
							{}+{}
						\|x_2\|^2
							{}={}
						p_{\rm min}\nrm{x_1}^2
							{}+{}
						\|x_2\|^2
							{}\geq{}
						\min\{1, p_{\rm min}\} \nrm{x}^2 > 0.
					\)
				\else
					\begin{align*}
						\inner{x, \projR x} 
							{}={}&
						\inner{x, P x}
							{}+{}
						\inner{x, \proj_{\ker{\M}} x}
							{}={}
						\inner{x_1, U_1 \Lambda_1 \tp U_1 x_1}
							{}+{}
						\|x_2\|^2\\
							{}\geq{}&
						p_{\rm min}\inner{x_1, \proj_{\range{\M}} x_1}
							{}+{}
						\|x_2\|^2
							{}={}
						p_{\rm min}\nrm{x_1}^2
							{}+{}
						\|x_2\|^2
							{}\geq{}
						\min\{1, p_{\rm min}\} \nrm{x}^2 > 0.
					\end{align*}
				\fi
			\item \ref{it:lem:PQ:properties:decomposition}: 
			The equivalence $\tp U P U = \tp U \projR U$ follows from \Cref{it:lem:PQ:properties:proj,it:lem:PQ:properties:equivalences} using $\tp U U = \I$. Positive definiteness follows from \Cref{it:lem:PQ:properties:pd} and the fact that $U$ has full column rank \cite[Prop. 8.1.2(xiii)]{bernstein2009matrix}.
			\qedhere
		\end{proofitemize}
	\end{proof}
\end{lemma}

\begin{proposition}[resolvent of a comonotone operator]\label{prop:comon:fulldom}
    Let $T: \R^n \rightrightarrows \R^n$ be a $\DRSRho$\hyp{}comonotone operator and suppose that \(\M\) is selected according to \Cref{ass:PPPA:2}.
    Then, the preconditioned resolvent \(\PRES\) has full domain if and only if \(T\) is maximal.
    \begin{proof}
        Since $T$ is $\DRSRho$\hyp{}comonotone, it holds by definition that 
        \(
            T^{-1} - \DRSRho 
        \)
        is monotone and thus also that
        \(
            \tilde T \coloneqq 
            U^\top T^{-1} U - U^\top \DRSRho U
        \)
        is monotone.
        Furthermore, it holds that \((U^\top \M U)^{-1} + U^\top \DRSRho U \succ 0\) by \Cref{ass:PPPA:2,rem:ass:PPPA:eigcond}.
        Denote 
        \(
            \hat T \coloneqq U^\top T^{-1} U + (U^\top \M U)^{-1}
        \). 
        Owing to \((U^\top \M U)^{-1} + U^\top \DRSRho U \succ 0\) and Minty's theorem \cite[Thm. 21.1]{bauschke2017Convex}, the operator 
        \[
            \bigl((U^\top \M U)^{-1} + U^\top \DRSRho U\bigr)^{-1}\hat T 
                {}={}
            \bigl((U^\top \M U)^{-1} + U^\top \DRSRho U\bigr)^{-1}\tilde T
                {}+{}
            \id
        \]
        has full range iff \(\bigl((U^\top \M U)^{-1} + U^\top \DRSRho U\bigr)^{-1}\tilde T\) is maximally monotone.
        Since \((U^\top \M U)^{-1} + U^\top \DRSRho U \succ 0\), this implies that \(\hat T\) has full range iff \(\tilde T\) is maximally monotone, i.e., iff \(T\) is maximally $\DRSRho$\hyp{}comonotone.

        Therefore, it only remains to show that \(\hat T\) has full range iff 
        \(\dom{\PRES} = \R^n\).
        First suppose that \(\hat T\) has full range.
        Let $m$ denote the dimension of the range of \(\M\) and fix \(\hat{y}\in\R^m\).
        Then, by the full range assumption for \(\hat T\), there exists some \(\hat{x}\in\R^m\) such that \(\hat{y} \in \hat{T} \hat{x}\). 
        Letting \(\hat{z} = \hat{y} - (U^\top \M U)^{-1} \hat{x}\), this implies that 
        \begin{align*}
            \hat{y} \in U^\top T^{-1} U \hat{x} + (U^\top \M U)^{-1} \hat{x}
                {}\iff{}
            \hat{z} \in U^\top T^{-1} U \bigl(U^\top \M U (\hat{y} - \hat{z})\bigr)
            \numberthis\label{eq:prop:comon:fulldom:equiv:1}
        \end{align*}
        Multiplying by $U$, this implies that \(U \hat{z} \in \projQ T^{-1} \bigl(\M U (\hat{y} - \hat{z})\bigr)\).
        Defining \(x' = U \hat{x}\), \(y' = U \hat{y}\) and \(z' = U \hat{z}\), this corresponds to 
        \(
            z' \in \projQ T^{-1} \bigl(\M (y' - z')\bigr).
        \)
        Consequently, \(T^{-1} \bigl(\M (y' - z')\bigr)\) is nonempty and there exists a $z'' \in \ker{\M} \subseteq \R^n$ so that
        \(
            z'' \in \proj_{\ker{\M}} T^{-1} \bigl(\M (y' - z')\bigr).
        \)
        Summing both inclusions, it holds that
        \begin{align*}
            z' + z'' \in T^{-1} \bigl(\M (y' - z')\bigr)
                {}\iff{}
            \M (y' - z') \in T (z' + z'')
                {}\iff{}
            z' + z'' \in \big(\M + T\big)^{-1} \M y'
            \numberthis\label{eq:prop:comon:fulldom:equiv:2}
        \end{align*}
        Since the choice of \(y'\) is arbitrary in \(\range{\M}\), \(\bigl(\M + T\bigr)^{-1} \M\) has full domain. 

        The converse follows similarly.
        Assume that \(\dom{\bigl(\M + T\bigr)^{-1} \M}=\R^n\).
        Then, for any fixed \(\hat{y}\in\R^m\) and \(y' = U \hat{y}\) there exists some \(z\in\R^n\) such that \(z \in (T + \M)^{-1}\M y\).
        Without loss of generality, decompose \(z \coloneqq U \hat{z} + z''\), where $\hat z \in \R^m$ and \(z'' \in \ker{P}\). Then, by \eqref{eq:prop:comon:fulldom:equiv:2} it follows that
        \(
            U \hat{z} + z'' \in T^{-1} \bigl(\M U (\hat{y} - \hat{z})\bigr)
        \).
        Multiplying by $U^\top$, this implies that
        \(
            \hat{z} \in U^\top T^{-1} U \bigl(U^\top \M U (\hat{y} - \hat{z})\bigr)
        \).
        Defining $\hat{x} = U^\top P U (\hat{y} - \hat{z})$ and using \eqref{eq:prop:comon:fulldom:equiv:1}, we obtain \(\hat{y}\in \hat{T}\hat{x}\), thus establishing that \(\hat T\) has full range.
    \end{proof}
\end{proposition}

\begin{fact}[Solution of quadratic inequality]\label{lem:quadratic}
    Let \(a,b,c\in\R\) and suppose that $b > -2\sqrt{ac}$ when $ac \geq 0$.
    Then, the following hold.
    \begin{enumerate}
        \item \label{it:lem:quadratic:existence} There exists a $\gamma > 0$ satisfying
        \(
            a\gamma^2 + b\gamma + c > 0
        \)
        if and only if it holds that
        \begin{align*}
            b > 2\sqrt{ac}
            \quad
            \text{when}
            \quad
            a < 0 \text{ and } c < 0.\numberthis\label{eq:cor:quadratic:existence:condition}
        \end{align*}
        \item \label{it:lem:quadratic:solutions} 
        \iftables
            If \eqref{eq:cor:quadratic:existence:condition} holds, then all strictly positive solutions $\gamma > 0$ of \eqref{eq:lem:quadratic} are given by \Cref{tab:gamma:general}.
        \else
            If \eqref{eq:cor:quadratic:existence:condition} holds, then $\gamma > 0$ satisfies 
            \(
                a\gamma^2 + b\gamma + c > 0
            \)
            \ifspringer
                iff
            \else
                if and only if
            \fi
            \(
                \smash{
                \gamma \in 
                \left(
                    \tfrac{2[-c]_{+}}{b + \sqrt{b^2 - 4ac}}, \tfrac{b + \sqrt{b^2 - 4ac}}{2[-a]_{+}}
                \right)
                }
            \),
            which is nonempty.
        \fi
    \end{enumerate}
    \iftables
        \begin{table}
            \centering
            \begin{tabular}{c}
                {
                \begin{tblr}{
                    hlines,
                    vlines,
                    vline{1-1} = {1-1}{0pt},
                    hline{1-1} = {1-1}{0pt},
                    colspec={M{2cm}M{4.25cm}M{4.25cm}},
                    cell{1-1}{2-3}={c,gray!20},
                    cell{2-3}{1-1}={c,gray!20}
                    }
                    &$a \geq 0$ & $a < 0$\\
                    $c \geq 0$&$\gamma \in \left(0,+\infty\right)$&$\gamma \in \left(0, -\tfrac{b + \sqrt{b^2 - 4ac}}{2a}\right)$\\
                    $c < 0$&$\gamma \in \left(-\tfrac{2c}{b + \sqrt{b^2 - 4ac}},+\infty\right)$&$\gamma \in \left(-\tfrac{2c}{b + \sqrt{b^2 - 4ac}},-\tfrac{b + \sqrt{b^2 - 4ac}}{2a}\right)$
                \end{tblr}
                }
            \end{tabular}
            \caption{Range of $\gamma$ satisfying \eqref{eq:lem:quadratic}.}
            \label{tab:gamma:general}
        \end{table}
    \fi
\end{fact}


\begin{fact}\label{prop:parsum:order}
    Let $a, b \in \R$ such that $a + b > 0$. Then, $a \Box b < \min\{a,b\}$.
\end{fact}

\begin{lemma}[implications of \Cref{ass:DRS}]\label{lem:implication:assDRS:assPPPA}
    If \Cref{ass:DRS} holds, then \Crefrange{ass:PPPA:0}{ass:PPPA:1} hold for operator $\Tpd$ from \eqref{eq:primaldual} and preconditioner $\M$ as in \eqref{eq:M}.
    \begin{proof}
        First, outer semicontinuity of $\Tpd$ follows from that of $A$ and $B$ \cite[Thm. 5.7(a)]{rockafellar2009Variational}, showing \ref{ass:PPPA:0}. Second, \ref{ass:PPPA:0.5} holds since $J_{\gamma A}$ and $J_{\gamma B}$ having full domain implies that 
        the preconditioned resolvent $(\M + \Tpd)^{-1}\M$ has full domain, owing to \cite[Lem. 12.14]{rockafellar2009Variational}.
        Finally, statement \ref{ass:PPPA:1} is immediate, completing the proof.
    \end{proof}
\end{lemma}

\begin{lemma}[implications of \Cref{ass:SWMVIDRSsemimonotone:maximal}]\label{prop:SWMVIDRSsemimonotone:maximal:fulldom}
    Suppose that \Cref{ass:SWMVIDRSsemimonotone:maximal} holds and that $\gamma$ is selected according to \Cref{tab:gamma:DRS:semi}, where $\gamma_{-}$ and $\gamma_{+}$ are defined as in \Cref{eq:eigcondDRS:gamma:semi}.
    Then, 
        \Cref{ass:DRS:1,ass:DRS:2} hold
    and the resolvents $J_{\gamma A}$ and $J_{\gamma B}$ are Lipschitz continuous.
    \begin{proof}
        First, observe that $[\monA \Box \monB]_-[\comA \Box \comB]_- < \tfrac{1}{4}$, which by \Cref{prop:parsum:order}
        and \Cref{it:prop:semi:young:2} implies that
        $[\monA]_-[\comA]_- < \tfrac{1}{4}$ and $[\monB]_-[\comB]_- < \tfrac{1}{4}$.
                Therefore, owing to \Cref{prop:semimon:osc}, it follows that $A$ and $B$ are outer semicontinuous, showing that \Cref{ass:DRS:1} holds.
            Furthermore,
        owing to \Cref{cor:semimon:resolvent}, it follows that the resolvents $J_{\gamma A}$ and $J_{\gamma B}$ are Lipschitz continuous and have full domain if $\gamma$ complies with \Cref{tab:gamma:DRS:semi:fulldomain}, where
        \begin{align}\label{eq:SWMVIDRSsemimonotone:maximal:fulldom:bounds}
            \bar\gamma_{-} = 
                \begin{cases}
                    \tfrac{-2\comA}{1+\sqrt{1-4\monA\comA}}, & \quad \text{if } \comA < 0,\\
                    \tfrac{-2\comB}{1+\sqrt{1-4\monB\comB}}, & \quad \text{otherwise},
                \end{cases}
            \quad \text{and} \quad
            \bar\gamma_{+} =
                \begin{cases}
                    -\tfrac{1+\sqrt{1-4\monA\comA}}{2\monA}, & \quad \text{if } \monA < 0,\\
                    -\tfrac{1+\sqrt{1-4\monB\comB}}{2\monB}, & \quad \text{otherwise}.
                \end{cases}
        \end{align}
        As $\gamma$ is selected according to \Cref{tab:gamma:DRS:semi}, it simply remains to prove that these ranges are a subset of the ranges for $\gamma$ provided in \Cref{tab:gamma:DRS:semi:fulldomain}, which can be verified directy using \Cref{prop:parsum:order}.
        \begin{table}
            \centering
            \caption{Ranges of $\gamma$ for which the resolvents $J_{\gamma A}$ and $J_{\gamma B}$ have full domain.}
            \label{tab:gamma:DRS:semi:fulldomain}
            \begin{tabular}{c}
                {
                \ifspringer
                    \begin{tblr}{
                        hlines,
                        vlines,
                        vline{1-2} = {1-2}{0pt},
                        hline{1-2} = {1-2}{0pt},
                        colspec={M{1.6cm}M{1.3cm}M{2.6cm}M{2.1cm}M{2.6cm}},
                        cell{1-2}{3-5}={c,gray!20},
                        cell{3-5}{1-2}={c,gray!20}
                        }
                        &&\SetCell[r=2]{c}$\monA = \monB = 0$ & \SetCell[c=2]{c}$\monA + \monB > 0$\\
                        &&& $\monA\monB \geq 0$ & $\monA\monB < 0$\\
                        \SetCell[c=2]{c}$\comA = \comB = 0$&&\SetCell[c=2, r=2]{c}$\gamma \in \left(0,+\infty\right)$&&$\gamma \in \left(0, -\tfrac{1}{\min\{\monA,\monB\}}\right)$\\
                        \SetCell[r=2]{c}$\comA + \comB > 0$&$\comA\comB \geq 0$&&&$\gamma \in \left(0, \bar\gamma_{+}\right)$\\
                        &$\comA\comB < 0$&$\gamma \in \left({\scriptstyle-\min\{\comA,\comB\}}, +\infty\right)$&$\gamma \in \left(\bar\gamma_{-}, +\infty\right)$&$\gamma \in \left(\bar\gamma_{-}, \bar\gamma_{+}\right)$
                    \end{tblr}
                \else
                    \begin{tblr}{
                        hlines,
                        vlines,
                        vline{1-2} = {1-2}{0pt},
                        hline{1-2} = {1-2}{0pt},
                        colspec={M{1.8cm}M{1.4cm}M{3.1cm}M{2.2cm}M{2.75cm}},
                        cell{1-2}{3-5}={c,gray!20},
                        cell{3-5}{1-2}={c,gray!20}
                        }
                        &&\SetCell[r=2]{c}$\monA = \monB = 0$ & \SetCell[c=2]{c}$\monA + \monB > 0$\\
                        &&& $\monA\monB \geq 0$ & $\monA\monB < 0$\\
                        \SetCell[c=2]{c}$\comA = \comB = 0$&&\SetCell[c=2, r=2]{c}$\gamma \in \left(0,+\infty\right)$&&$\gamma \in \left(0, -\tfrac{1}{\min\{\monA,\monB\}}\right)$\\
                        \SetCell[r=2]{c}$\comA + \comB > 0$&$\comA\comB \geq 0$&&&$\gamma \in \left(0, \bar\gamma_{+}\right)$\\
                        &$\comA\comB < 0$&$\gamma \in \left({\scriptstyle-\min\{\comA,\comB\}}, +\infty\right)$&$\gamma \in \left(\bar\gamma_{-}, +\infty\right)$&$\gamma \in \left(\bar\gamma_{-}, \bar\gamma_{+}\right)$
                    \end{tblr}
                \fi
                }
            \end{tabular}
        \end{table}
    \end{proof}
\end{lemma}


\begin{proposition}[normal cone of a box]\label{prop:box}
    The normal cone operator $N_C: \R^n \rightrightarrows \R^n$ of the $n$\hyp{}dimensional box 
    \(
        C \coloneqq \set{x\in \R^n}[l_i \leq x_i \leq u_i, \; i = 1, \dots, n]
    \) is 
    $
        \Bigl(\min\Bigl\{\frac{|\tilde v_1|}{u_1 - l_1}, \dots, \frac{|\tilde v_n|}{u_n - l_n}\Bigr\}, 0\Bigr)
    $-semimonotone at $(\tilde x, \tilde v) \in \gph N_C$.
    \begin{proof}
        By monotonicity of the normal cone $N_C$ \cite[Ex. 20.26]{bauschke2017Convex}, it holds that
        \begin{align*}
            \inner{\tilde v, \tilde x - x}
                {}={}&
            \sum_{i=1}^n
            |\tilde v_i||\tilde x_i - x_i|
                {}\geq{}
            \sum_{i=1}^n
            \tfrac{|\tilde v_i|}{u_i - l_i}|\tilde x_i - x_i|^2
                {}\geq{}
            \min\Bigl\{\tfrac{|\tilde v_1|}{u_1 - l_1}, \dots, \tfrac{|\tilde v_n|}{u_n - l_n}\Bigr\}
            \nrm{\tilde x - x}^2,
        \end{align*}
        where the first inequality holds since \(|\tilde x_i - x_i| \leq u_i-\ell_i\) for all $i \in \set{1,\hdots,n}$.
    \end{proof} 
\end{proposition}

  \section{Omitted proofs}
    \subsection{Douglas-Rachford splitting}

\begin{appendixproof}{lem:equivalence:PPPA-DRS}
    The update rule \eqref{eq:PPPA-intro} applied to \eqref{eq:primaldual} with preconditioner \eqref{eq:M}  can be equivalently written as 
    \begin{subequations}
        \begin{align}
            \bar{x}^k
                {}\in{}&
            J_{\gamma A}(x^k - \gamma y^k)\label{eq:lem:equivalences:PPPA-DRS:1}
            \\
            \bar{y}^k 
                {}\in{}&
            (\gamma \id + B^{-1})^{-1}(2\bar{x}^k - x^k + \gamma y^k)
            \Leftrightarrow
            2\bar{x}^k - x^k + \gamma y^k - \gamma\bar{y}^k
                {}\in{}
            J_{\gamma B}(2\bar{x}^k - x^k + \gamma y^k)\label{eq:lem:equivalences:PPPA-DRS:2}\\
            x^{k+1} {}={}& x^k + \lambda_k(\bar{x}^k - x^k)\label{eq:lem:equivalences:PPPA-DRS:3}\\
            y^{k+1} {}={}& y^k + \lambda_k(\bar{y}^k - y^k),\label{eq:lem:equivalences:PPPA-DRS:4}
        \end{align}
    \end{subequations}
    where \cite[Lem. 12.14]{rockafellar2009Variational} was used in \eqref{eq:lem:equivalences:PPPA-DRS:2}.
    We proceed by induction. 
    The base case $k = 0$ is vacuously satisfied through initialization $s^0 = x^0 - \gamma y^0$, and by using \eqref{eq:lem:equivalences:PPPA-DRS:1} and \eqref{eq:lem:equivalences:PPPA-DRS:2}, noting that for any \(\bar x^0\), \(\bar y^0\) there exists \(u^0 = \bar x^0 \in J_{\gamma A}(s^0)\), \(v^0= 2u^0-s^0 - \gamma \bar y^0 \in J_{\gamma B}(2u^0 - s^0)\) (and vice versa). 
    Suppose now that the claim holds for some $k \geq 0$. In particular, \(s^k = x^k - \gamma y^k\), $u^k = \bar{x}^k$ and $v^k = 2 u^k - s^k - \gamma\bar{y}^k$. 
    According to update rule \eqref{eq:DRS}, it holds that 
    \ifspringer
        \(
            s^{k+1} 
            {}={}
            s^k + \lambda_k (v^k-u^k)
            {}={}
            x^k - \gamma y^k + \lambda_k (\bar x^k - \gamma\bar{y}^k - x^k + \gamma y^k)
            {}={}
            x^{k+1} - \gamma y^{k+1},
        \)
    \else
        \begin{align*}
            s^{k+1} 
            {}={}
            s^k + \lambda_k (v^k-u^k)
            {}={}
            s^k + \lambda_k (\bar x^k - s^k - \gamma\y^k)
            {}={}
            x^k - \gamma y^k + \lambda_k (\bar x^k - \gamma\y^k - x^k + \gamma y^k)
            {}={}
            x^{k+1} - \gamma y^{k+1},
        \end{align*}
    \fi
    where the induction hypothesis was used in the first two equalities, and  \eqref{eq:lem:equivalences:PPPA-DRS:3} and \eqref{eq:lem:equivalences:PPPA-DRS:4} were used in the last equality. 
    Having established \(s^{k+1} = x^{k+1} - \gamma y^{k+1}\), the claim for the correspondences between \(\bar x^{k+1}, \bar y^{k+1}\) and \(u^{k+1}, v^{k+1}\) follow from \eqref{eq:lem:equivalences:PPPA-DRS:1} and \eqref{eq:lem:equivalences:PPPA-DRS:2} as argued in the base case.  
\end{appendixproof}

\begin{appendixproof}{lem:SWMI:primal-dual}
    Let $(x_\A, y_\A) \in \graph \A$ and $(x_\B, y_\B) \in \graph \B$. Then, by definition of the primal-dual operator and using that $(y_\B, x_\B) \in \graph \B^{-1}$ it holds that 
    \(
        (z, v) \coloneqq \Bigl((x_\A,y_\B), (y_\A + y_\B, x_\B - x_\A)\Bigr) \in \graph \Tpd.
    \)
    Therefore, by \Cref{ass:SWMVIDRS} it holds for all $z^\star = (x^\star, y^\star) \in \pazocal{S}^\star$ that \(\inner*{v, z - z^\star} \geq \qindef{v}{\DRSRho}\), i.e.,
    \begin{align*}
        \inner*{y_\A + y_\B, x_\A - x^\star} + \inner*{x_\B - x_\A, y_\B - y^\star} \geq \DRSrhocom\nrm{y_\A + y_\B}^2 + \DRSrhomon\nrm{x_\B - x_\A}^2.\numberthis\label{eq:WMI:primal-dual}
    \end{align*}

    Let $x_\A = x_\B \eqqcolon x_P$. Then, from \eqref{eq:WMI:primal-dual} it follows that
    \(
        \inner*{y_\A + y_\B, x_P - x^\star} \geq \DRSrhocom\nrm{y_\A + y_\B}^2
    \) 
    for all $x^\star \in \pazocal{S}^\star_P$.
    Since
    \(
        (x_P, y_\A + y_\B) \in \graph \Tp
    \)
    it is of immediate verification that $\Tp$ has $\DRSrhocom$-weak Minty solutions at $\pazocal{S}^\star_P$.

    Similarly, let $-y_\A = y_\B \eqqcolon y_D$. Then, due to \eqref{eq:WMI:primal-dual}
    \(
        \inner*{x_\B - x_\A, y_D - y^\star} \geq \DRSrhomon\nrm{x_\B - x_\A}^2
    \)
    for all $
        y
^\star \in \pazocal{S}^\star_D$ and the claim for $\Td$ follows analogously from the fact that
    \(
        (y_D, x_\B - x_\A) \in \graph \Td
    \)%
    , completing the proof.
\end{appendixproof}

\begin{appendixproof}{cor:DRS:CP:lin}
    \phantom{.}\newline
    Since operator \(B\) is piecewise polyhedral so is its inverse, and in turn so is the mapping \(\hat T = A \times B^{-1}\). Therefore, \(T_{\rm PD}\) being the sum of \(\hat T\) and a linear mapping is also piecewise polyhedral. Since the inverse of a piecewise polyhedral mapping is piecewise polyhedral, it follows from \cite[3H.1 and 3H.3]{dontchev2009implicit} that \(\Tpd\) is metrically subregular at any \(z\) for any \(v\) with \((z,v)\in \graph \Tpd\). The claimed rates for \ref{eq:DRS} are an immediate consequence of \cref{thm:PPPA:linear}, owing to the relations \(U_1^\top z^k = s^k\) and \eqref{eq:DRSequivalence:seminorm}. 
\end{appendixproof}

    \subsection{Semimonotone operators}

\begin{appendixproof}{prop:semi:young}
        Consider $(x, y),(\other{x}, \other{y}) \in \graph A$ and let $\epsilon > 0$. Then, by the Fenchel-Young inequality, it holds that
        \begin{align*}
            -\epsilon\nrm{x - \other{x}}^2 - \epsilon^{-1}\nrm{y - \other{y}}^2 \leq 2\inner{x-\other{x},y-\other{y}} \leq \epsilon\nrm{x-\other{x}}^2 + \epsilon^{-1}\nrm{y-\other{y}}^2,
            \numberthis\label{eq:prop:semi:young}
        \end{align*}
        \vspace{-\baselineskip}
        \begin{proofitemize}
            \item\ref{it:prop:semi:young:1}: The first inequality of \eqref{eq:prop:semi:young} implies that \eqref{eq:obliquequasisemimonotonicity:matrix} is always satisfied for $\mon \leq -\tfrac{\epsilon}{2}$ and $\com \leq -\tfrac{1}{2\epsilon}$, i.e., for $[\mon]_{-}[\com]_{-} \geq \tfrac{1}{4}$, which completes the proof.
            \item\ref{it:prop:semi:young:2}: The second inequality of \eqref{eq:prop:semi:young} implies that \eqref{eq:obliquequasisemimonotonicity:matrix} can never be satisfied when $\mon > \tfrac{\epsilon}{2}$ and $\com > \tfrac{1}{2\epsilon}$, i.e., when $[\mon]_+[\com]_+ > \tfrac{1}{4}$.
            \item\ref{it:prop:semi:young:3}: The second inequality of \eqref{eq:prop:semi:young} implies that \eqref{eq:obliquequasisemimonotonicity:matrix} can only be satisfied with equality when $\mon = \tfrac{\epsilon}{2} > 0$ and $\com = \tfrac{1}{2\epsilon} > 0$, i.e.,
            \ifspringer
                when $[\mon]_+[\com]_+ = \tfrac{1}{4}$. Hence, in this case it holds for all $(x, y),(\other{x}, \other{y}) \in \graph A$ that
                \(
                    \inner*{x - \other{x}, y - \other{y}} = \mon\nrm{x - \other{x}} + \tfrac{1}{4\mon}\nrm{y - \other{y}}^2,
                \)
            \else
                \begin{equation*}
                    \inner*{x - \other{x}, y - \other{y}} = \mon\nrm{x - \other{x}} + \tfrac{1}{4\mon}\nrm{y - \other{y}}^2, \qquad \text{for all $(x, y),(\other{x}, \other{y}) \in \gph A$},
                \end{equation*}
            \fi
            which is equivalent to
            \ifspringer
                \(
                    \mon\nrm{x - \other{x} - \tfrac{1}{2\mon}(y - \other{y})}^2 = 0.
                \)
            \else
                \begin{equation*}
                    \mon\norm{x - \other{x} - \tfrac{1}{2\mon}(y - \other{y})}^2 = 0, \qquad \text{for all $(x, y),(\other{x}, \other{y}) \in \gph A$}.
                \end{equation*}
            \fi
            By construction, this may only hold if $\graph A \subseteq \graph T$ for some $c \in \R^n$.
            Since $\gph T$ cannot be enlarged by another point $(\hat{x}, \hat{y})$ while maintaining \((\mon, \com)\)\hyp{}semimonotonicity, $T$ is maximally \((\mon, \com)\)\hyp{}semimonotone, and the proof is complete.
            \qedhere
        \end{proofitemize}
\end{appendixproof}

\begin{appendixproof}{it:prop:semi:lin:2}
    From \Cref{it:prop:semi:lin:1} and using that $A = \alpha \I_n$, it follows that $A$ is \(\left(\mon, \com\right)\)\hyp{}semimonotone if and only if \(\alpha \geq \mon + \com \alpha^2\), which is trivially satisfied when \(\left(\mon, \com\right) = \left(\alpha(1-c\alpha), c\right)\) for any $c \in \R$.
    Now, let $(\hat x, \hat y) \in \R^n \times \R^n$ and suppose that
    \begin{align*}
        \inner{x-\hat x, \alpha x - \hat y} \geq \alpha(1-c\alpha)\nrm{x - \hat x}^2 + c\nrm{\alpha x - \hat y}^2, \qquad \forall x \in \R^n.\numberthis\label{eq:id:maximal}
    \end{align*}
    To show when $A$ is maximally \(\left(\alpha(1-c\alpha), c\right)\)\hyp{}semimonotone, we must show when it is guaranteed that $(\hat x, \hat y) \in \graph A$, i.e., that $\hat y = \alpha \hat x$.
    Through elementary algebra, it can be shown that \eqref{eq:id:maximal} is equivalent to
    \begin{align*}
        0 \geq (1-2c\alpha)\inner{x,\hat y -\alpha \hat x} + \alpha (1-c\alpha)\nrm{\hat x}^2 - \inner{\hat x, \hat y} + c \nrm{\hat y}^2, \qquad \forall x \in \R^n.\numberthis\label{eq:id:maximal:cond}
    \end{align*}
    Due to the dependence of the first term on $x$, necessarily either $\hat y = \alpha \hat x$ or $c = \tfrac{1}{2\alpha}$. Plugging in $\hat y = \alpha \hat x$ into \eqref{eq:id:maximal:cond}, we get $0 \geq 0$ (as expected). Plugging in $c = \tfrac{1}{2\alpha}$ into \eqref{eq:id:maximal:cond}, we get
    \begin{align*}
        0 
            {}\geq{}
        \tfrac12\alpha\nrm{\hat x}^2 - \inner{\hat x, \hat y} + \tfrac{1}{2\alpha} \nrm{\hat y}^2
            {}={}
        \tfrac{1}{2\alpha}\nrm{\alpha \hat x - \hat y}^2 
        , \qquad \forall x \in \R^n.
    \end{align*}
    Therefore, \eqref{eq:id:maximal:cond} holds iff either (i) $\hat y = \alpha \hat x$ or (ii) $c = \tfrac{1}{2\alpha}$ and $\alpha < 0$, completing the proof.
\end{appendixproof}

\begin{appendixproof}{lem:sum}
    \begin{proofitemize}
        \item \ref{lem:sum:1}: First, consider the case where $A$ and $B$ are only semimonotone at a single point.
        Let $(x, y_A) \in \graph A$ and $(x, y_B) \in \graph B$. Then,
        it holds that
        \begin{align*}
            \inner{x-\other{x},y_A + y_B-(\other{y}_A+\other{y}_B)} 
                {}={}&
            \inner{x-\other{x}, y_A-\other{y}_A} + \inner{x-\other{x}, y_B-\other{y}_B}\\
                {}\geq{}&
            (\monA+\monB) \nrm{x-\other{x}}^2 + \comA\nrm{y_A-\other{y}_A}^2 + \comB\nrm{y_B-\other{y}_B}^2.
        \end{align*}
        Consequently, by definition $A+B$ is $(\monA+\monB,\com)$\hyp{}semimonotone at $(\other{x}, \other{y}_A+\other{y}_B)$ if and only if
        \ifspringer
            \begin{align*}
                \com\nrm{y_A + y_B - (\other{y}_A + \other{y}_B)}^2 
                    {}={}&
                \com\nrm{y_A-\other{y}_A}^2 + 2\com\inner{y_A-\other{y}_A,y_B-\other{y}_B} + \com\nrm{y_B-\other{y}_B}^2\\
                    {}\leq{}&
                \comA\nrm{y_A-\other{y}_A}^2 + \comB\nrm{y_B-\other{y}_B}^2,
            \end{align*}
        \else
            \begin{align*}
                \com\nrm{y_A + y_B - (\other{y}_A + \other{y}_B)}^2 
                    {}={}&
                \com\nrm{y_A-\other{y}_A}^2 + 2\com\inner{y_A-\other{y}_A,y_B-\other{y}_B} + \com\nrm{y_B-\other{y}_B}^2
                    {}\leq{}
                \comA\nrm{y_A-\other{y}_A}^2 + \comB\nrm{y_B-\other{y}_B}^2,
            \end{align*}
        \fi
        which is equivalent to the LMI 
        $
            \left(\begin{smallmatrix} 
                \comA - \com & -|\com|\\
                -|\com| & \comB - \com
            \end{smallmatrix}\right) \succeq 0.
        $
        This is satisfied if and only if the determinant and trace are nonnegative, i.e.,
        \begin{align*}
            \comA\comB - \com(\comA+\comB) \geq 0
            \;\; \text{and} \;\;
            \comA + \comB - 2\com \geq 0
            \quad\Leftrightarrow\quad 
            \com \leq 
            \begin{cases}
                \tfrac{\comA\comB}{\comA+\comB}, & \quad \text{if } \comA + \comB > 0,\\
                0, & \quad \text{if } \comA = \comB = 0.
            \end{cases}\numberthis\label{eq:proof:lem:sum:suff}
        \end{align*}
        Since by definition the largest possible value for $\com$ satisfying \eqref{eq:proof:lem:sum:suff} is given by $\comA\Box\comB$ (if either $\comA+\comB > 0$ or $\comA=\comB=0$), the first part of the proof is completed.

        If $A$ and $B$ are semimonotone then by definition $A$ is $(\monA,\comA)$\hyp{}semimonotone at all $(\other{x}_A,\other{y}_A) \in \graph A$ and $B$ is $(\monB,\comB)$\hyp{}semimonotone at all $(\other{x}_B,\other{y}_B) \in \graph B$. Hence, using the previous result we know that $A + B$ is $(\monA+\monB,\comA \Box \comB)$\hyp{}semimonotone at all points in the set 
        \ifspringer
            \(
                \set{(\other{x}, \other{v}_A + \other{v}_B)}[\other{x} \in \dom A \cap \dom B,\; \other{v}_A \in A\other{x},\; \other{v}_B \in B\other{x}].
            \)
        \else
            \[
                \set{(\other{x}, \other{v}_A + \other{v}_B)}[\other{x} \in \dom A \cap \dom B,\; \other{v}_A \in A\other{x},\; \other{v}_B \in B\other{x}].
            \]
        \fi
        Since this set is equal to $\graph (A+B)$, the proof is completed.
        \item \ref{lem:sum:2}: This is an immediate consequence of the definition of the parallel sum, owing to \Cref{prop:calculus:inverse} and \Cref{lem:sum:1}.
        \qedhere 
    \end{proofitemize}
\end{appendixproof}

\begin{appendixproof}{cor:sum:identity}
    \begin{proofitemize}
        \item \ref{it:cor:sum:identity:1}: First, consider the assertion where the semimonotonicity holds only at a single point.
        Operator $A$ is $(\monA, \comA)$\hyp{}semimonotone at $(\other{x}, \other{y}) \in \graph A$ if and only if \(\inner{x-\other{x},y-\other{y}} \geq \monA \nrm{x-\other{x}}^2 + \comA \nrm{y-\other{y}}^2\) for all $(x,y) \in \graph A$. Defining $p \coloneqq y + \alpha x$ and $\other{p} \coloneqq \other{y} + \alpha \other{x}$, 
        so that $(x,p), (\other{x}, \other{p}) \in \graph T$,        
        this is equivalent to
        \begin{align*}
            (1 + 2\comA\alpha)\inner{x - \other{x},p-\other{p}} 
                {}\geq{}
            (\monA + \alpha(1 + \comA\alpha))\nrm{x - \other{x}}^2 + \comA\nrm{p - \other{p}}^2, \quad \forall (x, p) \in \graph T.
        \end{align*} 
        Dividing both sides by $1 + 2\comA\alpha > 0$, 
        the claimed equivalence between $(\monA, \comA)$\hyp{}semimonotonicity of $A$ at $(\other{x}, \other{y})$ and $\left(\tfrac{\monA + \alpha(1 + \comA\alpha)}{1 + 2\comA\alpha}, \tfrac{\comA}{1 + 2\comA\alpha}\right)$\hyp{}semimonotonicity of $T$ at $(\other{x}, \other{y} + \alpha\other{x})$ is established.
        The assertion for (global) semimonotonicity follows analogously by considering all points $(\other{x}, \other{y}) \in \graph A$.
        
        Finally, the assertion for maximality holds since \(\graph A\) can be enlarged with a point \((\hat x,\hat y)\) without breaking \((\monA, \comA)\)\hyp{}semimonotonicity if and only if \(\graph T\) can be enlarged with \((\hat x,\hat y + \alpha \hat x)\) without violating 
        $\left(\tfrac{\monA + \alpha(1 + \comA\alpha)}{1 + 2\comA\alpha}, \tfrac{\comA}{1 + 2\comA\alpha}\right)$\hyp{}semimonotonicity.
        \item \ref{it:cor:sum:identity:2}: This result immediately follows from \Cref{lem:sum:1} and \Cref{it:prop:semi:lin:2}.
        \qedhere
    \end{proofitemize}
\end{appendixproof}

\begin{appendixproof}{prop:semimon:osc}
    First, assume that $\com\mon < \tfrac{1}{4}$ and 
    let $\xi$ and $\nu$ be defined as in \eqref{eq:XiBeta}.
    Let $(x^k, y^k) \in \graph A$, where \(x^k \rightarrow \bar{x}\) and $y^k \rightarrow \bar{y}$.
    Then, it holds by construction that $x^k - \xi(y^k - \nu x^k) \in \Bigl[(A^{-1} - \nu \id)^{-1} - \xi\id\Bigr](y^k - \nu x^k)$.
    Owing to \Cref{cor:generation semimonotone operator},
    $(A^{-1} - \nu \id)^{-1} - \xi\id$ is maximally monotone and thus also outer semicontinuous \cite[Ex. 12.8]{rockafellar2009Variational}. Consequently, it follows that $\bar{x} - \xi(\bar{y} - \nu \bar{x}) \in \Bigl[(A^{-1} - \nu \id)^{-1} - \xi\id\Bigr](\bar{y} - \nu \bar{x})$, which in turn implies that \(\bar{x} \in A \bar{y}\), showing outer semicontinuity of \(A\). Finally, in the case where $[\mon]_+[\com]_+ = \nicefrac14$, maximal monotonicity of $A$ follows by combining \Cref{it:prop:semi:young:3:maximal,it:prop:semi:lin:2,prop:calculus:scaling}, which implies its outer semicontinuity \cite[Ex. 12.8]{rockafellar2009Variational}.
\end{appendixproof}

\begin{appendixproof}{cor:semimon:fulldom}    
    \begin{proofitemize}
      \item \ref{cor:semimon:resolvent}:
        First, assume that $\mon\com < \tfrac14$. Define $\theta \coloneqq \sqrt{1-4\com\mon}$ and let \(\nu \coloneqq \tfrac{2\mon}{1+\theta}\), \(\xi \coloneqq \tfrac{\com}{\theta}\) and 
        \(\omega \coloneqq \tfrac{\gamma}{1+\gamma\nu}\).
        Using the upper bound $\gamma < \tfrac{1+\theta}{2[-\mon]_{+}}$ from \eqref{eq:cor:fulldom:gamma}, it follows that
        \(
        1+ \gamma \nu
            {}={}
        1+ \gamma \tfrac{2\mon}{1+\theta}
            {}\geq{}
        1 - \tfrac{2\gamma[-\mon]_+}{1+\theta}
            {}>{}
        0
        \).
        Consequently, \(\omega + \xi > 0\) since
        \begin{align*}
            \tfrac{\gamma}{1+\gamma\nu} + \tfrac{\com}{\theta} > 0
                {}\stackrel{1+\gamma\nu \,>\, 0}{\Longleftrightarrow}{}
            \gamma + \Bigl(1+\gamma\tfrac{2\mon}{1+\theta}\Bigr)\tfrac{\com}{\theta} > 0
                {}\stackrel{\theta\,>\,0}{\Longleftrightarrow}{}
            \Bigl(2\theta + \tfrac{4\mon\com}{1+\theta}\Bigr)\gamma + 2\com > 0
                {}\Longleftrightarrow{}
            (1 + \theta)\gamma + 2\com > 0
        \end{align*}
        which is satisfied by the lower bound $\gamma > \tfrac{2[-\com]_{+}}{1+\theta}$ from 
        \iftables
            \Cref{tab:semimon:domain}.
        \else
            \eqref{eq:cor:fulldom:gamma}.
        \fi
        Denote 
        \(
            \hat A \coloneqq (A-\nu \id)^{-1} + \omega \id
        \). 
        Owing to \(\omega + \xi > 0\) and Minty's theorem \cite[Thm. 21.1]{bauschke2017Convex}, the operator 
        \(
            \tfrac{1}{\omega+\xi}\hat A 
                {}={}
            \tfrac{1}{\omega+\xi}
            \left(%
                (A-\nu\id)^{-1} - \xi\id
            \right)
                {}+{}
            \id
        \)
        has full range iff \(\tfrac{1}{\omega+\xi}\big((A-\nu\id)^{-1} - \xi\id\big)\) is maximally monotone.
        Since \(\omega+\xi>0\), this implies that \(\hat A\) has full range iff the operator \(\tilde A \coloneqq(A-\nu\id)^{-1} - \xi\id\) is maximally monotone. 

        On the other hand, by definition 
         \(
            1 + 2 \xi \nu 
            {}={}
            1+ \tfrac{1-\theta^2}{\theta(1+\theta)}
            {}={}
            \tfrac{1}{\theta} 
            {}>{}
            0.
         \)
         Hence, using \cref{it:cor:gen:1}, the operator  $\tilde A\coloneqq(A-\nu\id)^{-1} - \xi\id$ is maximally monotone if and only if $A$ is maximally $\Bigl(\tfrac{\nu(1 + \xi\nu)}{1 + 2\xi\nu}, \tfrac{\xi}{1 + 2 \xi\nu}\Bigr) = (\mu, \rho)$ semimonotone, where the equality is the result of substituting \(\xi\) and \(\nu\). 
         Therefore, it suffices to show \(\hat A\) has full range iff  \(\range{\id+\gamma A}=\R^n\), \ie, that the resolvent has full domain. First suppose that \(\hat A\) has full range. Fix \(y\in\R^n\). Then, by the full range assumption for \(\hat A\), there exists some \(x\in\R^n\) such that 
        \(
            \tfrac{\omega}{\gamma} y\in (A - \nu \id)^{-1}x+\omega x.
        \) 
        Letting \(z = \omega(\tfrac{1}{\gamma}y-x)\) we have 
        \begin{align*}
            \tfrac{1}{\gamma} y - \tfrac1{\omega}z \in (A - \nu \id)z  
                {}\iff{}
            \tfrac{1}{\gamma} y \in \big(A + (\tfrac1{\omega} - \nu) \id\big)z      
                {}\iff{}
            y \in \big(\gamma A +  \id\big)z,    
            \numberthis\label{eq:semimon:resolvent:1}  
        \end{align*}
        where in the last equivalence we used the fact that \(\gamma = \tfrac{\omega}{1-\omega \nu}\). Since the choice of \(y\) was arbitrary, \(J_{\gamma A}\) has full domain. 
         The converse follows similarly. Assume that  \(\range{\I+\gamma A}=\R^n\). For any fixed \(y\in\R^n\) there exists some \(z\in\R^n\) such that \(y\in \gamma Az + z\). Therefore, by \eqref{eq:semimon:resolvent:1} and letting \(x = \tfrac{1}{\gamma}y - \tfrac{1}{\omega}z\) we obtain \(\tfrac{\omega}{\gamma}y\in\hat Ax\), thus establishing that \(\range{\hat A}=\R^n\).          
        

        Finally, consider the case \([\mon]_+[\com]_+ = \tfrac14\). 
        The range for the stepsize in \eqref{eq:cor:fulldom:gamma} then reduces to $\gamma \in (0, +\infty)$. 
        First, suppose that $J_{\gamma A}$ has full domain, \ie \(\range{\id + \gamma A} = \R^n\). Fix \(v\in \R^n\). Then, there exists \((x,y)\in \graph A\) such that \(v = x + \gamma y\). By  \Cref{it:prop:semi:young:3:regular}, \(y = 2\mu x + c\) for some \(c\in \R\)  and some \(x\in \dom A\), and therefore, \(v = (1+ 2\gamma \mu) x +c \gamma\). Since \(1+ 2\gamma \mu>0\) and the choice of \(v\) was arbitrary, necessarily \(\dom A = \R^n\), and as such \(A\) is maximally \((\mu,\rho)\)-semimonotone, see \Cref{it:prop:semi:young:3:maximal}. 
        Now, suppose that $A$ is maximally $(\mon, \com)$\hyp{}semimonotone. Then, by \Cref{it:prop:semi:young:3:maximal} it holds that \(A : x \mapsto 2\mon x + c\) for some $c \in \R^n$, which implies that $A$ is maximally monotone. Therefore by Minty's theorem \cite[Thm. 21.1]{bauschke2017Convex} the resolvent $J_{\gamma A}$ has full domain.  

    \item \ref{cor:semimon:Lip}: 
          Take $(u,v), (\other{u}, \other{v}) \in \graph J_{\gamma A}$. By semimonotonicity of \(A\) and using the fact that  
    \(
        \tfrac1{\gamma}(u-v) \in Av,
    \) 
    \(
        \tfrac1{\gamma}(\other{u}-\other{v}) \in A\other{v}
    \), we have 
    \begin{align*}
        \mon \nrm{v-\other{v}}^2 + \com \gamma^{-2} \nrm{(u-v)-(\other{u}-\other{v})}^2 &\leq \gamma^{-1}\inner{(u-v)-(\other{u}-\other{v}), v-\other{v}}.
    \end{align*}
    Expanding the last term of the left-hand side and rearranging, we get
    \begin{equation}\label{eq:cor:semimon:inequality}
        \begin{aligned}
            \left(1 + \mon\gamma + \com \gamma^{-1}\right)\nrm{v-\other{v}}^2 + \com \gamma^{-1}\nrm{u-\other{u}}^2 &\leq \left(1 + 2 \com \gamma^{-1}\right)\inner{u-\other{u},v-\other{v}}\\
            & \leq \left|1 + 2 \com \gamma^{-1}\right|\nrm{u-\other{u}}\nrm{v-\other{v}}.
        \end{aligned}
    \end{equation}
    By \Cref{it:lem:quadratic:solutions} it follows that
    \iftables
        \Cref{tab:semimon:domain}
    \else
        \eqref{eq:cor:fulldom:gamma}
    \fi
    provides all strictly positive solutions for $\gamma$ of \(1 + \mon\gamma + \com \gamma^{-1} > 0\). Hence, if \(u=\other{u}\) then necessarily \(v=\other{v}\), and the Lipschitz inequality holds trivially. Else, divide both sides by \(\|u-\other{u}\|^2\) to obtain 
    \ifspringer
        \(
            \left(1 + \mon\gamma + \com \gamma^{-1}\right)\omega^2 - \left|1 + 2 \com \gamma^{-1}\right| \omega + \com \gamma^{-1}\leq0,
        \)
    \else
        \begin{equation*}
            \left(1 + \mon\gamma + \com \gamma^{-1}\right)\omega^2 - \left|1 + 2 \com \gamma^{-1}\right| \omega + \com \gamma^{-1} \leq 0,
        \end{equation*}
    \fi
    where  
        \(
        \omega \coloneqq \tfrac{\nrm{v-\other{v}}}{\nrm{u-\other{u}}}    
    \).
    Solving this inequality for $\omega$ yields the claimed Lipschitz modulus.
    \qedhere
    \end{proofitemize}
\end{appendixproof}

    \subsection{Examples of functions with semimonotone subdifferentials}\label{subsec:auxiliary:functions}

\begin{appendixproof}{prop:curvature}
    \begin{proofitemize}
        \item \ref{prop:hypo:cvx}: Follows by definition.  
        \item \ref{prop:hypo:ccv}: 
        By \cite[Prop. 12.60]{rockafellar2009Variational} and strong convexity of $-f$, the conjugate $(-f)^*$ is $-\tfrac{1}{\ell}$\hyp{}smooth, i.e.
        \begin{align*}
            \inner{\nabla (-f)^*(x) - \nabla (-f)^*(\other{x}), x - \other{x}} \leq -\tfrac{1}{\ell} \nrm{x - \other{x}}^2, \quad \forall x, \other{x} \in \R^n.\numberthis\label{eq:cor:strcvx:1}
        \end{align*}
        Moreover, since $-f$ is continuous and convex, we have by 
        \cite[Prop. 11.3]{rockafellar2009Variational} that 
        \(
            \nabla (-f)^* = (\partial (-f))^{-1}
        \),
        so that \eqref{eq:cor:strcvx:1} is equivalent to
        \begin{align*}
            \inner{y - \other{y}, x - \other{x}} \leq -\tfrac{1}{\ell} \nrm{y - \other{y}}^2, \quad \forall (x, y), (\other{x}, \other{y}) \in \graph \partial (-f).\numberthis\label{eq:cor:strcvx:2}
        \end{align*}
        Finally, since $-f$ is finite and convex, it follows from \cite[Cor. 9.21]{rockafellar2009Variational} that \(\partial f(\bar{x}) \subseteq - \partial (-f)(\bar{x})\), \(\forall \bar{x} \in \R^n\), so that
        \begin{align*}
            \inner{y - \other{y}, x - \other{x}} \geq \tfrac{1}{\ell} \nrm{y - \other{y}}^2, \quad \forall (x, y), (\other{x}, \other{y}) \in \graph \partial f \subseteq - \graph \partial (-f).
        \end{align*}

        \item \ref{prop:ccv}: Consider the function 
        \(\varphi = f - \tfrac{\alpha}2\|\cdot\|^2\). Then, by definition \(\varphi\in\mathcal F^{\ell - \alpha}(\R^n)\), implying that \(\partial \varphi\) is \((0, \tfrac{1}{\ell - \alpha})\)-semimonotone owing to \Cref{prop:hypo:ccv} since \(\alpha>\ell\). Since \(\partial f= \partial \varphi + \alpha \id\), \cite[Ex. 8.8(c)]{Rockafellar1970Convex},  it follows from \Cref{cor:sum:identity} that for \(c > \tfrac{1}{\alpha - \ell}\)  the claimed semimonotonicity holds. 

        \item \ref{prop:hypo:dif+}, \ref{prop:hypo:dif-}:
        Suppose that \(\ell>\sigma\), since if \(\ell = \sigma\), then, \(\nabla f = \ell \id\) and the claims follow trivially.
        Consider the function \(\varphi = f - \tfrac{\sigma}2\|\cdot\|^2 \in \mathcal F^{\ell - \sigma}_0(\R^n)\).
        Since $\varphi$ is finite and convex its subdifferential 
        $\partial \varphi(x)$ is nonempty for all $x \in \R^n$ \cite[Thm. 23.4]{Rockafellar1970Convex}. Consequently, $\varphi$ is continuously differentiable due to \cite[Lem. 2.1]{themelis2020DouglasRachford} and $\nabla \varphi$ is \((0, \tfrac{1}{\ell-\sigma})\)\hyp{}semimonotone due to \cite[Thm. 2.1.10]{nesterov2004introductory}.
        Therefore, by \cref{prop:calculus:inverse}, there exists a monotone mapping \(\tilde A\) for which
        \(
        (\nabla \varphi)^{-1} = \tilde A + \tfrac{1}{\ell-\sigma}
        \) and thus, owing to $\nabla f = \nabla \varphi + \sigma \id$ \cite[Ex. 8.8(c)]{rockafellar2009Variational} it holds that 
        \ifspringer
            \(
                \nabla f 
                {}={}
                \nabla \varphi + \sigma \id 
                {}={}
                \big(\tilde A + \tfrac{1}{\ell-\sigma}\id \big)^{-1} + \sigma \id.
            \)
        \else
            \[
                \nabla f 
                {}={}
                \nabla \varphi + \sigma \id 
                {}={}
                \big(\tilde A + \tfrac{1}{\ell-\sigma}\id \big)^{-1} + \sigma \id.
            \]
        \fi
        The claimed results from \Cref{prop:hypo:dif+,prop:hypo:dif-} then follow from \cref{cor:sem2mon}, where we let $\xi = \tfrac{1}{\ell - \sigma}$ and $\nu = \sigma$. This is because \(\ell - \sigma > 0\) and therefore $1 + 2\xi\nu > 0$ (resp. $\leq 0$) if and only if $\sigma + \ell > 0$ (resp. $\leq 0$).
        \qedhere
    \end{proofitemize}
\end{appendixproof}

\begin{appendixproof}{lem:pointwise:semi}
    Since \(f_i \in \mathcal F^{\ell_{f_i}}(\R^n)\), it follows that 
    \(\tfrac{\bar \ell}2\|\cdot\|^2 - f_i\) is convex. Therefore, $\tfrac{\bar \ell}2\|x\|^2 - f(x) = \max\left\{\tfrac{\bar \ell}2\|x\|^2 - f_1(x),\hdots,\tfrac{\bar \ell}2\|x\|^2 - f_N(x)\right\}$ is convex being the pointwise max of convex functions and thus $f \in \mathcal F^{\bar \ell}(\R^n)$. The claimed result then follows immediately from \Cref{prop:hypo:ccv,prop:ccv}.
\end{appendixproof}

\begin{appendixproof}{lem:infimal}
    \begin{proofitemize} 
        \item \ref{it:infimal:subdif}: Let \(f(w) = f_1(w_1) + f_2(w_2)\), and \(L\) be the linear mapping \(Lw = w_1+w_2\). 
        Then, the infimal convolution may be viewed as the parametric minimization of \(\tilde f(w,s) = f(w) + \delta_{\{0\}}(Lw-s)\) with respect to the variable \(w\), i.e., \(\varphi(s) = \min_w \tilde f(w,s)\).
        Let 
        \(
            M(s) \coloneqq \argmin_{(w_1,w_2)} \set{f_1(w_1) + f_2(w_2)}[w_1+w_2=s].
        \) 
        Then, owing to \eqref{eq:lem:infimal:uniform-level-bounded} it holds by \cite[Thm. 10.13]{rockafellar2009Variational} for all \(\bar s \in \R^n\) that
        \begin{align*}
            \partial \varphi(\bar s) 
                {}\subseteq{}
            \bigcup_{\bar x \in M(\bar s)} \set{y\in R^n}[(0,y)\in \partial \tilde f(\bar x,\bar s)].
            \numberthis\label{eq:subdif:epicompos:proof}
        \end{align*}
        Consider \(\bar f (w,u) = f(w) + \delta_{\{0\}}(u)\) and the linear transformation \(K:(w,s)\mapsto(w, Lw-s)\). Then, 
        \(\tilde f = \bar f \circ K\). Since \(K\) has full rank, by \cite[Ex. 10.7]{rockafellar2009Variational} for any pair \((\bar x,\bar s)\) such that \(\bar x\in M(\bar s)\) it holds that
        \ifspringer
            \(
                \partial \tilde f(\bar x, \bar s) 
                    {}={}
                (K^\top \circ \partial \bar f \circ K) (\bar x, \bar s)
            \),
            where
            \(
                \partial \bar f (w,u) = (\partial f(w), N_{\{0\}}(u))
            \).
        \else 
            \begin{align*}
                \partial \tilde f(\bar x, \bar s) 
                    {}={}&
                (K^\top \circ \partial \bar f \circ K) (\bar x, \bar s)
                \quad\text{where}\quad
                \partial \bar f (w,u) = (\partial f(w), N_{\{0\}}(u)).
            \end{align*}
        \fi
        Using this in \eqref{eq:subdif:epicompos:proof} yields
        \ifspringer
            \begin{align*}
                \partial \varphi(\bar s) 
                    {}\subseteq{}&
                \bigcup_{\bar x\in M(\bar s)} \set{y\in R^n}[(0,y)\in \left(\partial f(\bar x) + \tp L N_{\{0\}}(L\bar x - \bar s), - N_{\{0\}}(L\bar x - \bar s)\right)]\\
                    {}={}&
                \bigcup_{(\bar x_1,\bar x_2)\in M(\bar s)} \partial f_1(\bar x_1) \cap \partial f_2(\bar x_2).
            \end{align*}
        \else
            \begin{align*}
                \partial \varphi(\bar s) 
                    {}\subseteq{}&
                \bigcup_{\bar x\in M(\bar s)} \set{y\in R^n}[(0,y)\in \left(\partial f(\bar x) + \tp L N_{\{0\}}(L\bar x - \bar s), - N_{\{0\}}(L\bar x - \bar s)\right)]
                    {}={}
                \bigcup_{(\bar x_1,\bar x_2)\in M(\bar s)} \partial f_1(\bar x_1) \cap \partial f_2(\bar x_2).
            \end{align*}
        \fi
        Hence, the claimed result follows from \cite[Prop. 25.30(i)]{bauschke2017Convex}.
        \item \ref{it:infimal:semi}: Follows immediately from \Cref{it:infimal:subdif} and \Cref{lem:sum:2}.
        \qedhere
    \end{proofitemize}
\end{appendixproof}

\begin{appendixproof}{prop:saddleOp}
    Consider the saddle envelope 
    \begin{equation*}
        \varphi_\delta(x,y) \coloneqq \min_{u\in \R^n}\max_{v\in \R^m} \big\{\varphi(u,v) + \tfrac1{2\delta}\|u-x\|^2- \tfrac1{2\delta}\|v-y\|^2\big\},     
    \end{equation*}
    along with the associated resolvent $J_{\delta T_{\varphi}} \coloneqq (\id + \delta T_{\varphi})^{-1}$. For any  \((z,{w})\in \graph J_{\delta T_{\varphi}}\) with the decomposition \(z= (x,y)\) and \(w = (u,v)\) we have by construction that
    \ifspringer
        \(
            0 
                {}={}
            \nabla_x \varphi({w}) + \tfrac1{\delta}(u-x)
                {}={}
            \nabla_x \varphi({w}) - \nabla_x \varphi_{\delta}(z)
        \)
        and
        \(
            0
                {}={}
            \nabla_y \varphi({w}) - \tfrac1{\delta}(v-y)
                {}={}
            \nabla_y \varphi({w}) - \nabla_y \varphi_{\delta}(z)
        \).
    \else
        \begin{align*}
            0 
            {}={}&
            \nabla_x \varphi({w}) + \tfrac1{\delta}(u-x)
            {}={}
            \nabla_x \varphi({w}) - \nabla_x \varphi_{\delta}(z)
            \quad \text{and} \quad
            0
            {}={}
            \nabla_y \varphi({w}) - \tfrac1{\delta}(v-y)
            {}={}
            \nabla_y \varphi({w}) - \nabla_y \varphi_{\delta}(z).
        \end{align*}
    \fi
    Hence,
    \begin{equation}\label{eq:optimality:grad}
    T_{\varphi}({w}) = \big(\nabla_x \varphi({w}), -\nabla_y \varphi({w})\big)= \big(\nabla_x \varphi_\delta(z), -\nabla_y \varphi_\delta(z)\big) \eqqcolon T_{\varphi_\delta}(z).
    \end{equation}
    On the other hand, it follows from \cite[Prop. 1]{grimmer2022landscape} that 
    \begin{equation}\label{eq:mon:F_env}
    \tfrac{1}{\delta} \I \succeq \nabla_{xx}\varphi_{\delta}(z)\succeq \tfrac{\alpha}{1+\delta\alpha}\I,
    \quad 
    \tfrac{1}{\delta} \I \succeq -\nabla_{yy}\varphi_{\delta}(z)\succeq \tfrac{\alpha}{1+\delta\alpha}\I, \quad \forall z,
    \end{equation}
    thus implying that the operator $T_{\varphi_\delta}$ is $\tfrac{\alpha}{1+\delta\alpha}$-monotone. 
    Let $w_1,w_2\in\R^n$ and $z_{i}=w_{i} + \delta T_{\varphi}w_i$ for $i=1,2$. Then, \((z_i,w_i)\in \graph J_{\delta T_{\varphi}}\) for \(i=1,2\) and we have that
    \begin{align*}
        \langle T_{\varphi}{w_1}-T_{\varphi}{w_2},{w_1}-{w_2}\rangle 
            {}={}&
        \langle T_{\varphi}{w_1}-T_{\varphi}{w_2},z_{1}-z_{2}\rangle
            {}-{}
        \delta\|T_{\varphi}{w_1}-T_{\varphi}{w_2}\|^{2}
        \\%
        \dueto{\eqref{eq:optimality:grad}}
            {}={}&
        \langle T_{\varphi_{\delta}}z_{1}-T_{\varphi_{\delta}}z_{2},z_{1}-z_{2}\rangle-\delta\|T_{\varphi}{w_1}-T_{\varphi}{w_2}\|^{2}
        \\%
        \dueto{\eqref{eq:mon:F_env}}
            {}\geq{}&
        \tfrac{\alpha}{1+\alpha\delta}\|z_{1}-z_{2}\|^{2}-\delta\|T_{\varphi}{w_1}-T_{\varphi}{w_2}\|^{2}.\numberthis\label{eq:cohy}
    \end{align*}
    The second assertion follows from the above inequality noting that \(1+\alpha\delta>0\). 
    Moreover, plugging $z_{i}=w_{i} + \delta T_{\varphi}w_i$ in \eqref{eq:cohy} and multiplying by $(1+\alpha \delta)$ yields 
    \ifspringer
        \(
            \left(1-\alpha\delta\right)\langle T_{\varphi}w_1- T_{\varphi}w_2,w_1-w_2\rangle 
                {}\geq{}
            {\alpha}\|{w_1}-{w_2}\|^{2}-\delta\|T_{\varphi}{w_1}-T_{\varphi}{w_2}\|^{2},
        \)
    \else
        \begin{align*}
            \left(1-\alpha\delta\right)\langle T_{\varphi}w_1- T_{\varphi}w_2,w_1-w_2\rangle 
                {}\geq{}
            {\alpha}\|{w_1}-{w_2}\|^{2}-\delta\|T_{\varphi}{w_1}-T_{\varphi}{w_2}\|^{2},
        \end{align*}
    \fi
    establishing the first assertion. 
\end{appendixproof}

    \subsection{Douglas-Rachford splitting for semimonotone operators}

\begin{appendixproof}{prop:semi:WMVI}
    Let $(x_\A, y_\A) \in \graph \A$ and $(x_\B, y_\B) \in \graph \B$.
        Additionally, if \Cref{ass:SWMVIDRSsemimonotone:maximal} holds, let $(\other{x}_\A, \other{y}_\A) \in \graph \A$ and $(\other{x}_\B, \other{y}_\B) \in \graph \B$.
        Otherwise, if \Cref{ass:SWMVIDRSsemimonotone} holds, let $(x^\star, y^\star) \in \pazocal{S}^\star$ and define $(\other{x}_\A, \other{y}_\A) \coloneqq (x^\star,-y^\star) \in \graph A$
        and $(\other{x}_\B, \other{y}_\B) \coloneqq (x^\star,y^\star) \in \graph B$.
    Then, by definition of the primal-dual operator from \eqref{eq:primaldual}, it holds that 
        \ifspringer
            \(
                (z, v) \coloneqq \Bigl((x_\A,y_\B), (y_\A + y_\B, x_\B - x_\A)\Bigr) \in \graph \Tpd
            \)
            and
            \(
                (\other{z}, \other{v}) \coloneqq \Bigl((\other{x}_\A,\other{y}_\B), (\other{y}_\A + \other{y}_\B, \other{x}_\B - \other{x}_\A)\Bigr) \in \graph \Tpd.
            \)
        \else
            \[
                (z, v) \coloneqq \Bigl((x_\A,y_\B), (y_\A + y_\B, x_\B - x_\A)\Bigr) \in \gph \Tpd
                \quad\text{and}\quad
                (\other{z}, \other{v}) \coloneqq \Bigl((\other{x}_\A,\other{y}_\B), (\other{y}_\A + \other{y}_\B, \other{x}_\B - \other{x}_\A)\Bigr) \in \graph \Tpd.
            \]
        \fi
        Owing to the semimonotonicity assumptions and \Cref{prop:calculus:inverse}, this implies that
        \begin{align*}
        \inner*{v - \other{v}, z - \other{z}} 
            {}={}&
        \inner*{y_\A + y_\B - (\other{y}_\A + \other{y}_\B), x_\A - \other{x}_\A} + \inner*{x_\B - x_\A - (\other{x}_\B - \other{x}_\A), y_\B - \other{y}_\B}\\
            {}={}&
        \inner*{y_\A - \other{y}_\A, x_\A - \other{x}_\A} + \inner*{x_\B - \other{x}_\B, y_\B - \other{y}_\B}\\
            {}\geq{}&
        \monA \nrm{x_\A - \other{x}_\A}^2 + \comA \nrm{y_\A - \other{y}_\A}^2 + \comB \nrm{y_\B - \other{y}_\B}^2 + \monB \nrm{x_\B - \other{x}_\B}^2.\numberthis\label{eq:proof:semi:WMVI:1}
        \end{align*}
        On the other hand, it follows from the Cauchy-Schwarz inequality that
        \begin{align*}
            \qindef{v - \other{v}}{\DRSRho} 
                {}={}&
            \DRSrhocom\nrm{y_\A + y_\B - (\other{y}_\A + \other{y}_\B)}^2 + \DRSrhomon\nrm{x_\B - x_\A - (\other{x}_\B - \other{x}_\A)}^2\\
                {}\leq{}& 
            \begin{multlined}[t]
                \DRSrhocom\nrm{y_\A - \other{y}_\A}^2 + \DRSrhocom\nrm{y_\B - \other{y}_\B}^2 + 2|\DRSrhocom|\nrm{y_\A - \other{y}_\A}\nrm{y_\B - \other{y}_\B}\\
                + \DRSrhomon\nrm{x_\B - \other{x}_\B}^2 + \DRSrhomon\nrm{x_\A - \other{x}_\A}^2 + 2|\DRSrhomon|\nrm{x_\B - \other{x}_\B} \nrm{x_\A - \other{x}_\A}.
            \end{multlined}\numberthis\label{eq:proof:semi:WMVI:2}
        \end{align*}
        Combining \eqref{eq:proof:semi:WMVI:1} and \eqref{eq:proof:semi:WMVI:2}, it follows 
        that the primal-dual operator $\Tpd$ 
        has $\DRSRho$\hyp{}oblique weak Minty solutions at $\pazocal{S}^\star$ if \Cref{ass:SWMVIDRSsemimonotone} holds and that it
        is $\DRSRho$\hyp{}comonotone if \Cref{ass:SWMVIDRSsemimonotone:maximal} holds,
    provided that \(X \coloneqq \left(\begin{smallmatrix} \monA - \DRSrhomon & -|\DRSrhomon|\\ -|\DRSrhomon| & \monB - \DRSrhomon\end{smallmatrix}\right) \succeq 0\) and \(Y \coloneqq \left(\begin{smallmatrix} \comB - \DRSrhocom & -|\DRSrhocom|\\ -\DRSrhocom & \comA - |\DRSrhocom|\end{smallmatrix}\right) \succeq 0\).
    As previously discussed in the proof of \Cref{lem:sum:1}, the largest possible value for $\DRSrhomon$ satisfying $X \succeq 0$ is given by $\monA\Box\monB$ if either $\monA+\monB > 0$ or $\monA=\monB=0$ (and analogously the largest value for $\DRSrhocom$ satisfying $Y \succeq 0$ is given by $\comA\Box\comB$ if either $\comA+\comB > 0$ or $\comA=\comB=0$), which completes the proof
    for \Cref{it:prop:semi:WMVI}.
    To show \Cref{it:prop:semi:WMVI:max}, it only remains to show maximality.
        To this end, consider an arbitrary stepsize $\gamma$ complying with
        \Cref{tab:gamma:DRS:semi}, where $\gamma_{-}$ and $\gamma_{+}$ are defined as in \Cref{eq:eigcondDRS:gamma:semi}.
        Then, by \Cref{prop:SWMVIDRSsemimonotone:maximal:fulldom}
        the resolvents $J_{\gamma A}$ and $J_{\gamma B}$ have full domain, which in turn implies that 
        the preconditioned resolvent $(\M + \Tpd)^{-1}\M$ has full domain, where $\M$ is given by \eqref{eq:M} \cite[Lem. 12.14]{rockafellar2009Variational}.
        The claim then follows from \cref{prop:comon:fulldom}.
\end{appendixproof}

\begin{appendixproof}{cor:DRS:semi}
    Consider the following implications.
        \begin{proofitemize}
        \item For \Cref{it:cor:DRS:semi:1}, it follows from \Cref{ass:SWMVIDRSsemimonotone} and 
                \Cref{it:prop:semi:WMVI}
        that \Cref{ass:SWMVIDRS} holds.
        \item For \Cref{it:cor:DRS:semi:2}, it follows from \Cref{ass:SWMVIDRSsemimonotone:maximal} and \Cref{prop:SWMVIDRSsemimonotone:maximal:fulldom} that 
                \Cref{ass:DRS:1,ass:DRS:2} hold
        and from 
        \Cref{ass:SWMVIDRSsemimonotone:maximal} and
                \Cref{it:prop:semi:WMVI:max}
            that $\Tpd$
            is maximally
            $
                \DRSRho
                    \coloneqq
                \blkdiag(\DRSrhocom \I_n, \DRSrhomon \I_n)
            $\hyp{}comonotone, where 
            \(
                [\DRSrhomon]_{-}[\DRSrhocom]_{-} < \tfrac{1}{4}
            \).
    \end{proofitemize}
        Moreover, the stepsize $\gamma$ satisfies \eqref{eq:thm:pppaDRS:gamma}, since
    the provided bounds on $\gamma$ from \Cref{tab:gamma:DRS:semi} are obtained by plugging in \eqref{eq:DRSRho:semimonotonicity} into \eqref{eq:thm:pppaDRS:gamma}.
        Based on these observations, the claims follow directly from \Cref{thm:pppaDRS} and \Cref{cor:DRS:lastiter}, concluding the proof.
\end{appendixproof}

  \section{Examples}\label{sec:app:examples}

\begin{appendixdetail}{ex:pppa:toy}
    \begin{proofitemize}
        \item \ref{it:ex:pppa:toy:f1:minty}: 
            Since \(T\) is equal to the linear mapping \(T = \begin{bsmallmatrix}
                b & a \\
                -a & b
            \end{bsmallmatrix}\), it is $\rho$-comonotone if and only if \(\tfrac{1}2(T+T^\top) \preceq \rho T^\top T\) which holds for
                $\rho \leq \nicefrac{b}{a^2 + b^2}$.
        \item \ref{it:ex:pppa:toy:f1:convergence}: 
        To show that the result from \Cref{it:pppa:full} is tight, define
        \ifspringer
            \begin{align*}
                H 
                    {}\coloneqq{}&
                \I_2 + \lambda\left((\I_2 + 
                \gamma
                T)^{-1}-\I_2\right)
                    {}={}
                \I_2 + \lambda\left(
                    \tfrac{1}{(1+
                    \gamma
                    b)^2+
                    \gamma^2
                    a^2}
                    \begin{bmatrix}
                        1+
                        \gamma
                        b & -
                        \gamma
                        a \\
                        \gamma
                        a & 1+
                        \gamma
                        b
                    \end{bmatrix}-\I_2\right),
            \end{align*}
        \else
            \begin{align*}
                H 
                    {}\coloneqq{}&
                    \I_2 + \lambda\left((\I_2 + 
                    \gamma
                    T)^{-1}-\I_2\right)
                    {}={}
                \I_2 + \lambda\left(
                    \begin{bmatrix}
                        1+
                        \gamma
                        b & 
                        \gamma
                        a \\
                        -
                        \gamma
                        a & 1+
                        \gamma
                        b
                    \end{bmatrix}^{-1}-\I_2\right)
                    {}={}
                \I_2 + \lambda\left(
                    \tfrac{1}{(1+
                    \gamma
                    b)^2+
                    \gamma^2
                    a^2}
                    \begin{bmatrix}
                        1+
                        \gamma
                        b & -
                        \gamma
                        a \\
                        \gamma
                        a & 1+
                        \gamma
                        b
                    \end{bmatrix}-\I_2\right),
            \end{align*}
        \fi
        such that the iteration of the proximal point algorithm with fixed relaxation parameter $\lambda$ can be expressed as the linear dynamical system $z^{k+1} = H z^k$.
        This system is globally asymptotically stable if and only if the spectral radius of
                \(H\), given by 
            \(
                \sqrt{1-\tfrac{\lambda (2 \bar \lambda - \lambda)(a^2 + b^2)\gamma^2}{(a^2 + b^2)\gamma^2 + 2 b \gamma + 1}}
            \),
            is strictly less than one, which holds iff \eqref{eq:ex:pppa:toy:lambda} is satisfied.
        \item \ref{it:ex:pppa:toy:f2:minty}: 
        Recall that for a given $\rho \in \R$, the $\rho$-weak Minty solutions of this problem are given by the set $\pazocal{S}^\star \subseteq \zer T$ for which it holds that
        \begin{align*}
            \inner{T(z), z - z^\star} \geq \rho \nrm{T(z)}^2,\quad \forall z \in \R^2, z^\star\in\pazocal{S}^\star.\numberthis\label{eq:example:pppa:WMVI}
        \end{align*}
        The analytical expression for $f$ as defined in \Cref{fig:ex:pppa:toy:f} is given by
        \begin{align*}
            f\Bigl(\nrm{z}\Bigr) \coloneqq 
            \begin{cases}
                \nrm{z}, & \quad \text{if } \nrm{z} \leq \nicefrac25,\\
                \nicefrac45-\nrm{z}, & \quad \text{if }\nrm{z} \in \bigl(\nicefrac25, \nicefrac45\bigr),\\
                0, & \quad \text{if } \nrm{z} \in \bigl[\nicefrac45, 1\bigr],\\
                \nicefrac52\bigl(\nrm{z} - 1\bigr), & \quad \text{if } \nrm{z} \in \bigl(1, \nicefrac75\bigr),\\1, &\quad \text{if } \nrm{z} \geq \nicefrac75.
            \end{cases}\numberthis\label{eq:example:pppa:f}
        \end{align*}
        Therefore, by definition of $T$, it holds that
        $\zer T = \set{z \in \R^2}[\nrm{z} \in \{0\} \cup{[\nicefrac45, 1]}]$ and that \eqref{eq:example:pppa:WMVI} corresponds to
        \ifspringer
            \begin{align*}
                f\Bigl(\nrm{z}\Bigr)\Bigl(b\nrm{z}^2 - b\bigl(z_1 z_1^\star + z_2 z_2^\star\bigr) + a\bigl(z_1z_2^\star - z_1^\star z_2\bigr)\Bigr) 
                    {}\geq{}
                \rho \bigl(a^2 + b^2\bigr) f\Bigl(\nrm{z}\Bigr)^2 \nrm{z}^2,\numberthis\label{eq:example:pppa:WMVI:T}
            \end{align*}
            for all $z \in \R^2$ and $z^\star\in\pazocal{S}^\star$.
        \else
            \begin{align*}
                f\Bigl(\nrm{z}\Bigr)\Bigl(b\nrm{z}^2 - b\bigl(z_1 z_1^\star + z_2 z_2^\star\bigr) + a\bigl(z_1z_2^\star - z_1^\star z_2\bigr)\Bigr) 
                    {}\geq{}
                \rho \bigl(a^2 + b^2\bigr) f\Bigl(\nrm{z}\Bigr)^2 \nrm{z}^2,\qquad \forall z \in \R^2, z^\star\in\pazocal{S}^\star.\numberthis\label{eq:example:pppa:WMVI:T}
            \end{align*}
        \fi
        This obviously holds for all $z \in \zer T$. For $z \notin \zer T$, we distinguish between the following two cases.
        \begin{proofitemize}
            \item $\nrm{z^\star} = 0$: Then, \eqref{eq:example:pppa:WMVI:T} holds iff for all $z \notin \zer T$ it holds that
            \(
                \rho
                    {}\leq{}
                \tfrac{b}{(a^2 + b^2) f\left(\nrm{z}\right)}.
            \)
            Since $f\left(\nrm{z}\right) \in (0,1]$, this condition is satisfied for
            \(
                \rho = \nicefrac{b}{a^2 + b^2}.
            \)
            \item $\nrm{z^\star} \in [\nicefrac45, 1]$: Let $z = \epsilon z^\star$ where $\epsilon \in (0,\nicefrac25)$. Then,  
            $f(\nrm{z}) = \epsilon\nrm{z^\star}$ and \eqref{eq:example:pppa:WMVI:T} holds iff
            \(
                \rho 
                    {}\leq{}
                \tfrac{
                    b(1 - \epsilon^{-1})
                }{                    
                    \epsilon(a^2 + b^2) \nrm{z^\star}
                }
            \)
            Taking the limit for $\epsilon$ to zero from above, we get that 
            $
                \lim_{\epsilon \rightarrow 0+} 
                \tfrac{b(1 - \epsilon^{-1})}{\epsilon(a^2 + b^2) \nrm{z^\star}}
                    {}={}
                - \infty
            $
            and thus \eqref{eq:example:pppa:WMVI:T} cannot be satisfied for any $\rho \in \R$.\qedhere
        \end{proofitemize}
    \end{proofitemize}
\end{appendixdetail}
    \ifspringer\else

\begin{appendixdetail}{ex:VonNeumann}
    We proceed by establishing three intermediate facts for the set $C \coloneqq \Delta_2 \times \Delta_2$ and the operators $F$ and $N_C$. 
    First, observe that by definition of the simplex
    \begin{align*}
        z \in C
        \iff 
        z = (s, 1-s, t, 1-t) 
        \;\text{for some}\; s, t \in [0,1].
        \numberthis\label{eq:VonNeumann:parametrization}
    \end{align*}
    Second, recall that
    \[
        F(z)
            {}={}
        \left(%
            F_1(z), F_2(z), F_3(z), F_4(z)
        \right)
            {}={}
        \left(
            \nabla_x f(x,y), - \nabla_y f(x,y)
        \right)
            {}={}
        \left(
            \tfrac{\inner{x, Sy} Ry - \inner{x, Ry} Sy}{\inner{x, Sy}^2},
            \tfrac{\inner{x, Ry} S^\top x - \inner{x, Sy} R^\top x}{\inner{x, Sy}^2}
        \right).
    \]
    Therefore, using the parametrization from \eqref{eq:VonNeumann:parametrization} it holds for any $z \in C$ that
    \begin{align*}
        \inner{ 
                F(z),
                z - z^\star
            }
            {}={}&
        \tfrac{                        
                \epsilon(s + t)
                -
                st(\epsilon + s)
            }{
                \nicefrac12(2 - s)^2
            },
            \quad
        F_1(z) - F_2(z)
            {}={}
        \tfrac{
                \epsilon - t(2 + \nicefrac{\epsilon}{2})
            }{
                \nicefrac12(2 - s)^2
            }
        \quad\text{and}\quad
        F_3(z) - F_4(z)
            {}={}
        \tfrac{\nicefrac{\epsilon}{2} + s}{\nicefrac12(2 - s)}.
        \numberthis\label{eq:VonNeumann:F}
    \end{align*}
    Finally, denoting the vector of all $1$s
    of dimension $n$ by $\mathbf{1}_n$, it holds owing to \cite[Ex. 5.2.6(c)]{Hiriart1996Convex} that
    \begin{align*}
        N_C(z)
            {}={}
        \begin{cases}
            \set{
            \alpha \otimes 
            \mathbf{1}_2
            -
            \beta
        \,\middle|\, \alpha \in \R^2,
        \beta \in \R_+^4,
        \inner{\beta, z} = 0}, &\qquad \text{if } z \in C,\\
        \emptyset, &\qquad \text{otherwise}.
        \end{cases}
        \numberthis\label{eq:VonNeumann:NC}
    \end{align*}
    We now proceed with the main argument.
    In \cite[Prop. 2(1)]{daskalakis2020independent}, it was shown that
    $
        z^\star =         
        (0,1,0,1)
    $ is the unique Nash equilibrium of \eqref{eq:ex:neumann} (noting that 
    $
        f(x, y) - 2 =\nicefrac{
            \inner{x,         
            \begin{bsmallmatrix}
                -1 & \nicefrac{\epsilon}{2}\\
                - \nicefrac{\epsilon}{2} & 0
            \end{bsmallmatrix}y
        }}{{\inner{x, Sy}}}
    $).
    This Nash equilibrium is a $\DRSrho_\epsilon$\hyp{}weak Minty solution of $T$ provided that
    \(
        \DRSrho_\epsilon
            {}={}
        \inf_{z \in C} \varrho_\epsilon(z)
    \)
    where
    \(
        \varrho_\epsilon(z)
            {}={}
        \inf_{v \in N_C(z)}
        \nicefrac{
            \inner*{F(z) + v,z-z^\star}
        }
        {
            \nrm{F(z) + v}^2
        }
    \).
    From now onwards, we tacitly assume that $z \in C$. Then, owing to \eqref{eq:VonNeumann:NC}, it holds that
    \begin{align*}
        \varrho_\epsilon(z)
            {}={}
        \inf\set{
        \tfrac{%
            \inner*{%
                F(z),
                z - z^\star
            } 
                +
            \inner*{                   
                \alpha \otimes 
                \mathbf{1}_2
                    -
                \beta,
                z - z^\star
            } 
        }
        {\nrm{
            F(z)
                +                         
            \alpha \otimes 
            \mathbf{1}_2
            -
            \beta
        }^2
        }%
        }
        [\alpha \in \R^2,\, \beta \in \R_+^4, 
        \;\inner{
            \beta,
            z
        } 
        = 0
        ].
    \end{align*}
    Owing to the parametrization from \eqref{eq:VonNeumann:parametrization} and the fact that    
    $
        z^\star =
            (0,1,0,1)
    $,
    it holds that
    \(
        \inner*{                   
                \alpha \otimes 
                \mathbf{1}_2,
                z - z^\star
            }
            {}={}
        0
    \) 
    and
    $\inner{\beta, z^\star} = \beta_2 + \beta_4$, so that
    \begin{align*}
        \varrho_\epsilon(z)
            {}={}& 
        \inf\set{
            \tfrac{
                \inner*{ 
                    F(z),
                    z - z^\star
                } 
                    +
                \beta_2
                    +
                \beta_4
            }
            {\nrm{
                F(z)
                    +                         
                \alpha \otimes 
                \mathbf{1}_2
                -
                \beta
            }^2
            }
        }
        [\alpha \in \R^2,\, \beta \in \R_+^4, 
        \,\inner{
            \beta,
            z
        } 
        = 0
        ].
    \end{align*}
    Note that by construction the sign of $\varrho_\epsilon(z)$ is equal to the sign of 
    $            
        \inner*{ 
            F(z),
            z - z^\star
        }
    $:
    if
    \(\inner*{ 
        F(z),
        z - z^\star
    }\)
    is positive, this is obvious since \(\beta_i \geq 0\); if it is negative, then this can be seen by taking the candidate solution $\beta_2 = \beta_4 = 0$.
    Therefore, it follows from \eqref{eq:VonNeumann:F} that $\varrho_\epsilon(z) < 0$
    if and only if
    \(
        \epsilon(s + t)
            -
        st(\epsilon + s)
            <
        0
    \), i.e., if and only if $z$ lies in the nonempty set
    \[
        D
            \coloneqq
        \set{(s, 1-s, t, 1-t)}
        [{s \in (\hspace{-2pt}\sqrt{\epsilon}, 1],\, t \in (\bar{t}(s), 1]}]
            \quad\text{where}\quad
        \bar{t}(s)
            {}\coloneqq{}
        \tfrac{\epsilon s}{s(\epsilon + s)-\epsilon}.
    \]
    Therefore, 
    \(
        \DRSrho_\epsilon
            {}={}
        \inf_{z \in C} \varrho_\epsilon(z)
            {}={}
        \inf_{z \in D} \varrho_\epsilon(z).
    \)
    Consider $z \in D \subset C$. Then, by construction $\varrho_\epsilon(z)$ is negative, which implies that \(\varrho_\epsilon(z)\) is minimized for \(\alpha\) that minimizes the denominator, given by 
    \[
        \argmin_{\alpha \in \R^2} \;
        \nrm{
            F(z)
                +                         
            \alpha \otimes 
            \mathbf{1}_2
            -
            \beta
        }^2
            {}={}
        \left(
            \tfrac12(F_1(z) - F_2(z) - \beta_1 + \beta_2), 
            \tfrac12(F_3(z) - F_4(z) - \beta_3 + \beta_4)\right),
    \]
    this implies that
    \begin{align*}
        \varrho_\epsilon(z)
            {}={}& 
        \inf\set{\tfrac{
            2\inner*{ 
                F(z),
                z - z^\star
            } 
                +
            2\beta_2
                +
            2\beta_4
        }
        {
            \bigl(
                F_1(z) - \beta_1 - F_2(z) + \beta_2
            \bigr)^2
                +
            \bigl(
                F_3(z) - \beta_3 - F_4(z) + \beta_4
            \bigr)^2
        }
        }[
            \beta \in \R_+^4,
            \,
            \inner{
                \beta,
                z
                } 
            = 0
        ].
    \end{align*}
    Plugging in the parametrization for $z$ from \eqref{eq:VonNeumann:parametrization}, using \eqref{eq:VonNeumann:F} and introducing $\bar{\beta} \coloneqq \nicefrac12(2 - s)^2\beta$, it follows that
    \begin{align*}
        \varrho_{\epsilon}(s,t)
            {}={}
        \inf\set{
        \tfrac{
            4(2 - s)^2
            \bigl(
                \epsilon(s + t)
                    -
                st(\epsilon + s)
                    +
                \bar{\beta}_2
                    +
                \bar{\beta}_4
            \bigr)
        }
        {
            \bigl(
                2\epsilon - t(4 + \epsilon)
                -
                2\bar{\beta}_1
                + 
                2\bar{\beta}_2
            \bigr)^2
                +
            \bigl(
                (\epsilon + 2s)(2 - s)
                -
                2\bar{\beta}_3
                +
                2\bar{\beta}_4
            \bigr)^2
        }
        }[
            \bar{\beta} \in \R_+^4,\,
            \inner{
            \bar{\beta},
            (s, 1-s, t, 1-t)
            } 
        = 0
        ],       
    \end{align*}
    \begin{figure}
        \centering
        \includetikz{Examples/VonNeumann/test}
        \caption{
            Comparison between the optimal weak Minty constant $\DRSrho_\epsilon$ and its lower bound $\tfrac{\epsilon - 1}{4}$ for \cref{ex:VonNeumann}.
            }
        \label{fig:VonNeumann:rho}
    \end{figure}%
    We claim that $\bar \beta = 0$ at the infimum for any $s \in (\hspace{-2pt}\sqrt{\epsilon}, 1]$ and $t \in (\bar{t}(s), 1]$.
    For $\bar \beta_1$ and $\bar \beta_3$, this follows immediately from the equality constraint and the fact that $s$ and $t$ are strictly positive, so that
    \begin{align*}
        \varrho_{\epsilon}(s,t)
            {}={}
        \inf\set{
        \tfrac{
            4(2 - s)^2
            \bigl(
                \epsilon(s + t)
                    -
                st(\epsilon + s)
                    +
                \bar{\beta}_2
                    +
                \bar{\beta}_4
            \bigr)
        }
        {
            \bigl(
                2\epsilon - t(4 + \epsilon)
                + 
                2\bar{\beta}_2
            \bigr)^2
                +
            \bigl(
                (\epsilon + 2s)(2 - s)
                +
                2\bar{\beta}_4
            \bigr)^2
        }
        }[
            (\bar{\beta}_2, \bar{\beta}_4) \in \R_+^2,\,
            \inner{
                (\bar{\beta}_2, \bar{\beta}_4),
                (1-s, 1-t)
                } 
            = 0
        ].
    \end{align*}
    If $s \in (\hspace{-2pt}\sqrt{\epsilon}, 1)$ and $t \in (\bar{t}(s), 1)$, the claim for $\bar \beta_2$ and $\bar \beta_4$ also follows from the equality constraint. Otherwise, consider the following cases.
    \begin{proofitemize}
        \item If $s = 1$ and $t \in (\bar{t}(1), 1) = (\epsilon, 1)$, then $\bar{\beta}_4 = 0$ due to the equality constraint and
        \begin{align*}
            \varrho_{\epsilon}(1,t)
                {}={}
            \inf_{
                \bar{\beta}_2 \in \R_+
            } \; 
            \tfrac{
                4
                (
                    \epsilon
                        -
                    t
                        +
                    \bar{\beta}_2
                )
            }
            {
                (
                    2\epsilon - t(4 + \epsilon)
                    + 
                    2\bar{\beta}_2
                )^2
                    +
                (
                    \epsilon + 2
                )^2
            }
                {}={}
            \inf_{
                \bar{\beta}_2 \in [0, t - \epsilon)
            } \; 
            \tfrac{
                4
                (
                    \epsilon
                        -
                    t
                        +
                    \bar{\beta}_2
                )
            }
            {
                (
                    2\epsilon - t(4 + \epsilon)
                    + 
                    2\bar{\beta}_2
                )^2
                    +
                (
                    \epsilon + 2
                )^2
            }.
        \end{align*}
        Since this function is monotonically increasing for $\bar{\beta}_2 \in [0, t - \epsilon)$, the infimum is attained for $\bar{\beta}_2 = 0$.
        \item If $s \in (\hspace{-2pt}\sqrt{\epsilon}, 1)$ and $t = 1$, then $\bar{\beta}_2 = 0$ due to the equality constraint and
        \begin{align*}
            \varrho_{\epsilon}(s,1)
                {}={}
            \inf_{
                \bar{\beta}_4 \in \R_+
            } \; 
            \tfrac{
                4(2 - s)^2
                (
                    \epsilon
                        -
                    s^2
                        +
                    \bar{\beta}_4
                )
            }
            {
                (
                    \epsilon - 4
                )^2
                    +
                \bigl(
                    (\epsilon + 2s)(2 - s)
                    +
                    2\bar{\beta}_4
                \bigr)^2
            }
                {}={}
            \inf_{
                \bar{\beta}_4 \in [0, s^2 - \epsilon)
            } \; 
            \tfrac{
                4(2 - s)^2
                (
                    \epsilon
                        -
                    s^2
                        +
                    \bar{\beta}_4
                )
            }
            {
                (
                    \epsilon - 4
                )^2
                    +
                \bigl(
                    (\epsilon + 2s)(2 - s)
                    +
                    2\bar{\beta}_4
                \bigr)^2
            }.
        \end{align*}
        Since this function is monotonically increasing for $\bar{\beta}_4 \in [0, s^2 - \epsilon)$, the infimum is attained for $\bar{\beta}_4 = 0$.
        \item If $s = 1$ and $t = 1$, then
        \begin{align*}
            \varrho_{\epsilon}(1,1)
                {}={}
            \inf_{
                (\bar{\beta}_2, \bar{\beta}_4) \in \R_+^2
            } \; 
            \tfrac{
                4
                (
                    \epsilon
                        -
                    1
                        +
                    \bar{\beta}_2
                        +
                    \bar{\beta}_4
                )
            }
            {
                (
                    \epsilon - 4
                    + 
                    2\bar{\beta}_2
                )^2
                    +
                (
                    \epsilon + 2
                    +
                    2\bar{\beta}_4
                )^2
            },
        \end{align*}
        of which the global minimum is attained for $\bar{\beta}_2 = \bar{\beta}_4 = 0$.
    \end{proofitemize}
    Having established that $\bar \beta = 0$ at the infimum, it follows that
    \begin{align*}
        \DRSrho_\epsilon
            {}={}& 
        \inf\set{
            \,\varrho_\epsilon(s, t)
        }[
            {
            s \in (\sqrt{\epsilon}, 1],\,
            t \in \left(
                    \tfrac{\epsilon s}{s(\epsilon + s)-\epsilon}, 1
                \right]
            }
        ]
            \quad\text{where}\quad
        \varrho_\epsilon(s, t)
            =
        \tfrac{
            4(2 - s)^2
            \bigl(
                \epsilon(s + t) - st(\epsilon + s)
            \bigr)
        }
        {
            \bigl(
                2\epsilon - t(4 + \epsilon)
            \bigr)^2
                +
            (\epsilon + 2s)^2(2 - s)^2
        }.
        \numberthis\label{eq:prop:ex:neumann:localMinty}
    \end{align*}
    Finally, the lower bound $\tfrac{\epsilon - 1}{4} \leq \DRSrho_\epsilon$ is obtained by solving this parametric optimization problem numerically (see also \Cref{fig:VonNeumann:rho}).
\end{appendixdetail}

\begin{appendixdetail}{example:ex5}[ (saddle point problem)]
    \begin{proofitemize}
        \item \ref{it:example5:1}: By defining
        \(\overbar H \coloneqq (\I + \gamma A)^{-1}\),
        \(\overhat H \coloneqq (\I + \gamma B)^{-1} (2 \overbar H - \I)\) and
        \(H \coloneqq \I + \lambda (\overhat H - \overbar H)\), \ref{eq:DRS} can be expressed as the linear dynamical system $s^{k+1} = H s^k$. This system is globally asymptotically stable if and only if the spectral radius of \(H\) is strictly less than one, which through elementary algebra can be shown to hold iff \eqref{eq:example5:lambda} is satisfied.
        \item \ref{it:example5:2}:
        The primal-dual operator and its inverse are given by 
        $$
            \Tpd 
                {}={}
            \begin{bmatrix}
                A & \I_2 \\ -\I_2 & B^{-1}
            \end{bmatrix}
            \quad\text{and}\quad
            \Tpd^{-1}
                {}={}
            \begin{bmatrix}
                (A + B)^{-1} & - (A + B)^{-1}B\\
                B (A + B)^{-1} & B - B (A + B)^{-1}B
            \end{bmatrix}
        $$
        where we used the Schur complement lemma.
        Consequently, when in \eqref{def:WMVI} the vector $v$ is restricted to $\range{\M}$ we may equivalently state \Cref{ass:SWMVIDRS} as 
        \begin{equation}
            \tp z \left(\tfrac{\Tpd + \tp \Tpd}{2} - \tp \Tpd \DRSRho \Tpd\right) z \geq 0,\qquad \text{for all $z \in \R^n : z \in \Tpd^{-1}\range{\M}$},\label{eq:example:linear:WMVI}
       \end{equation}
        where $\DRSRho$ is defined as in~\eqref{eq:WMSDRS}. Using that
        \(U = \begin{bsmallmatrix} \I_2 \\ -\gamma\I_2\end{bsmallmatrix}\)
        is an orthogonal basis for $\range{\M}$, it follows that \eqref{eq:example:linear:WMVI} is satisfied iff
        \[
            \bigl(\Tpd^{-1} U\bigr)^\top\left(\tfrac{\Tpd + \tp \Tpd}{2} - \tp \Tpd \DRSRho \Tpd\right) \Tpd^{-1} U \succeq 0
                \quad\Longleftrightarrow\quad
            \tfrac{1}{\gamma} \DRSrhocom + \gamma \DRSrhomon \leq \tfrac{b(a^2 \gamma + 1)}{\gamma(a^2 + b^2)}.
        \]
        Consequently, for the upper bound on $\lambda$ implied by \Cref{thm:pppaDRS}, denoted by $\lambda_{\rm{max}}$, it holds that
        \[
            \lambda_{\rm{max}} \coloneqq 2\left(1+\tfrac{1}{\gamma} \DRSrhocom + \gamma \DRSrhomon \right)
                {}\leq{}
            2 + 2\tfrac{b(a^2 \gamma + 1)}{\gamma(a^2 + b^2)}.
        \]     
        As this bound matches the upper bound on \(\lambda\) from \eqref{eq:example5:lambda}, the proof is completed.
        \item \ref{it:example5:3}: For example, consider \(a = 2, b = -1\), for which the trace of \(\Tp\),\(\Td\) and \(\Tpd\) are all equal to -2.
        \item  
        Finally, we will discuss the range of relaxation parameters $\lambda$ covered by \Cref{cor:DRS:semi}. To this end, the semimonotonicity parameters $\monA$, $\comA$, $\monB$, $\comB$ need to be selected such that the relaxation parameter is maximized, i.e.
        $$
        \begin{aligned}
            \lambda_{\rm{max}}^\mathrm{semi} \;\coloneqq \max_{\substack{\monA, \comA, \monB, \comB}} \quad& 2\left(1 + \tfrac{1}{\gamma}(\comA \Box \comB) + \gamma(\monA \Box \monB)\right)\\
            \stt \quad& \monA + a^2\comA \leq 0 \quad \text{and} \quad \monB + b^2 \comB - b \leq 0,
        \end{aligned}
        $$
        where the constraints are obtained using \Cref{it:prop:semi:lin:1} by enforcing $(\monA, \comA)$\hyp{} and $(\monB,\comB)$\hyp{}semimonotonicity of operators $A$ and $B$, respectively.
        Note that this optimization problem can be encoded as a convex SDP. Solving this problem numerically
        for $a = 10$, $b=-1$, the result visualized in \Cref{fig:ex5:semi vs spectral} is obtained. Hence, when characterizing the nonmonotonicity of the problem through the semimonotonicity properties of $A$ and $B$ separately rather than through the (oblique) weak Minty condition for the combined primal-dual operator, some looseness is inevitably introduced.
        \qedhere
        \begin{figure}
            \centering
            \includetikz{Examples/Saddle/saddle-thm-semi}
            \caption{Comparison between the upper bounds on the relaxation parameter $\lambda$ in function of the stepsize $\gamma$ for \Cref{example:ex5} with $a=10$, $b=-1$. $\bar \lambda$ is the bound given by \eqref{eq:example5:lambda} and $\lambda_\mathrm{max}^\mathrm{semi}$ is the bound implied by \Cref{cor:DRS:semi}.}
            \label{fig:ex5:semi vs spectral}
        \end{figure}
    \end{proofitemize}
\end{appendixdetail}

\begin{appendixdetail}{ex:nonsmooth:minimization}[ (nonsmooth optimization)]
    The corresponding limiting subdifferentials are given by
    \begin{equation}
        \begin{aligned}
            \partial f_A(x) =
            \begin{cases}
                -x-3, & \quad\text{if } x < -1,\\
                \{-2, -9\}, & \quad \text{if } x = -1,\\
                6x-3, & \quad\text{if } x \in (-1, 1), \\
                \{-4, 3\}, & \quad \text{if } x = 1,\\
                -x-3, & \quad\text{if } x > 1,
            \end{cases}
            \quad \text{and} \quad
            \partial f_B(x) =
            \begin{cases}
                4x - 1, & \quad\text{if } x < -1,\\
                \left[-5,1\right] & \quad\text{if } x=-1,\\
                2x + 3, & \quad\text{if } x \in (-1, 2),\\
                \left[7, 19\right], & \quad\text{if } x = 2,\\
                2x + 15, & \quad\text{if } x > 2.
            \end{cases}
        \end{aligned}
    \end{equation}
    \vspace{-\baselineskip}
    \begin{proofitemize}
        \item \Cref{ass:DRS:1}: Follows by definition of the limiting subdifferential \cite[Def. 8.3(b)]{Rockafellar1970Convex}.
        \item \Cref{ass:DRS:2}: It is easy to verify that for $\gamma \neq 1$ it holds that $\range{1 + \gamma \partial f_A} = \range{1 + \gamma \partial f_B} = \R$. Hence, the resolvents $J_{\gamma \partial f_A}$ and $J_{\gamma \partial f_B}$ have full domain.
        \item \Cref{ass:SWMVIDRSsemimonotone}: Let $(x^\star, y^\star) = (0, 3)$. Then, it is of immediate verification that the subdifferentials \(\partial f_A\) and \(\partial f_B\) are semimonotone at respectively $(x^\star, -y^\star)$ and $(x^\star, y^\star)$ with \(\monA = -1.2\), \(\comA = 0.2\), \(\monB = 1.6\) and \(\comB = 0.1\).
        \item Follows by applying \Cref{cor:DRS:semi} and using that $\monA \Box \monB = -\nicefrac{24}{5}$ and $\comA \Box \comB = \nicefrac{2}{30}$. 
        \qedhere
    \end{proofitemize}
\end{appendixdetail}

\begin{appendixdetail}{ex:nonsmooth:stationary}[ (stationary point)]
    Define $g_1(x) \coloneqq (\nicefrac{1-\sqrt{1-4\monA\comA}}{2\comA})x - 1$, $g_2(x) \coloneqq (\nicefrac{1+\sqrt{1-4\monA\comA}}{2\comA})x - 1$, $h_1(x) \coloneqq 2x+1$ and $h_2(x) \coloneqq \tfrac{1}{2}x+1$. Then, the limiting subdifferentials of $f_A$ and $f_B$ are given by 
    \begin{align*}
        \partial f_A(x)
            {}={}&
        \begin{cases}
            \Bigl\{g_1(x), g_2(x)\Bigr\}, & \quad\text{if } x \in \{-3, 3\}, \\
            g_1(x), & \quad\text{if } x \in (-3, 3), \\
            g_2(x), & \quad\text{otherwise},
        \end{cases}
        \quad \text{and} \quad
        \partial f_B(x)
            {}={}
        \begin{cases}
            \Bigl\{h_1(x),h_2(x)\Bigr\}, & \quad\text{if } x \in \{-1, 1\}, \\
            h_1(x), & \quad\text{if } x \in (-1, 1), \\
            h_2(x), & \quad\text{otherwise}.
        \end{cases}
    \end{align*}
    \vspace{-\baselineskip}
    \begin{proofitemize}
        \item \Cref{ass:DRS:1}: Follows by definition of the limiting subdifferential \cite[Def. 8.3(b)]{Rockafellar1970Convex}.
        \item \Cref{ass:DRS:2}: 
        The resolvents $J_{\gamma \partial f_A}$ and $J_{\gamma \partial f_B}$ have full domain iff $\range{1 + \gamma \partial f_A} = \range{1 + \gamma \partial f_B} = \R$.
        It is easy to verify that $\range{1 + \gamma \partial f_B} = \R$ for all $\gamma \in \R_{+}$. On the other hand, we have that
        \[
            (1 + \gamma \partial f_A)(x)
                {}={}
            \begin{cases}
                x - \gamma + \Bigl\{\tfrac{1-\sqrt{1-4\monA\comA}}{2\comA}\gamma x, \tfrac{1+\sqrt{1-4\monA\comA}}{2\comA}\gamma x\Bigr\}, & \quad\text{if } x \in \{-3, 3\}, \\
                x - \gamma + \tfrac{1-\sqrt{1-4\monA\comA}}{2\comA}\gamma x, & \quad\text{if } x \in (-3, 3), \\
                x - \gamma + \tfrac{1+\sqrt{1-4\monA\comA}}{2\comA}\gamma x, & \quad\text{otherwise}
            \end{cases}
        \]
        from which follows that $\range{1 + \gamma \partial f_A} = \R$ for all $\gamma \in \R_{+}$ iff $\gamma$ is selected according to the provided range.
        \item \Cref{ass:SWMVIDRSsemimonotone}: Let $(x^\star, y^\star) = (0, -1)$. Then, it is of immediate verification that the subdifferentials \(\partial f_A\) and \(\partial f_B\) are semimonotone at respectively $(x^\star, -y^\star)$ and $(x^\star, y^\star)$ for the given parameters \(\monA, \comA, \monB, \comB\). Furthermore, these parameters satisfy \Cref{ass:SWMVIDRSsemimonotone:params}.
        \item Follows by applying \Cref{cor:DRS:semi} and using that $\monA \Box \monB = -\nicefrac{6}{5}$ and $\comA \Box \comB = -\nicefrac{4}{30}$. 
        \qedhere
    \end{proofitemize}
\end{appendixdetail}
    \fi

  \end{appendix}

  \ifspringer
    \bibliographystyle{jabbrv_abbrv}%
    \bibliography{TeX/Minty-DRS}%
  \else
    \bibliographystyle{arxivplain}%
    \bibliography{TeX/Minty-DRS.bib}%
  \fi

\end{document}